\documentclass[a4paper,10pt]{amsart}
\usepackage[foot]{amsaddr}
\usepackage{geometry}
\theoremstyle{remark}

\usepackage{graphicx}
\usepackage{lineno}
\usepackage{bm}
\usepackage{amsfonts} 
\usepackage{amsmath}
\usepackage{amssymb}
\usepackage{framed} 
\usepackage{multicol} 
\usepackage{tabularx}
\usepackage{tabulary}
\usepackage{siunitx}
\usepackage{caption}
\usepackage{float}
\usepackage{url} 
\usepackage{multicol}
\usepackage{breakurl} 
\usepackage[breaklinks]{hyperref} 

\usepackage{subcaption}
\usepackage{xpatch}
\usepackage{nomencl} 
\makenomenclature
\makenomenclature 
\RequirePackage{ifthen}
\let\oldref\ref
\renewcommand{\ref}[1]{(\oldref{#1})}

 \graphicspath{
    {figures_ok/}
    {cyl/}
    {cavity/}
}

\usepackage{amssymb}
\providecommand{\norm}[1]{\lVert#1\rVert}

\DeclareMathAlphabet\mathbfcal{OMS}{cmsy}{b}{n}

\xpatchcmd{\thenomenclature}{%
  \section*{\nomname}
}{
}{\typeout{Success}}{\typeout{Failure}}

\usepackage{ifthen}
\renewcommand{\nomgroup}[1]{%
  \ifthenelse{\equal{#1}{A}}{\item[\textbf{Abbreviations}]}{%
    \ifthenelse{\equal{#1}{G}}{\item[\textbf{Symbols}]}{%
      \ifthenelse{\equal{#1}{C}}{\item[\textbf{Abbreviations}]}{%
        \ifthenelse{\equal{#1}{S}}{\item[\textbf{Subscripts}]}{%
          \ifthenelse{\equal{#1}{Z}}{\item[\textbf{Mathematical Symbols}]}{}
        }
      }
    }
  }
}

\ifpdf
\usepackage{pdfsync}
\fi

\usepackage{mathabx}

\begin{document}
\nomenclature{$\text{ROM}$}{Reduced Order Model}
\nomenclature{$\text{POD}$}{Proper Orthogonal Decomposition}
\nomenclature{$\text{RB}$}{Reduced Basis}
\nomenclature{$\text{FEM}$}{Finite Element Methods}
\nomenclature{$\text{FVM}$}{Finite Volume Methods}
\nomenclature{$\text{Mixed-ROM}$}{The mixed projection/data-driven based reduced order model developed in this work}
\nomenclature{$\text{RANS}$}{Reynolds Average Navier-Stokes}
\nomenclature{$\text{RBF}$}{Radial Basis Functions}
\nomenclature{$\text{EVM}$}{Eddy Viscosity Models}

\nomenclature[G]{$\bm{u}$}{velocity field}
\nomenclature[G]{${p}$}{pressure field}
\nomenclature[G]{${\nu}$}{dimensionless kinematic viscosity}
\nomenclature[G]{${N_u^h}$}{number of unknowns for velocity at full-order level}
\nomenclature[G]{${N_{\nu_t}^h}$}{number of unknowns for eddy viscosity at full-order level}
\nomenclature[G]{${N_p^h}$}{number of unknowns for pressure at full-order level}
\nomenclature[G]{$N_s$}{total number of snapshots}
\nomenclature[G]{$T$}{final time}
\nomenclature[G]{$\mathcal{P}$}{parameter space of dimension $q$}
\nomenclature[G]{$q$}{the dimension of the parameter space $\mathcal{P}$}
\nomenclature[G]{$\mathcal{P}_M$}{training set in the parameter space $\mathcal{P}$ with cardinality of $M$}
\nomenclature[G]{$\nu_t$}{eddy viscosity field}
\nomenclature[G]{${\Omega}$}{bounded domain}
\nomenclature[G]{${\Gamma}$}{boundary of $\Omega$}
\nomenclature[G]{${\Gamma_{In}}$}{the inlet boundary part of ${\Gamma}$}
\nomenclature[G]{${\Gamma_0}$}{the physical walls boundary part of ${\Gamma}$}
\nomenclature[G]{${\Gamma_{Out}}$}{the outlet boundary part of ${\Gamma}$}
\nomenclature[G]{$\otimes$}{Tensor product}

\nomenclature[G]{$\bm{n}$}{outward normal vector}
\nomenclature[G]{$\bm{\phi}_i$}{$i$-th POD basis function for velocity}
\nomenclature[G]{$\eta_i$}{$i$-th POD basis function for eddy viscosity}
\nomenclature[G]{$\chi_i$}{$i$-th POD basis function for pressure}
\nomenclature[G]{$\bm{M}$}{ROM mass matrix} 
\nomenclature[G]{$\bm{C}$}{ROM convection tensor}
\nomenclature[G]{$\bm{C_{T1}}$}{ROM turbulent tensor}
\nomenclature[G]{$\bm{C_{T2}}$}{ROM turbulent tensor}
\nomenclature[G]{$\bm{B_T}$}{ROM diffusion turbulent matrix}
\nomenclature[G]{$\bm{H}$}{ROM pressure gradient matrix}
\nomenclature[G]{$\bm{P}$}{ROM divergence matrix}

\nomenclature[G]{$\bm{B}$}{ROM diffusion matrix}
\nomenclature[G]{${M}$}{number of parameter samples in the training set $\mathcal{P}_M$}
\nomenclature[G]{$\bm{{\mathcal{S}_u}}$}{snapshots matrix for the velocity field}
\nomenclature[G]{$\bm{{\mathcal{S}_p}}$}{snapshots matrix for the pressure field}
\nomenclature[G]{$\bm{\mathcal{\tilde{U}}}$}{homogenized velocity snapshots matrix}
\nomenclature[G]{$\bm{\mathcal{\tilde{P}}}$}{homogenized pressure snapshots matrix}
\nomenclature[G]{$\bm{{\mathcal{S}_{\nu_t}}}$}{snapshots matrix for the eddy viscosity field}

\nomenclature[G]{$\bigtimes$}{Cartesian product}
\nomenclature[G]{$\bm{\nabla}\cdot$}{divergence operator}
\nomenclature[G]{$\bm{\nabla}$}{gradient operator}
\nomenclature[G]{$\tau$}{The penalization factor in the penalty boundary treatment method}
\nomenclature[G]{${\Gamma_D}_j$}{The $j$-th part of the boundary where Dirichlet conditions are imposed}
\nomenclature[G]{$N_{BC}$}{The number of scalar non-zero boundary conditions needed to be set at reduced order level}
\nomenclature[G]{$\bm{\phi_{L}}$}{The matrix of the lifting functions $\bm{{\phi_L}_{i,j}}$}
\nomenclature[G]{$\bm{{\phi_L}_{i,j}}$}{The velocity lifting function that has unitary value in its $i$-th component at ${\Gamma_D}_j$}

\nomenclature[G]{$\bm{a}$}{reduced vector of unknowns for velocity}
\nomenclature[G]{$\bm{b}$}{reduced vector of unknowns for pressure}
\nomenclature[G]{$\bm{g}$}{reduced vector of unknowns for eddy viscosity}
\nomenclature[G]{${(\cdot , \cdot)}_{L^2(\Omega)}$}{inner product in $L^2(\Omega)$}
\nomenclature[G]{$\overline{\bm{u}}$}{mean velocity field in RANS equations}
\nomenclature[G]{$\bm{u}^\prime$}{fluctuating velocity field in RANS equations}
\nomenclature[G]{$k$}{turbulence kinetic energy}
\nomenclature[G]{$\omega$}{specific turbulent dissipation rate}
\nomenclature[G]{$\epsilon$}{turbulent dissipation}

\nomenclature[G]{$N_u$}{number of modes used in the online phase for velocity}
\nomenclature[G]{$N_S$}{number of modes used in the online phase for the supremizer}
\nomenclature[G]{$L$}{The component that is parallel to the lift direction of the force which acts on a certain surface}

\nomenclature[G]{$N_p$}{number of modes used in the online phase for pressure}
\nomenclature[G]{$N_{\nu_t}$}{number of modes used in the online phase for eddy viscosity}
\nomenclature[G]{$\bm{C^u}$}{correlation matrix of the velocity field snapshot matrix}
\nomenclature[G]{$\bm{V^u}$}{eigenvectors matrix of the correlation matrix of the velocity field snapshot matrix}
\nomenclature[G]{$\bm{\lambda^u}$}{eigenvalues matrix of the correlation matrix of the velocity field snapshot matrix}
\nomenclature[G]{$X_{\bm{\mu},t}$}{the combined set of samples and time instants at which snapshots are taken}
\nomenclature[G]{$\bm{x}_{\bm{\mu},t}^i$}{$i$-th member of $X_{\bm{\mu},t}$}

\nomenclature[G]{$\bm{\mu}^*$}{The sample parameter introduced to the ROM in the online stage}
\nomenclature[G]{$\bm{z}$}{The generic parameter-time vector which lives in $\mathbb{R}^{q+1}$}
\nomenclature[G]{$\bm{z}^*$}{The sample parameter combined with the time instant at which Mixed-ROM solution is desired in the online stage}
\nomenclature[G]{$\zeta_{i,j}$}{The RBF functions used in interpolating the $i$-th eddy viscosity coefficient in the expansion and centered at the $j$-th element of $X_{\bm{\mu},t}$}
\nomenclature[G]{$w_{i,j}$}{The weight of the $j$-th RBF used in interpolating the $i$-th eddy viscosity coefficient in the expansion.}

\nomenclature[G]{$\mathbb{V}_{POD}$}{POD space for velocity}
\nomenclature[G]{$V_i$}{control volume in the mesh}
\nomenclature[G]{$\bm{u}_f$}{the velocity vector evaluated at the centre of each face of the control volume}
\nomenclature[G]{$(\bm{\nabla} \bm{u})_f$}{the gradient of $\bm{u}$ at the faces}
\nomenclature[G]{$\bm{u}_N$}{the velocity at the centre of one neighboring cell}
\nomenclature[G]{$\bm{u}_P$}{the velocity at the centre of one neighboring cell}
\nomenclature[G]{$\bm{Y}_L$}{the vector of observed outputs (the $L^2$ projection coefficients of the $L$-th viscosity mode onto the snapshots) for the interpolation procedure}

\nomenclature[G]{$\epsilon_u$}{the $L^2$ relative error between the FOM velocity field and a reduced order velocity field}
\nomenclature[G]{$\epsilon_p$}{the $L^2$ relative error between the FOM velocity field and a reduced order pressure field}
\nomenclature[G]{$\epsilon_{C_L}$}{the $L^2$ relative error between the FOM lift coefficient time signal and the reduced order reconstructed one}
\nomenclature[G]{$C_l$}{the lift coefficient which corresponds to the force component in the lift direction $L$}

\nomenclature[G]{$\bm{U_{BC}}$}{the vector of non-zero boundary velocity values which are imposed as non-homogeneous Dirichlet conditions at the Dirichlet boundary}

\nomenclature[G]{$\bm{D^k}$}{a vector in the penalty method for treatment of the boundary conditions at reduced order level, that contains the values of the $L^2$ norms of the velocity POD modes at the ${\Gamma_D}_k$ in the Dirichlet boundary}

\nomenclature[G]{$\bm{E^k}$}{a matrix in the penalty method for treatment of the boundary conditions at reduced order level, that contains the values of the $L^2$ scalar products of the velocity POD modes at the ${\Gamma_D}_k$ in the Dirichlet boundary}
\nomenclature[G]{$\bm{\delta}$}{a matrix calculated in the offline stage that represents the contribution of viscous forces acting on a surface in the domain}

\nomenclature[G]{$\bm{\theta}$}{a matrix calculated in the offline stage that represents the contribution of pressure forces acting on a surface in the domain}

\newcommand{\dif}{\mbox{d}}

\title[]{Data-Driven POD-Galerkin reduced order model
for turbulent flows}

\author{Saddam Hijazi\textsuperscript{1,*}}
\thanks{\textsuperscript{*}Corresponding Author.}
\address{\textsuperscript{1}SISSA, International School for Advanced Studies, Mathematics Area, mathLab Trieste, Italy.}
\email{shijazi@sissa.it}

\author{Giovanni Stabile\textsuperscript{1}}
\email{gstabile@sissa.it}
\author{Andrea Mola\textsuperscript{1}}
\email{amola@sissa.it}
\author{Gianluigi Rozza\textsuperscript{1}}
\email{grozza@sissa.it}


\date{}

\dedicatory{}


\begin{abstract} 
In this work we present a Reduced Order Model which is specifically designed to deal with turbulent flows in a finite volume setting. The method used to build the reduced order model is based on the idea of merging/combining projection-based techniques with data-driven reduction strategies. In particular, the work presents a mixed strategy that exploits a data-driven reduction method to approximate the eddy viscosity solution manifold and a classical POD-Galerkin projection approach for the velocity and the pressure fields, respectively. The newly proposed reduced order model has been validated on benchmark test cases in both steady and unsteady settings with Reynolds up to $Re=O(10^5)$.
\end{abstract}

\maketitle

%
%
%
%
\section{Introduction}\label{sec:intro}

A large part of physical problems (fluid dynamics, mechanics and heat transfer, ...) in relevant engineering and physics applications is governed by conservation laws. Over the years, several different numerical methods have been developed to solve the systems of Partial Differential Equations (PDEs) resulting from these conservation laws. Among these we mention the finite difference (FDM), the finite element (FEM), the finite volume (FVM), and the spectral element method (SEM). In particular, the Finite Volume one \cite{Moukalled:2015:FVM:2876154,Versteeg1995AnIT} is very often used to solve fluid dynamics and more in general hyperbolic problems.\par

Despite the recent increase of available computational power and new computational methods, the resolution of the governing equations, using one of the classical discretization methods previously mentioned, may become for several reasons not convenient. This is evident in common situations such as real-time control problems, where a small computational time is a major requirement or in a multi-query contest (e.g. optimization, uncertainty quantification, repetitive computational environment), where one needs to compute a certain output of interest for a large number of different input settings. This makes the cost of resorting to standard numerical methods (that will be referred as the Full Order Model (FOM)) prohibitive. These challenges in simulating computational problems has pushed the scientific community to seek techniques which could reduce the computational cost. Reduced Order Methods (ROMs) have been successful in meeting the needs of reducing the computational time offering high speed up rates. For a comprehensive review on ROMs, the reader may refer to \cite{hesthaven2015certified,quarteroniRB2016,bennerParSys,Benner2015,Bader2016}.\par


Projection based ROMs \cite{Balajewicz2012,Amsallem2012}, on which this article is focused,  have been applied in several scientific contributions dealing with laminar fluid dynamics problems and the methodology is already well established. On the other side, for what concerns turbulent flows, there are still several issues that need to be addressed. For instance, it is well known that projection based ROMs of turbulent flows suffer from energy stability issues \cite{Carlberg2017}. This is due to the fact that the POD retrieves the modes which are biased toward large, high-energy scales, but the turbulent small scales are the responsible scales for the dissipation of the turbulent kinetic energy \cite{Moin1998}.\par

Several strategies have been proposed to stabilize ROMs for turbulent flows and here a brief overview of the possible strategies is outlined. A possible approach suggests to include dissipation via a closure model, see \cite{wang_turb,Aubry1988}. In\cite{COUPLET2003}, it has been theoretically and numerically shown that the POD modes have similar energy transfer to the one of the Fourier modes. This suggests that the use of Large Eddy Simulations (LES) at the full order level could be beneficial in the case of POD-Galerkin-based ROMs.

Another possible approach \cite{Iollo2000} to obtain more dissipative ROMs, and justified by the fact that small scale modes have $H^1$ norm value that is higher than their $L^2$ norm value, proposes the usage of the $H^1$ inner product instead of the $L^2$ one in order to compute the POD modes.

Efforts to reduce CFD problems for turbulent flows include also employing minimum residual formulation in the reduced order model \cite{Carlberg2010,Carlberg2013,Tallet2015} or the use of the Dynamic Mode Decomposition (DMD) \cite{Alla,demo2018shape,dmdTezzele,le2017higher}. 

Recently in \cite{Fick2018}, the authors proposed a constrained formulation to deal with long time instabilities. In the latter work, a constrained Galerkin formulation is proposed in order to correct the standard Galerkin approach. The reduced order model in \cite{Fick2018} was generated using $H^1$-POD-$h$Greedy strategy, which is a simplified version of the $h$-type Greedy \cite{Eftang2011}. In \cite{Rebollo2017} the authors presented a reduced order model (based on the FEM) for the Smagorinsky turbulence model \cite{SMAGORINSKY1963} for steady flows. Their approach consisted into the approximation of the non-linear eddy diffusion term using the Empirical Interpolation Method. The contribution \cite{Ullman2010APR} presents a reduced order model which is designed also to deal with Smagorinsky turbulence model. The authors in \cite{Ullman2010APR} proposed a model which solves for the degrees of freedom of the velocity components and does not take into account pressure, the matrix coefficients, which come from the projection of the eddy viscosity term onto the velocity POD modes, have been assumed to be time dependent, and thus, these coefficients were dynamically updated during the time integration of the momentum equation at reduced order level. Additional works on Smagorinsky ROMs are presented in \cite{noackTR,Wang2011,noack2002}. In the context of ROMs for turbulent flows it is also worth mentioning the Variational Multi-Scale (VMS) method \cite{Bergmann2009, Stabile2019}. Smagorinsky VMS-ROMs are proposed in \cite{ChacnRebollo2018,BallarinChaconDelgadoGomezRozza2019}. In \cite{Yano2019} started from a Discontinuous Galerkin formulation to inherits the stability of the full order discretization. \par 

Most of the works mentioned above make use of Projection-based methods. However, ROMs can also be obtained by data-driven approaches  \cite{Ionita2014,Peherstorfer2015,Loiseau2018,ROWLEY2009,Kaiser2014,Guo2018,Hesthaven2018,NOACK2003}. A recent work on data-driven reduced order modeling for time-dependent problems can be found in  \cite{Guo2019}, where the authors proposed a regression based model to approximate the maps between the time-parameter values and the projection coefficients onto the reduced basis.\par

Since the final aim is to develop ROMs for flows with high Reynolds number, at the FOM level a Direct Numerical Simulation (DNS) is not affordable and thus we have to introduce turbulence modeling. In the FVM setting the most used techniques to introduce turbulence modeling are based on the Reynolds Averaged Navier--Stokes (RANS) equations and on the Large Eddy Simulation (LES) method. In this work, the RANS approach is considered. In order to solve the RANS equations a turbulence closure model that describes the effect of sub grid scales is required. In order to approximate the Reynolds stress tensor, we analyzed eddy viscosity closure models for both steady parametrized flows and unsteady flows. We considered closure models with both $k-\epsilon$ and SST $k-\omega$   \cite{Menter1994,Jones1972} which are two equations models, in which the eddy viscosity $\nu_t$ depends algebraically on two variables $k$ and $\epsilon$ or $\omega$. These variables stand respectively for the turbulent kinetic energy, turbulent dissipation and the specific turbulent dissipation rate. An additional PDE is solved for each of the turbulence variables.\par  

In this work we present a mixed approach between projection-based ROMs and data-driven-based ROMs, for some references on hybrid projection/data-driven ROM see \cite{Xie2018,GALLETTI2004,Couplet2005,Lu2017,noack2005,Peherstorfer2016}. In \cite{Xie2018} the authors presented a combination of projection based ROM with a Data Driven Filtering technique. In particular the work proposed to modify the standard Galerkin ROM by introducing a correction term which models the interaction between resolved modes and truncated modes. The authors used data driven modeling only to approximate the correction term, and tested the ROM on a 2D channel flow past a circular cylinder at Reynolds number of $100$, $500$ and $1000$.\par 
In \cite{GALLETTI2004,Couplet2005}, calibration methods have been constructed for the goal of reducing the Navier--Stokes equations, the authors used POD-Galerkin projection strategy and then they utilized data-driven techniques for calibrating the reduced order models. In \cite{GALLETTI2004}, this is done by assuming that the term which contains the  pressure gradient (in the projected momentum equation) is modeled by the product of a calibration matrix and the reduced vector of velocity coefficients. Afterwards, the calibration matrix entries can be found by minimizing a functional that depends on the values of the interpolated velocity vector of $L^2$ projection coefficients. In \cite{Couplet2005}, the calibration is done by finding the polynomial function that sets up the reduced dynamical system for the velocity coefficients as the solution to an optimization problem, where the functional which has to be minimized has two weighted terms. The first term measures the error between the values of the projection coefficients obtained from the data and the reduced solution of the dynamical system. The second imposes a cost for the difference between the original polynomial of the reduced dynamical system and the new one that determines the calibrated system.\par 

In \cite{noack2005}, the hybrid approach is similar to \cite{GALLETTI2004}, where an empirical pressure model is used to approximate the pressure term in the projected momentum equation. The data-driven approach utilized is a linear regression which fitted a set of coefficients in the empirical model from the data. In the last mentioned works, the hybrid/mixed approaches include modeling projected terms at the reduced order level and modifying the reduced order matrices entries. We mentioned only  works which focus on reducing the  the Navier--Stokes equations in both laminar and turbulent settings. Since such works were focused on reconstructing the velocity field of Direct Navier--Stokes resolutions, we here stress that the corresponding reduced model did not include the pressure field nor any turbulence associated field. In the present reduced approach we instead aim at reconstructing both the velocity and pressure fields and also consider the turbulent viscosity field $\nu_t$. This is motivated by the fact that the eddy viscosity is used at the reduced level to stabilize the momentum equation, as is the case for any FOM employing one or two equations turbulence model based on the Boussinesq eddy viscosity assumption. In fact, including the eddy viscosity in the ROM formulation introduces consistency with the FOM. Furthermore, the motivation behind the computation of a reduced version of the pressure field is that in several applications, important performance parameters not only depend on the velocity field, but also on the pressure one. Among these performance parameters, we mention for instance the fluid dynamics forces acting on the surface of a certain body. Thus, the ROM approach developed aims at approximating the fluid dynamics variables $\bm{u}$, $p$ and $\nu_t$. For such reason, separate sets of ROM coefficients are employed for the reduced order expansion of the $\bm{u}$, $p$ and $\nu_t$ fields. Yet, if the
pressure and velocity coefficients are determined through a well assessed projection methodology, the correct identification of the ROM coefficient for the turbulent variable $\nu_t$ is less obvious. Ideally, a proper projection procedure requires that the specific turbulence model equations used in the FOM solver must be taken into account. Unfortunately, given the wealth of one and two equations turbulence models of common use in the engineering community, their several variants and the even higher number of closure coefficients to be tracked at the ROM level, this approach appears not suitable for versatile ROMs which aim at being applied to FOM results obtained with different solvers. For instance, for solvers included in the OpenFOAM\textsuperscript{\textregistered} (OF) \cite{weller1998tensorial} library --- which are employed in this work --- users can access to about 20 RANS one or two equation turbulence models. This would not only require the development of a projection strategy for each turbulence equation encountered, but would also force constant monitoring of the FOM solvers libraries updates. In fact, even minimal changes in the turbulence models closure coefficients would make the results of the ROM solvers inaccurate. Hence, a versatile ROM solver, that can be employed in the every day virtual prototyping work by design engineers, should ideally be sensitive to the turbulence models used at the FOM level, but its implementation should not be dependent on their smallest details and intricacies. For such reason, we decided to use data-driven techniques for the computation of the reduced order coefficients of $\nu_t$, while still resorting to reduced order expansion of the $\bm{u}$ and $p$ fields. Indeed, such approach is able to reproduce differences due to changes in the particular turbulence model employed in the FOM simulations, while avoiding the increased ROM complexity due to the projection of the specific turbulence equations.\par

As a result, the approach developed in this work exploits the traditional projection methods in the part that computes the degrees of freedom for the reduced velocity and pressure fields. On the other hand, it uses a data-driven technique for the computation of the reduced coefficients of the eddy viscosity field. This is done by means of an interpolation process with Radial Basis Functions (RBF). The approach in the offline stage involves the construction of a RBF interpolant function (with Gaussian kernel functions) based on the set of samples used to train the ROM. In the general case of parametrized unsteady flows, both the coefficients obtained by the $L^2$ projection of the velocity snapshots (obtained by different values of the parameters and/or acquired at different time instants) onto the spatial modes of the velocity, as well as their vector derivatives, will be used to compute the weights of the RBF interpolant function. In the online stage, the values of the eddy viscosity coefficients are obtained by interpolation.  The dynamical system resulted from the projection step can be solved to obtain POD coefficients of the pressure and velocity expansion. To summarize, this approach is based on two main ideas. The first one is to approximate the solution manifold of the eddy viscosity field by means of an interpolation based approach. The second idea is to still exploit projection based methods to determine the expansion coefficients for velocity and pressure.

The work is organized as follows: section \ref{sec:FOM} deals with the description of the full order model and of the numerical methods used to solve the incompressible Navier--Stokes equations. Section \ref{sec:ROM} presents the methodologies used in this work to assemble the reduced order model. A review of projection based ROMs is outlined in \ref{sec:ProjectionROMs}, then the POD-Galerkin projection method is addressed in \ref{sec:POD}. Subsection \ref{sec:DDROM} focuses on the mixed projection-based/data-driven reduced order model. Subsection \ref{sec:bc} addresses how boundary conditions are treated at the reduced order level. The numerical examples are presented in \ref{sec:results} with two benchmark test cases which are the steady case of the backstep and the unsteady case of the flow past a circular cylinder. Conclusions and perspectives follow.\par



\section{The full order model (FOM)}\label{sec:FOM}
The present section is devoted to a description of the governing equations of the full order fluid dynamic model. Thus, the parametrized incompressible Navier--Stokes equations will be presented, along with details of their finite volumes discretization. Finally
the Reynolds Averaged Navier--Stokes equations will be presented, including  some relevant aspects of the turbulence modeling considered in this
work.

\subsection{The mathematical problem: parametrized Navier-Stokes equations}\label{sec:NS}
In this subsection, the strong form of the mathematical problem of interest is recalled. Given a parameter vector $\bm{\mu} \in \mathcal{P} \subset \mathbb{R}^q$, where $\mathcal{P}$ is a $q$-dimensional parameter space. The Navier-Stokes equations parametrized by $\bm{\mu}$ read as follows :
\begin{equation}\label{eq:navstokes}
\begin{cases}
\frac{\partial\bm{u}(t,\bm{x};\bm{\mu})}{\partial t}+ \bm{\nabla} \cdot (\bm{u}(t,\bm{x};\bm{\mu}) \otimes \bm{u}(t,\bm{x};\bm{\mu})) - \bm{\nabla} \cdot \nu \left(\bm{\nabla}\bm{u}(t,\bm{x};\bm{\mu})+\left(\bm{\nabla}\bm{u}(t,\bm{x};\bm{\mu})\right)^T\right)=-\bm{\nabla}p(t,\bm{x};\bm{\mu}) &\mbox{ in } \Omega \times [0,T],\\
\bm{\nabla} \cdot \bm{u}(t,\bm{x};\bm{\mu})=\bm{0} &\mbox{ in } \Omega \times [0,T],\\
\bm{u} (t,\bm{x};\bm{\mu}) = \bm{f}(\bm{x},\bm{\mu}) &\mbox{ on } \Gamma_{In} \times [0,T],\\
\bm{u} (t,\bm{x};\bm{\mu}) = \bm{0} &\mbox{ on } \Gamma_{0} \times [0,T],\\ 
(\nu\bm{\nabla} \bm{u} - p\bm{I})\bm{n} = \bm{0} &\mbox{ on } \Gamma_{Out} \times [0,T],\\ 
\bm{u}(0,\bm{x})=\bm{R}(\bm{x}) &\mbox{ in } (\Omega,0),\\            
\end{cases}
\end{equation}
where $\Gamma = \Gamma_{In} \cup \Gamma_0 \cup \Gamma_{Out}$ is the boundary of the fluid domain $\Omega \in \mathbb{R}^d$, with $d=1,$ $2$ or $3$. The boundary is formed by three different parts $\Gamma_{In}$, $\Gamma_{Out}$ and $\Gamma_0$ , which correspond respectively to the inlet boundary, the outlet boundary and the physical walls. $\bm{u}$ is the flow velocity vector field, $t$ is the time, $\nu$ is the fluid kinematic viscosity, and $p$ is the normalized pressure field, which is divided by the fluid density $\rho_f$, $\bm{f}$ is a generic function that describe the velocity on the inlet $\Gamma_{In}$ and it is parametrized through $\bm{\mu}$. $\bm{R}$ is the initial velocity field and $[0,T]$ is the time window under consideration. We remark that in this work the parameter $\bm{\mu}$ is always a physical parameter.\par

\subsection{The finite volume discretization}\label{sec:FV}

\begin{figure}
\centering
{
  \ifpdf
  \resizebox{0.5\textwidth}{!}{
    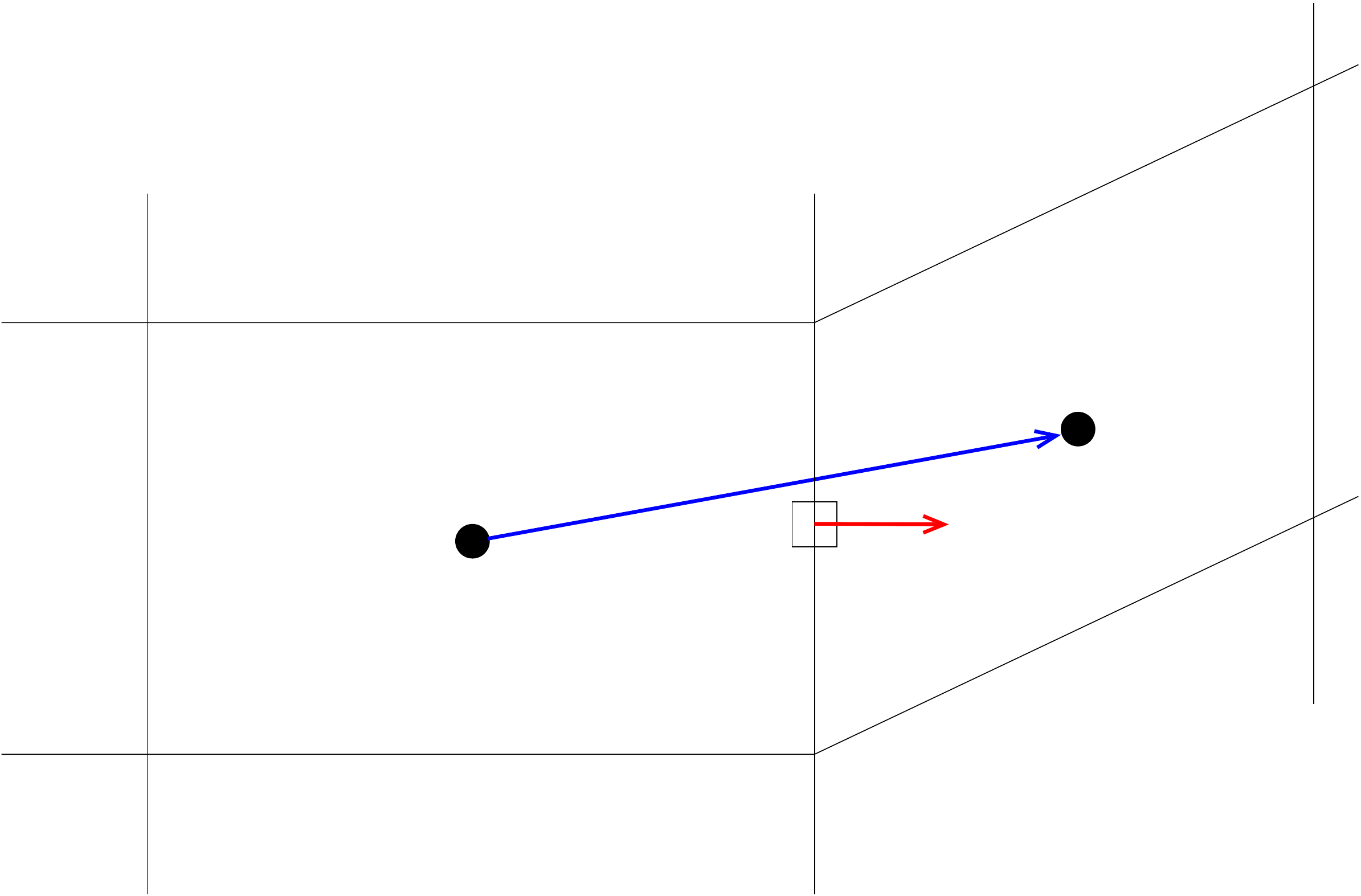
  }
  \else
  \resizebox{0.5\textwidth}{!}{
    \input{./figures_ok/fvm.pstex_t}
  }
  \fi
}
\caption{Sketch of a finite volume in 2 dimensions}\label{fig:FV_skatch}
\end{figure}

The governing equations of \eqref{eq:navstokes} are discretized using the FVM \cite{Moukalled:2015:FVM:2876154}. After choosing an appropriate polygonal tessellation, one can write the system of partial
differential equations \eqref{eq:navstokes} in integral form over each control volume. In the present work $2$-dimensional tessellations are considered. The number of degrees of freedom of the discretized problem represents the dimension of the full order model (FOM) which is denoted by $N_h$. In the next subsections, the discretization methodology of the momentum and continuity equations is addressed. In particular the momentum and continuity equations are solved using a segregated approach in the spirit of Rhie and Chow interpolation. The discretization starts  writing the momentum equation in integral form for each control volume $V_i$ as follows:
\begin{equation}\label{eq:mom_eq}
\int_{V_i} \frac{\partial}{\partial t}  \bm{u}dV + \int_{V_i} \bm{\nabla} \cdot (\bm{u} \otimes \bm{u}) dV - \int_{V_i} \bm{\nabla} \cdot \nu \left(\bm{\nabla}\bm{u}+\left(\bm{\nabla}\bm{u}\right)^T\right)
 dV +\int_{V_i} \bm{\nabla} p dV  = 0 . 
\end{equation}
We define then a generic cell center $P$ and a set of neighboring points around it $N$ (\autoref{fig:FV_skatch}). For each cell $P$, the discretized form of the momentum equation is then written as:
\begin{equation}\label{eq:disc_mom}
a_P^{\bm{u}} \bm{u}_P +\sum_N a_N^{\bm{u}} \bm{u}_N  = - \bm{\nabla} p ,
\end{equation}
where $\bm{u}_N$ and $\bm{u}_P$ are the velocities at the centers of two neighboring cells, $a_P^{\bm{u}}$ is the vector of diagonal coefficients of the equations and  $a_N^{\bm{u}}$ is the vector that consists off diagonal coefficients. \autoref{eq:disc_mom} is rewritten for all the cells in matrix form as:
\begin{equation}\label{eq:disc_mom_matr}
\mathcal{A} \bm u = \mathcal{H} - \bm{\nabla} p.
\end{equation}
In the above expression the terms $\mathcal{H} = -\sum_N a_N^{\bm{u}} \bm{u}_N$ and $\bm{\nabla} p$ are evaluated in an explicit manner based on previous tentative values of the velocity and pressure fields or on the values converged at the previous iteration or at the previous time step. 
The $\mathcal{A}$ matrix is a diagonal matrix and can be easily inverted and therefore \autoref{eq:disc_mom_matr} can be easily solved:
\begin{equation}\label{eq:vel_mom}
\bm u = \mathcal{A}^{-1} \mathcal{H} - \mathcal{A}^{-1} \bm{\nabla} p.
\end{equation}
If we apply the divergence operator and then exploit the continuity equation ($\bm{\nabla} \cdot \bm u = 0$) we obtain a Poisson equation for pressure:
\begin{equation}\label{eq:poisson_disc}
\bm{\nabla} \cdot (\mathcal{A}^{-1} \bm{\nabla} p) = \bm{\nabla} \cdot (\mathcal{A}^{-1} \mathcal{H}).
\end{equation}
The equation for pressure can be solved and used together with \autoref{eq:vel_mom} and with the discretized version of the continuity equation to update $F_f$ (the mass flux through each face of the control volume):
\begin{equation}\label{eq:mass_flux}
F_f = \bm u_f \cdot \bm S_f = - \mathcal{A}^{-1} \bm{S}_f \cdot \bm{\nabla} p + \mathcal{A}^{-1} \bm S_f \cdot \mathcal{H},
\end{equation}
where $\bm{S_f}$ is the area vector of each face of the control volume and $\bm{u_f}$ is the velocity vector evaluated at the center of each face of the control volume. The procedure used to discretize all the different terms inside the Navier-Stokes equations is explained in what follows. The pressure gradient term is discretized with the use of Gauss's theorem:
\begin{equation}\label{eq:grad_pressure}
\int_{V_i} \bm{\nabla} p dV = \int_{S_i} p d \bm{S} \approx \sum_f \bm{S_f} p_f,
\end{equation}
where $p_f$ is the value of pressure at the center of the faces (\autoref{fig:FV_skatch}).

Again using Gauss's theorem, the convective term can be discretized as follows:
\begin{equation}\label{eq:con_term} 
\int_{V_i} \bm{\nabla} \cdot (\bm{u} \otimes \bm{u}) dV  = \int_{S_i} ( d\bm{S} \cdot ( \bm{u} \otimes \bm{u}))  \approx  \sum_f  \bm{S_f} \cdot  \bm{u_f} \otimes \bm{ u_f} = \sum_f F_f \bm{u_f}.
\end{equation}
We remark that the velocity unknowns in the discretized form of the equations are always computed at the center of the faces. Therefore these values must be interpolated using the values at the cell centers. Several interpolation schemes are available such as central, upwind, second order upwind and blended differencing schemes. ${F_f}$ is the mass flux through each face of the control volume and, in order to remove the non-linearity, it is computed using the previous converged velocity and updated with~\autoref{eq:mass_flux} .\par 
The diffusion term is discretized as follows:
\begin{equation}\label{eq:diff_term}
\int_{V_i} \bm{\nabla} \cdot \nu \left(\bm{\nabla}\bm{u}+\left(\bm{\nabla}\bm{u}\right)^T\right)
 dV = \int_{S_i} d\bm{S} \cdot  \nu \left(\bm{\nabla}\bm{u}+\left(\bm{\nabla}\bm{u}\right)^T\right) \approx \sum_f \nu \bm{S_f} \cdot (\bm{\nabla} \bm{u})_f,    
\end{equation}
where $(\bm{\nabla} \bm{u})_f$ is the gradient of $\bm{u}$ at the faces. A procedure similar to the one described for pressure
in \eqref{eq:grad_pressure} is used to compute the value of $(\bm{\nabla} \bm{u})_f$. As for computing the term $\bm{S_f} \cdot (\bm{\nabla} \bm{u})_f$ in \eqref{eq:diff_term}, its value depends on whether the mesh is orthogonal or non-orthogonal. The mesh \ref{fig:FV_skatch} is orthogonal if the line that connects two cell centers is orthogonal to the face that divides these two cells. For orthogonal meshes the term $\bm{S_f} \cdot (\bm{\nabla} \bm{u})_f$ is computed as follows :
\begin{equation}
\bm{S_f} \cdot (\bm{\nabla} \bm{u})_f = |{\bm{S_f}}| \frac{\bm{u}_N - \bm{u}_P}{\bm{|d|}},
\end{equation}
where $\bm{u}_N$ and $\bm{u}_P$ are the velocities at the centers of two neighboring cells and $\bm{d}$ is the distance vector connecting the two cell centers see \autoref{fig:FV_skatch}. If the mesh is not orthogonal, then a correction term has to be added to the above equation. In that case, one has to consider computing a non-orthogonal term to account for the non-orthogonality of the mesh \cite{jasak1996error} as given by the following equation:
\begin{equation}
\bm{S_f} \cdot (\bm{\nabla} \bm{u})_f = |\bm{\Delta}| \frac{\bm{u}_N - \bm{u}_P}{\bm{|d|}} + \bm{J} \cdot (\bm{\nabla} \bm{u})_f ,
\end{equation}
where the following relation holds $\bm{S_f} = \bm{\Delta} + \bm{J}$. The first vector $\bm{\Delta}$ is chosen parallel to $\bm{S_f}$. The term $(\bm{\nabla} \bm{u})_f$ is obtained through interpolation of the the values of the gradient at the cell centers $(\bm{\nabla} \bm{u})_N$ and $(\bm{\nabla} \bm{u})_P$ in which the subscripts $N$ and $P$ indicate the values at the center of the cells of the two neighboring cells.
The coupled system of discretized equations given by ~\autoref{eq:vel_mom} and~\autoref{eq:poisson_disc} is solved by a segregated approach and specifically using the the SIMPLE \cite{PATANKAR19721787} algorithm for the steady case and the PIMPLE \cite{Moukalled:2015:FVM:2876154} algorithm for the unsteady case that merges the PISO \cite{Issa198640} and the SIMPLE \cite{PATANKAR19721787} algorithms.

\subsection{Turbulence modeling}\label{sec:turb}
Since the interest is to solve and and to reduce computational fluid dynamics problems characterized by high Reynolds numbers, the direct numerical resolution of the whole spectrum of temporal and spatial scales is not feasible.
In order to model turbulence without resolving all the temporal and spatial scales up to the Kolmogorov scale two main different approaches are typically used. The first approach --- the one considered in this work --- is based on Reynolds Averaged Navier-Stokes (RANS) equations and substantially consists into the decomposition of velocity and pressure fields into a mean part and a fluctuating part with zero mean. The decomposition of a generic scalar field $\sigma(x,t)$ will read as follows
\begin{equation}\label{eq:Reynolds_decompositon}
\sigma = \overline{\sigma} + \sigma^\prime, 
\end{equation}
where $\overline{\sigma}$ is the mean part and $\sigma^\prime$ is the fluctuating one. The RANS equations are obtained after introducing such decomposition for each scalar field (there are four scalar fields consisting of the three velocity components and the pressure field) into Navier-Stokes equations and time averaging them. In RANS, the approach is based on solving the equations for the mean part of each field after making use of the assumption that the fluctuating part has zero mean. A second possible approach --- which not considered in this work --- consists into Large Eddy Simulations (LES) \cite{berselli2005mathematics,sagaut2006large}. LES turbulence modeling is done by filtering and solving the Navier-Stokes equations just for specific scales which are the large scales.\par

\subsubsection{RANS equations}\label{sec:RANS}
In this subsection, the RANS equations will be presented in further detail. As mentioned earlier, in RANS the turbulence modeling starts by the Reynolds decomposition of the velocity and pressure fields into a mean part and a fluctuating one. These are denoted with  $\overline{\bm{u}}, \overline{{p}}$ for the mean part and $\bm{u}^\prime, {p}^\prime$ for the fluctuating part. Inserting the Reynolds decomposition into \eqref{eq:navstokes}, and time averaging the equations yields the so-called Reynolds Average Navier-Stokes (RANS) equations.\par 
Due to the the non-linearity of Navier-Stokes equations, the velocity fluctuations will not completely vanish in the time averaged equations. In particular the so called Reynolds stress tensor $\mathcal{R} = \overline{\bm{u}^\prime \bm{u}^\prime}$ is the single residual term in which the fluctuating components still appear after time averaging. Thus, such tensor must be expressed in terms of the mean part of the flow variables so as to obtain a closed problem for the latter unknowns. In this work we consider eddy viscosity models, that are based on Boussinesq assumption that the Reynolds stress tensor can be expressed by $\mathcal{R} =  \displaystyle{\frac{\nu_t}{2}[\bm{\nabla} \overline{\bm u} + (\bm{\nabla}\overline{\bm u})^T]}$. Different possibilities are available for the approximation of the additional coefficient $\nu_t$, which is named eddy viscosity \cite{boussinesq1877essa}. In the most effective cases the estimation of $\nu_t$ is based on the resolution of one or more additional transport-diffusion equations. We mention here the one equation Spalart–Allmaras (S–A) turbulence model \cite{SPALART1992} and the two equations $k-\epsilon$ \cite{Hanjalic1972} and SST $k-\omega$ turbulence models \cite{Menter1994}.\par

We here report the RANS equations for the $k-\omega$ turbulence model, which reads:

\begin{equation}\label{eq:RANS}
\begin{cases}
 \frac{\partial\overline{\bm{u}}}{\partial t}+ \bm{\nabla} \cdot (\overline{\bm{u}} \otimes \overline{\bm{u}})  = \bm{\nabla} \cdot \left[-\overline{p} \mathbf{I}+\left(\nu+\nu_t \right) \left(\bm{\nabla}\overline{\bm{u}}+\left(\bm{\nabla}\overline{\bm{u}}\right)^T\right)\right] &\mbox{ in } \Omega \times [0,T],\\
\bm{\nabla} \cdot \overline{\bm{u}} = 0  &\mbox{ in } \Omega \times [0,T],\\
\overline{\bm{u}} (t,\bm{x}) = \bm{f}(\bm{x},\bm{\mu}) &\mbox{ on } \Gamma_{In} \times [0,T],\\
\overline{\bm{u}} (t,\bm{x}) = \bm{0} &\mbox{ on } \Gamma_{0} \times [0,T],\\ 
(\nu\bm{\nabla} \overline{\bm{u}} - \overline{p}\bm{I})\bm{n} = \bm{0} &\mbox{ on } \Gamma_{Out} \times [0,T],\\ 
\overline{\bm{u}}(0,\bm{x})=\bm{R}(\bm{x}) &\mbox{ in } (\Omega,0),\\           
\nu_t=F(k,\omega), &\mbox{ in } \Omega, \\ 
\mbox{Transport-Diffusion equation for $k$}, \\ 
\mbox{Transport-Diffusion equation for $\omega$},\\ 
\end{cases}
\end{equation}
where $F$ is the function that describes the algebraic relationship between $\nu_t$ and the turbulence variables $k$ and $\omega$.

\section{The reduced order model (ROM)}\label{sec:ROM}

The proposed reduced order model is an extension of the model introduced in \cite{Stabile2018}. In \ref{sec:ProjectionROMs} the main notions of projection-based ROMs are recalled. Subsection \ref{sec:POD} introduces the POD technique and the general procedure used to construct a POD-Galerkin ROM. Subsection \ref{sec:DDROM} addresses in details how data-driven techniques can be exploited to stabilize ROMs for turbulent flows. In particular, the subsection \ref{sec:DDROM} explains how the model in the online stage uses data acquired in the offline stage for approximating the Reynolds stress term. Finally subsection \ref{sec:bc} outlines the  treatment of non-homogeneous boundary conditions at the reduced order level.
\subsection{Projection based ROMs}\label{sec:ProjectionROMs}
In the context of this work, we aim to develop ROMs which are able to approximate the solutions of Parametrized PDEs (PPDEs) in turbulent fluid dynamic problems efficiently and accurately. Reduced order modeling for PPDEs is based on the assumption that the solution field lives in a low dimensional manifold \cite{hesthaven2015certified}. Based on this assumption any element of the solution manifold can be approximated by the linear combination of a reduced number of global basis functions. The velocity and pressure fields can be approximated as a linear combination of the dominant modes (basis functions) multiplied by scalar coefficients. The modes are assumed to be dependent on space variables only, while the coefficients are allowed to have temporal and/or parameter dependency. The last statement leads to the following approximation of the fields:
\begin{equation}\label{eq:decompose}
\bm{u}(\bm{x},t; \bm{\mu}) \approx  \sum_{i=1}^{N_u} a_i(t;\bm{\mu}) \bm{\phi}_i(\bm{x}), \quad p(\bm{x},t; \bm{\mu})\approx \sum_{i=1}^{N_p} b_i (t;\bm{\mu}) {\chi_i}(\bm{x}),
\end{equation}   
where $\bm{\phi}_i(\bm{x})$ and ${\chi_i}(\bm{x})$ (which do not depend on $\bm{\mu}$ and $t$) are the spatial modes for velocity and pressure, respectively, $a_i(t;\bm{\mu})$ and $b_i (t;\bm{\mu})$ are temporal coefficients which depend on time $t$ and on the parameter vector $\bm{\mu}$. The reduced basis spaces $ \mathbb{V}_{rb}=\mbox{span}\left\{\bm{\phi}_i\right\}_{i=1}^{N_u}$ and $ Q_{rb}=\mbox{span}\left\{\chi_i\right\}_{i=1}^{N_p}$ can be obtained either by Reduced Basis (RB) method with a greedy approach \cite{hesthaven2015certified}, using Proper Orthogonal Decomposition (POD) \cite{Stabile2018}, by the Proper Generalized Decomposition (PGD) \cite{Dumon20111387,Chinesta2011}, or by Dynamic Mode Decomposition (DMD) \cite{Schmid2010}. For unsteady PPDEs, a POD-Greedy approach (POD in time and RB method with greedy algorithm in parameter space) can be used as in \cite{Haasdonk2008} or a nested POD can be used where POD is applied on time and later on parameter space. In this work, in order to calculate the reduced basis functions, we rely on a POD approach applied onto the full snapshots matrices formed by the fields obtained for different values of the parameters as well as for different time instants.

\subsection{POD-Galerkin projection method for laminar flows}\label{sec:POD}
One of the most used approaches to construct reduced order spaces is the proper orthogonal decomposition (POD) \cite{volkwein2013proper,Bergmann2009516,Baiges2014189,Burkardt2006,Ballarin2016}. The POD is a method to compress a set of numerical realizations (in the time or parameter space) into a reduced number number of orthogonal basis (modes) that capture the most important information when suitably combined. As mentioned above, in this work the POD is applied on a group of different realizations which are called snapshots. The POD modes are optimal in the sense that, for every number of chosen modes, the difference between the $L^2$ projection of the snapshots onto the modes and the snapshots themselves is minimized. In this setting, it has to be remarked that the FOM presented in \ref{sec:FOM} is solved for each value of the vector parameter $\bm{\mu} \in \mathcal{P}_M = \{\bm{\mu}_1,...\bm{\mu}_M\} \subset \mathcal{P}$ where $\mathcal{P}_M$ is a finite set of samples inside the parameter space $\mathcal{P}$. In case of non-stationary problems, while generating the snapshots for constructing the reduced order space,  one has to consider time and parameter dependencies. Consequently, for each parameter value one has the time instants $\{t_1,t_2,...,t_{N_T}\} \subset [0,T]$ at which snapshots are taken. For this reason, there will be a total number of snapshots $N_s = M * N_T$. The snapshots matrices $\bm{{\mathcal{S}_u}}$ and $\bm{{\mathcal{S}_p}}$, for velocity and pressure respectively, will be given by:
\begin{equation}\label{eq:snapU}
\bm{{\mathcal{S}_u}} = \{\bm{u}(\bm{x},t_1;\bm{\mu}_1),...,\bm{u}(\bm{x},t_{N_T};\bm{\mu}_M)\} ~\in~\mathbb{R}^{N_u^h\times N_s},
\end{equation}
\begin{equation}\label{eq:snapP}
\bm{{\mathcal{S}_p}} = \{p(\bm{x},t_1;\bm{\mu}_1),...,p(\bm{x},t_{N_T};\bm{\mu}_M)\} ~\in~\mathbb{R}^{N_p^h\times N_s},
\end{equation}
where $N_u^h$ and $N_p^h$ are the degrees of freedom for velocity and pressure fields, respectively.
The POD space for velocity is constructed by solving the following optimization problem:
\begin{equation}\label{eq:min1}
\mathbb{V}_{POD}= \mbox{arg min} \frac{1}{N_{s}} \sum_{n=1}^{N_s} ||  \bm{u}_n- \sum_{i=1}^{N_u}(\bm{u}_n, \bm{\phi}_i )_{L^2(\Omega)} \bm{\phi}_i ||_{L^2(\Omega)}^2,
\end{equation}
where $\bm{u}_n$ is a general snapshot of the velocity field which is obtained for any  value of the parameter $\bm{\mu}$ and acquired at any time instant $t_i$. It can be shown that solving \eqref{eq:min1} is equivalent to solve the following eigenvalue problem \cite{Kunisch2002492} :
\begin{equation}
\bm{C^u}\bm{V^u}=\bm{V^u}\bm{\lambda^u},
\end{equation}
where $\bm{C^u}~\in~\mathbb{R}^{N_s\times N_s}$ is the correlation matrix of the velocity field snapshot matrix $\bm{{\mathcal{S}_u}}$, $\bm{V^u}~\in~\mathbb{R}^{N_s\times N_s}$ is the matrix whose columns are the eigenvectors, $\bm{\lambda^u}$ is a diagonal matrix whose diagonal entries are the eigenvalues. The entries of the correlation matrix are defined as follows:
\begin{equation}\label{eq:cor_matrix}
(\bm{C^u})_{ij}=\left(\bm{u_i},\bm{u_j}\right)_{L^2(\Omega)}.
\end{equation}
One can compute the velocity POD modes as follows \cite{Stabile2017},
\begin{equation}
\bm{\phi}_i=\frac{1}{N_s \lambda^u_{i}} \sum_{j=1}^{N_s} \bm{u}_j V^u_{ij},
\end{equation}
similar procedure can be followed for the computation of the POD pressure modes $[{\chi_i}(\bm{x})]_{i=1}^{N_p}$.\par
After computing the POD modes of velocity and pressure, one can perform a Galerkin projection of the governing equations onto the POD space. Projecting the momentum equation of \eqref{eq:navstokes} onto the POD space spanned by the velocity POD modes yields:
\begin{equation}\label{eq:l2proj_vel}
\left( \bm{\phi}_i,\frac{\partial\bm{u}}{\partial t} +  \bm{\nabla} \cdot (\bm{u} \otimes \bm{u}) - \bm{\nabla} \cdot \nu \left(\bm{\nabla}\bm{u}+\left(\bm{\nabla}\bm{u}\right)^T\right) +\bm{\nabla} p \right)_{L^2(\Omega)} = 0 .
\end{equation}
Inserting the approximations \eqref{eq:decompose} into \eqref{eq:l2proj_vel} gives the following system:
\begin{equation}\label{eq:dyn_sys_lam}
\bm{\dot{a}}= \nu\bm{B} \bm{a} - \bm{a}^T \bm{C} \bm{a}-\bm{H}\bm{b},
\end{equation}
where $\bm{a}$ and $\bm{b}$ are the reduced vectors of coefficients $a_i(t;\bm{\mu})$ and $b_i(t;\bm{\mu})$, respectively, while the reduced matrices $\bm B, \bm C$ and $\bm H$ are computed as follows:
\begin{align}
& (\bm{B})_{ij}=\left( \bm{\phi}_i ,\bm{\nabla} \cdot \bm{\nabla}\bm{\phi}_j\right)_{L^2(\Omega)}, \\
& (\bm{C})_{ijk}=\left( \bm{\phi}_i , \bm{\nabla} \cdot (\bm{\phi}_j \otimes \bm{\phi}_k)\right)_{L^2(\Omega)} , \label{eq:div_phi}\\
& (\bm{H})_{ij} = \left( \bm{\phi}_i , \bm{\nabla} \chi_j \right)_{L^2(\Omega)}.
\end{align}

In \cite{Stabile2018}, one can find more details on the treatment of the non-linear term in Navier-Stokes equations. An important remark is that the system \eqref{eq:dyn_sys_lam} has $N_u + N_p$ unknowns but just $N_u$ equations. Therefore one must seek $N_p$ additional equations in order to close the system. It is not possible to directly exploit the  continuity equation at this stage because the velocity snapshots are divergence free and so are the velocity POD modes. The additional equations could be obtained by the usage of a Poisson equation for pressure also at the reduced order level, see \cite{Stabile2017}. Another possible approach is to employ a supremizer enrichment technique \cite{Ballarin2015,Rozza2007} where the velocity POD space is enriched with additional, non divergence-free modes in order to satisfy a reduced version of the inf-sup condition. We refer to \cite{Stabile2018} for the implementation of this approach in the finite volume setting. There exist also other approaches to obtain pressure-stable ROMs, for example the use of Pressure Stabilized Petrov-Galerkin (PSPG) methods during the online procedure \cite{Baiges2014189,Caiazzo2014598} or ROMs based on the assumption that velocity and pressure expansions share the same scalar coefficients \cite{Bergmann2009516,Lorenzi2016}.\par
In this work, the supremizer stabilization method has been chosen. This approach will ensure that velocity POD modes are not all divergence free so one can project the continuity equation onto the space spanned by the POD pressure modes. This will give the following reduced system: 
\begin{eqnarray}\label{eq:reduced_system}
\left\{
\begin{matrix}
\bm{M} \bm{\dot{a}}= \nu\bm{B} \bm{a} - \bm{a}^T \bm{C} \bm{a}-\bm{H}\bm{b},
\\
\bm{P} \bm{a}= \bm{0},
\end{matrix}
\right.
\end{eqnarray}
where the new reduced matrices $\bm{M}$ and $\bm{P}$ are the mass matrix, that due to the additional supremizer modes is not anymore unitary, and the matrix associated with the continuity equation. The entries of the two additional matrices are given by:
\begin{align}
& (\bm{M})_{ij}=\left( \bm{\phi}_i , \bm{\phi}_j \right)_{L^2(\Omega)},\\
& (\bm{P})_{ij}=\left(\chi_i, \bm{\nabla} \cdot  \bm{\phi}_j \right)_{L^2(\Omega)}.
\end{align}

\subsection{POD-Galerkin Mixed-ROM for turbulent flows}\label{sec:DDROM}
In this subsection, the attention will be shifted to flows characterized by high Reynolds number. As mentioned earlier, turbulence modeling at the full order level is resolved using the RANS equations with a proper closure model \eqref{eq:RANS}. This motivated the development of a reduced order model specifically tailored to turbulent flows. This model will be referred to from now on as Mixed-ROM.\par

A possible approach could consist into a POD procedure applied also onto the additional turbulence variables ($k, \omega, \epsilon$) as it was done with the velocity and pressure fields in \eqref{eq:decompose}. This phase should be followed by a POD-Galerkin projection of the additional transport diffusion equations that define the specific turbulence model in order to obtain a reduced version of the equations. The last step, in the online phase, would consist into the coupling of all the "reduced" equations and into their simultaneous resolution. The reduced equations come from the momentum equation, the continuity equation and the additional PDEs of the turbulent model. However, this approach has some drawbacks:
\begin{itemize}
\item it implies that the reduced order model needs to be customized to the specific turbulence model used during the offline stage;
\item since it requires also the projection of the PDEs of the turbulent equations, the effort in the generation of the reduced order model and the number of reduced unknowns is increased.
\end{itemize} 
Since one of the aims of this work is to develop a ROM which is "independent" from the turbulence model used to generate the FOM snapshots the latter approach is ruled out.
The chosen approach involves the extension of the assumption of the reduced order expansion only for the eddy viscosity without considering the additional turbulence variables ($k, \epsilon$ or $\omega$). In more details, this means introducing the reduced order eddy viscosity as a sum of eddy viscosity POD modes multiplied by temporal or parameter dependent coefficients. The eddy viscosity modes are computed using a POD approach and, during the online stage, the scalar coefficients of the POD expansion are computed with a data-driven approach that uses interpolation with Radial Basis Functions (RBF) \cite{Lazzaro2002,Micchelli1986}, thus the reduced order viscosity reads as follows:
\begin{equation}\label{eq:visc}
\nu_t(\bm{x},t; \bm{\mu})\approx  \sum_{i=1}^{N_{\nu_t}} g_i (t,\bm{\mu}) {\eta_i}(\bm{x}),
\end{equation}
where ${\eta_i}(\bm{x})$ are the POD modes for the eddy viscosity field and $g_i (t,\bm{\mu})$ are the scalar coefficients of the POD expansion. One can see that the temporal coefficients in the above equation are not the same of neither the ones of the velocity $a_i(t,\bm{\mu})$ nor the ones of the pressure $b_i (t,\bm{\mu})$. The data-driven approach will be used for the computation of these coefficients. The momentum equation of the RANS \eqref{eq:RANS} is projected onto the spatial modes of velocity, inserting also the POD decomposition of the eddy viscosity field \eqref{eq:visc}. On the other hand, the continuity equation is projected onto the pressure modes with the usage of supremizer enrichment. The POD-Galerkin projection will result in the following reduced system:
\begin{equation}\label{eq:PODsys}
\begin{cases}
\bm{M} \bm{\dot{a}}= \nu (\bm{B}+\bm{B_T})\bm{a}-\bm{a^T}\bm{C}\bm{a} +\bm{g^T}(\bm{C_{T1}}+\bm{C_{T2}})\bm{a}-\bm{H}\bm{b},\\
\bm{P} \bm{a}= \bm{0},\\ 
\end{cases}
\end{equation}
where $\bm{g}$ is the vector of the coefficients $[g_i (t,\bm{\mu})]_{i=1}^{N_{\nu_t}}$, and the new terms with respect to the dynamical system in \eqref{eq:reduced_system} are computed as follows: 
\begin{align}
& (\bm{B_T})_{ij}=\left( \bm{\phi}_i ,\bm{\nabla} \cdot (\bm{\nabla} \bm{\phi}_j^T)\right)_{L^2(\Omega)},\\
& (\bm{C_{T1}})_{ijk}=\left( \bm{\phi}_i ,\eta_j  \bm{\nabla} \cdot \bm{\nabla}\bm{\phi}_k  \right)_{L^2(\Omega)} , \\
& (\bm{C_{T2}})_{ijk}=\left( \bm{\phi}_i ,\bm{\nabla} \cdot \eta_j (\bm{\nabla} \bm{\phi}_k^T)\right)_{L^2(\Omega)}. 
\end{align}
As one can notice, system \eqref{eq:PODsys} has more unknowns $\bm{a}$, $\bm{b}$ and $\bm{g}$ than the available equations. This problem can be resolved by finding a proper way to compute the coefficients of the eddy viscosity POD expansion $\bm{g}$. This is carried out with the usage of a POD-I approach \cite{Wang2012,Walton2013,Salmoiraghi2018} using radial basis functions.\par
Before explaining more details about the used methodology we fix a set of notations and conventions. Let $X_{\bm{\mu},t}$ be the set defined as follows:
\begin{equation}\label{eq:set}
X_{\bm{\mu},t} = \mathcal{P}_M \bigtimes \{t_1,t_2,...,t_{N_T}\},
\end{equation}
$X_{\bm{\mu},t}$ is the Cartesian product of the discretized parameter set and the set of time instants at which snapshots were taken. This set has a cardinality of $N_s$ and its $i$-th member will be referred as $\bm{x}_{\bm{\mu},t}^i$. We remark that for each term $\bm{x}_{\bm{\mu},t}^i$ there is a corresponding unique snapshot (for $\bm{u}$, $p$ and $\nu_t$) that is used to compute the reduced basis for each variable in the offline stage. On the other hand, we define the parameter sample $\bm{\mu}^*$ as the one introduced to the reduced order model in the online stage. A remark has to be made that $\bm{\mu}^*$ should be close enough in the parameter space to the parameter samples used in the offline stage that will assure an accurate ROM result. Also we define $t^*$ as the time instant at which the Mixed-ROM solution is sought, where  $t_1 \leqslant t^* \leqslant t_{N_T}$. The last statement essentially means that currently it is not possible to extrapolate in time. Also we define $\bm{z}^* = (t^*,\bm{\mu}^*)$ as the combination of the online parameter sample and the time instant at which the Mixed-ROM solution is desired.\par
As done for velocity and pressure in \eqref{eq:snapU} and \eqref{eq:snapP}, respectively, we define a matrix of snapshots for the eddy viscosity field as follows:
\begin{equation}
\bm{{\mathcal{S}_{\nu_t}}} = \{\nu_t(\bm{x},t_1;\bm{\mu}_1),...,\nu_t(\bm{x},t_{N_T};\bm{\mu}_M)\} ~\in~\mathbb{R}^{N^h_{{\nu_t}} \times N_s},
\end{equation}  
where the $i$-th column of the $\bm{{\mathcal{S}_{\nu_t}}}$ represents an eddy viscosity snapshot and is denoted by $\bm{{\mathcal{S}^i_{\nu_t}}}$. We define $g_{r,l}$ as the coefficient computed from the $L^2$ projection of the $r$-th eddy viscosity snapshot $\bm{{\mathcal{S}^r_{\nu_t}}}$ onto the $l$-th eddy viscosity mode $\eta_l$.
\begin{equation}\label{eq:visc_coeff}
g_{r,l} = (\bm{{\mathcal{S}^r_{\nu_t}}},\eta_l)_{L^2(\Omega)}, \quad \text{for} \quad r=1,2,...,N_s \quad \text{and} \quad l =1,2,...,N_{\nu_t}.
\end{equation}
The interpolation statement will be the following: given the set $X_{\bm{\mu},t}$, the corresponding eddy viscosity snapshots $[\bm{{\mathcal{S}^i_{\nu_t}}}]_{i=1}^{N_s}$ and the coefficients $[g_{r,l}]_{r=1,l=1}^{N_s,N_{\nu_t}}$, predict the value of the vector $\bm{g}$ in \eqref{eq:PODsys} for the vector $\bm{z}^*$ defined earlier. The goal can be split to each of the scalar coefficients $[g_i (t^*,\bm{\mu}^*)]_{i=1}^{N_{\nu_t}}$. Meaning that the interpolation will be done separately $N_{\nu_t}$ times for each one of the scalar coefficients. From now on, we will include the dependency as follows $\bm{g}(\bm{z}^*)$ or $[g_i(\bm{z}^*)]_{i=1}^{N_{\nu_t}}$.\par
The interpolation procedure will be carried out for each mode separately, therefore one could fix the viscosity mode in \eqref{eq:visc_coeff} to be $\eta_{L}$, and then the vector $\bm{Y}_L = [g_{r,L}]_{r=1}^{N_s} \in \mathbb{R}^{N_s}$ is considered as the set of observations. The next step is to consider the pair of data $(X_{\bm{\mu},t},\bm{Y}_L)$ which is obtained in the offline stage by doing the computations in \eqref{eq:visc_coeff}. The objective is to approximate the value of the scalar coefficient $g_L(\bm{z}^*)$.\par  
The interpolation using RBF functions is based on the following formula :
\begin{equation}
G_L(\bm{z}) = \sum_{j=1}^{N_s} w_{L,j} \zeta_{L,j}(\norm{\bm{z}-\bm{x}_{\bm{\mu},t}^j}_{L^2(\mathbb{R}^{q+1})}), \quad \text{for} \quad L = 1,2,...,N_{\nu_t},
\end{equation}
where $\bm{z} = (t,\bm{\mu})$ with $\bm{\mu} \in \mathcal{P}$ and $t \in [0,T]$, $w_{L,j}$ are some appropriate weights and $\zeta_{L,j}$ for $j=1,...,N_s$ are the RBF functions which are chosen to be Gaussian functions, $\zeta_{L,j}$ is centered in $\bm{x}_{\bm{\mu},t}^j$. For the computation of the weights, the following property has to be used, which essentially comes from the data of the FOM:
\begin{equation}
G_L(\bm{x}_{\bm{\mu},t}^i) = g_{i,L} , \quad \text{for} \quad i = 1,2,...,N_s,
\end{equation}
and then it follows that,
\begin{equation}
\sum_{j=1}^{N_s} w_{L,j} \zeta_{L,j}(\norm{\bm{x}_{\bm{\mu},t}^i-\bm{x}_{\bm{\mu},t}^j}_{L^2(\mathbb{R}^{q+1})}) = g_{i,L} , \quad \text{for} \quad i = 1,2,...,N_s.
\end{equation}
The last equation can be rewritten as a linear system, namely:
\begin{equation}
\bm{A}_L^{\zeta} \bm{w}_L = \bm{Y}_L,
\end{equation}
where $(\bm{A}_L^{\zeta})_{ij} = \zeta_{L,j}(\norm{\bm{x}_{\bm{\mu},t}^i-\bm{x}_{\bm{\mu},t}^j}_{L^2(\mathbb{R}^{q+1})}$, one can solve the latter linear system to obtain the weights $\bm{w}_L$, which will be stored to be then used in the online stage.\par 
In the \textbf{Online Stage}, as \textbf{Input} we have the new time-parameter vector $\bm{z}^*$ and the goal is to compute $\bm{g}(\bm{z}^*) = [g_i(\bm{z}^*)]_{i=1}^{N_{\nu_t}}$, which is done simply by:
\begin{equation}
g_i(\bm{z}^*) \approx G_i(\bm{z}^*) = \sum_{j=1}^{N_s} w_{i,j} \zeta_{i,j}(\norm{\bm{z}^*-\bm{x}_{\bm{\mu},t}^j}_{L^2(\mathbb{R}^{q+1})}), \quad \text{for} \quad i=1,2,...,N_{\nu_t}.
\end{equation}
To summarize the procedure, the interpolation using RBF is done in the online stage. The procedure consists into separated $N_{\nu_t}$ times interpolation tasks for the interpolation of the elements of the vector $\bm{g}$, which appears in \eqref{eq:PODsys} for some value of the combined time-parameter vector $\bm{z}^*$.\par
The interpolation problem has as input a set of known data called $X_{\bm{\mu},t}$ with cardinality of $N_s$. A member in that set is a vector called $\bm{x}_{\bm{\mu},t}^i$ and lies in $\mathbb{R}^{q+1}$, where one can see that basically time has been treated as another parameter. The other discrete set of outputs (which has the same cardinality $N_s$) is the set of the coefficients obtained by the projection mentioned in \eqref{eq:visc_coeff} with the viscosity mode being fixed. At the end, based on the observations given in the offline stage, the value of the coefficient $g_i(\bm{z}^*)$ (the interpolant) will be approximated.\par 
The approach above is general for unsteady parametrized cases, a description of the same approach but just for steady cases can be found in \cite{FEF}. A modified version of this approach for unsteady flows may involve the splitting of the eddy viscosity field into two parts. The first part, for each individual parameter sample describes the time averaged viscosity field and the second one contains the time varying contribution. In other words we assume that the eddy viscosity field can be rewritten as follows:
\begin{equation}\label{eq:visc_ave}
\nu_t(\bm{x},t; \bm{\mu}) =  \overline{\nu_t}(\bm{x}; \bm{\mu}) + \nu^\prime_t(\bm{x},t;\bm \mu).
\end{equation}
Such decomposition is justified by the fact that usually the part which is largely affected by parameter changes is the mean contribution. Small fluctuations are in fact poorly affected by parameter variations and in our numerical example we have noticed that excluding the parameter value from the RBF approximation, $\nu^\prime_t(\bm{x},t;\bm \mu) \approx \nu^\prime_t(\bm{x},t)$, does not lead to a degradation of the accuracy. 
Therefore, this approach with such an approximation we have the advantage of splitting the time and parameter contributions into two separate terms which will ease the function approximation. At this point the reduced order approximation of the eddy viscosity field will be modified as follows:
\begin{equation}\label{eq:visc_ave_ROM}
\nu_t(\bm{x},t; \bm{\mu})\approx \sum_{i=1}^{M} \overline{g}_i (\bm{\mu}) \overline{\eta}_i(\bm{x}) + \sum_{i=1}^{N_{\nu_t}} g_i (t) {\eta_i}(\bm{x}),
\end{equation}
where the first sum approximates the averaged part while the second one approximates the time varying contribution. The fields $[\overline{\eta}_i]_{i=1}^{M}$ are given by the time averaged eddy viscosity fields of the $M$ different parameter samples (each field is computed as the time averaged field of only the eddy viscosity snapshots corresponding to one parameter sample), while $[\overline{g}_i]_{i=1}^{M}$ are the parameter dependent coefficients which, for a parameter sample outside of the training set $\bm{\mu}^*$, can be approximated by an interpolation procedure. The dynamical system \ref{eq:PODsys} is modified as follows:
\begin{equation}\label{eq:PODsysAve}
\begin{cases}
\bm{M} \bm{\dot{a}}= \nu (\bm{B}+\bm{B_T})\bm{a}-\bm{a^T}\bm{C}\bm{a} + \overline{\bm{g}}^T(\overline{\bm{C}}_{T1}+\overline{\bm{C}}_{T2})\bm{a} +\bm{g^T}(\bm{C_{T1}}+\bm{C_{T2}})\bm{a}-\bm{H}\bm{b},\\
\bm{P} \bm{a}= \bm{0},\\ 
\end{cases}
\end{equation}
where there are two new tensor terms defined as follows:
\begin{align}
& ({\overline{\bm{C}}_{T1}})_{ijk}=\left( \bm{\phi}_i ,\overline{\eta}_j  \bm{\nabla} \cdot \bm{\nabla}\bm{\phi}_k  \right)_{L^2(\Omega)} , \\
& ({\overline{\bm{C}}_{T2}})_{ijk}=\left( \bm{\phi}_i ,\bm{\nabla} \cdot \overline{\eta}_j (\bm{\nabla} \bm{\phi}_k^T)\right)_{L^2(\Omega)}. 
\end{align}
\par
One of the drawbacks of the approach in the current setting is that for unsteady cases one can not extrapolate in time. In order to address this issue the RBF interpolation can be rewritten in a different fashion. The idea is to change the independent variable of the RBF interpolation from being the time value $t^*$ to the combination of the reduced order velocity coefficients vectors $\bm{a}$ and $\bm{\dot{a}}$. The motivation comes from the fact that the eddy viscosity field $\nu_t$ at time $t^n$, denoted hereafter by ${\nu^n_t}$, is a function of the time history of the velocity field $\bm{u}$, in other words ${\nu^n_t} = {\nu_t}(\bm{u}^1,\bm{u}^2,...,\bm{u}^n)$, in the last formula $\bm{u}^n$ is the FOM velocity field obtained at time $t^n$. This allows us to write the eddy viscosity coefficients vector in the expansion \ref{eq:visc_ave_ROM} as follows\footnote{This expression aim to mimic the dependency between the eddy viscosity field and the velocity field. At the FOM level this dependency is described by a PDE. Therefore the expression is an approximation that could be extended with further terms in order to get closer to the map described by the PDE. However, for the numerical examples considered in this work such an approximation turned out to be sufficient.}:
\begin{equation}\label{eq:inter_v}
\bm{g}^n = \bm{g}^n(t^n) \approx \bm{g}^n(\bm{a}^n,\bm{\dot{a}}^n).
\end{equation}

The training phase of the RBF in the offline stage is done with the $L^2$ projection coefficients of the velocity modes (excluding the supremizer modes) onto the snapshots as well as the vector of time derivatives of these coefficients. In order to establish a clear idea of the training methodology, the following notation will be used:
\begin{equation}
\bm{{\mathcal{S}_u}} = \begin{bmatrix} 
    \bm{{\mathcal{S}_{\bm{\mu}_1,u}}}  \\
    \bm{{\mathcal{S}_{\bm{\mu}_2,u}}}  \\
    \vdots  \\
    \bm{{\mathcal{S}_{\bm{\mu}_M,u}}} 
    \end{bmatrix},
    \bm{{\mathcal{S}_p}} = \begin{bmatrix} 
    \bm{{\mathcal{S}_{\bm{\mu}_1,p}}}  \\
    \bm{{\mathcal{S}_{\bm{\mu}_2,p}}}  \\
    \vdots  \\
    \bm{{\mathcal{S}_{\bm{\mu}_M,p}}} 
    \end{bmatrix},
    \bm{{\mathcal{S}_{\nu_t}}} = \begin{bmatrix} 
    \bm{{\mathcal{S}_{\bm{\mu}_1,{\nu_t}}}}  \\
    \bm{{\mathcal{S}_{\bm{\mu}_2,{\nu_t}}}}  \\
    \vdots  \\
    \bm{{\mathcal{S}_{\bm{\mu}_M,{\nu_t}}}} 
    \end{bmatrix},
\end{equation}
where the snapshots matrices for all the variables have been written as $M$ vertically aligned submatrices with each one of the submatrices containing the time snapshots corresponding to a single sample. The next step is to define the $L^2$ velocity projection coefficients denoted by $\bm{a}_{\bm{\mu}_k,L^2}^r  ~\in~\mathbb{R}^{N_u}$:
\begin{equation}\label{eq:vec_coeff}
\bm{a}_{\bm{\mu}_k,L^2}^r = [ (\bm{{\mathcal{S}^r_{\bm{\mu}_k,u}}}, \bm{\phi}_1)_{L^2(\Omega)},..., (\bm{{\mathcal{S}^r_{\bm{\mu}_k,u}}}, \bm{\phi}_{N_u})_{L^2(\Omega)} ], \quad \text{for} \quad r=1,2,...,N_T, \quad k=1,2,...,M.
\end{equation}
Let
\begin{equation}
\bm{{\mathcal{A}_{1,k}}} = \begin{bmatrix} 
    \bm{a}_{\bm{\mu}_k,L^2}^1  \\
    \bm{a}_{\bm{\mu}_k,L^2}^2  \\
    \vdots  \\
    \bm{a}_{\bm{\mu}_k,L^2}^{N_T-1} 
    \end{bmatrix} ~\in~\mathbb{R}^{(N_T-1)\times N_u},
    \bm{{\mathcal{A}_{2,k}}} = \begin{bmatrix} 
    \bm{a}_{\bm{\mu}_k,L^2}^2  \\
    \bm{a}_{\bm{\mu}_k,L^2}^3  \\
    \vdots  \\
    \bm{a}_{\bm{\mu}_k,L^2}^{N_T}
    \end{bmatrix} ~\in~\mathbb{R}^{(N_T-1)\times N_u},
\end{equation}
then the needed time derivative vectors for the RBF interpolation are simply computed by the backward differentiation scheme as follows:
\begin{equation}\label{eq:vec_coeff}
\bm{\dot{a}}_{\bm{\mu}_k,L^2}^r = \frac{\bm{a}_{\bm{\mu}_k,L^2}^r-\bm{a}_{\bm{\mu}_k,L^2}^{r-1}}{\Delta t_{\bm{\mu}_k}} , \quad \text{for} \quad r=2,3,...,N_T, \quad k=1,2,...,M.
\end{equation}
In the formula above $\Delta t_{\bm{\mu}_k}$ is the time step at which snapshots were acquired for the parameter sample ${\bm{\mu}_k}$. As a result, the following matrix of time derivative velocity coefficients is formed
\begin{equation}
\bm{{\mathcal{\dot{A}}_k}} = \frac{\bm{{\mathcal{A}_{2,k}}}-\bm{{\mathcal{A}_{1,k}}}}{\Delta t_{\bm{\mu}_k}} = \begin{bmatrix} 
    \bm{\dot{a}}_{\bm{\mu}_k,L^2}^2  \\
    \bm{\dot{a}}_{\bm{\mu}_k,L^2}^3  \\
    \vdots  \\
    \bm{\dot{a}}_{\bm{\mu}_k,L^2}^{N_T} 
    \end{bmatrix} ~\in~\mathbb{R}^{(N_T-1)\times N_u}.
\end{equation}
Finally, merging the $L^2$ projection coefficients of velocity starting from the second time snapshot with the time derivative coefficients will yield the following matrix
\begin{equation}
\tilde{\bm{A}}_k = \begin{bmatrix} 
    \bm{{\mathcal{A}_{2,k}}} & \bm{{\mathcal{\dot{A}}_k}} \\
    \end{bmatrix} ~\in~\mathbb{R}^{(N_T-1)\times 2 N_u}.
\end{equation}
On the other hand, the projection coefficients of the eddy viscosity modes onto the snapshots are given by:
\begin{equation}\label{eq:visc_coeff_L2}
g^r_{\bm{\mu}_k,i,L^2} = (\bm{{\mathcal{S}^r_{\bm{\mu}_k,{\nu_t}}}},\eta_i)_{L^2(\Omega)}, \quad \text{for} \quad r=2,3,...,N_T \quad \text{,} \quad i =1,2,...,N_{\nu_t} \quad \text{and} \quad k=1,2,...,M.
\end{equation}
One may define the vector $\tilde{\bm{G}}_{i,k} ~\in~\mathbb{R}^{(N_T-1)}$ as the vector containing the coefficients in \eqref{eq:visc_coeff_L2} for a fixed $i$ and $k$. The combined matrices and vectors for all parameter samples will be called $\tilde{\bm{A}}$ and $\tilde{\bm{G}}_i$, respectively, which are defined as follows:
\begin{equation}
\tilde{\bm{A}} = \begin{bmatrix} 
    \tilde{\bm{A}}_1  \\
    \tilde{\bm{A}}_2  \\
    \vdots  \\
    \tilde{\bm{A}}_M 
    \end{bmatrix} ~\in~\mathbb{R}^{(N_s-M)\times 2N_u},
    \tilde{\bm{G}}_i = \begin{bmatrix} 
    \tilde{\bm{G}}_{i,1}  \\
    \tilde{\bm{G}}_{i,2}  \\
    \vdots  \\
    \tilde{\bm{G}}_{i,M}
    \end{bmatrix} ~\in~\mathbb{R}^{(N_s-M)},
\end{equation}
 At this point, the goal of the interpolation will be to approximate the maps $[f_i]_{i=1}^{N_{\nu_t}}$ in $g_i = f_i(\bm{a},\bm{\dot{a}})$, where:
\begin{equation}
f_i : \mathbb{R}^{2 N_u} \rightarrow \mathbb{R}.
\end{equation}
This approximation is based on the interpolation points given in each row of the  matrix $\tilde{\bm{A}}$ and the vector $\tilde{\bm{G}}_i$.
\subsection{Treatment of boundary conditions}\label{sec:bc}

In reduced order modeling, it is often the case that the parameterization is in the boundary conditions, and in particular at the inlet boundary. In this subsection, the available methodologies to tackle this aspect are presented. The main two methods to take into account boundary conditions at reduced order level are the penalty method \cite{Bizon2012,Babuka1973,Barrett1986,Kalashnikova2012,Sirisup2005} and the lifting function method \cite{Graham1999,Gunzburger2007,QUIET}.\par
Let $\Gamma_D$ be the Dirichlet boundary that might be composed by separate boundaries, i.e. $\Gamma_D = {\Gamma_D}_1 \cup {\Gamma_D}_2 ... \cup {\Gamma_D}_K$. Let $N_{BC}$ be the number of velocity boundary conditions we would like to impose on some parts of the Dirichlet boundary. We emphasize that, each non-zero scalar component value of the velocity field that has to be set at one part of the boundary, is counted as one boundary condition. As an example let $\bm{U_{Dir}} = (U_x,U_y)$ be the velocity vector that must be imposed at the Dirichlet boundary for the problem under interest. It is supposed that $U_x$ and $U_y$ are the values of the velocity components in the $x$ and $y$ directions, respectively, in this case there are two boundary conditions to set and thus $N_{BC} = 2$. Let $U_{BC,i,j}$ be the value of $i$-th component of the velocity to be imposed at reduced order level at the $j$-th part of the Dirichlet boundary ${\Gamma_D}_j$. We define $\bm{U_{BC}}$ as the vector of all scalar velocities $U_{BC,i,j}$, this vector has a dimension of $N_{BC}$, and ${U_{BC}}_k$ is the $k$-th element of $\bm{U_{BC}}$. 
\subsubsection{The penalty method}\label{sec:penalty}
In the penalty method, an additional term is added in the formulation of the dynamical system of the reduced order model. The added term represents a constraint that has to be satisfied at the reduced order level on certain parts of the boundary. The penalty method has been used for both laminar and turbulent reduced order models as presented in \cite{Lorenzi2016}. If we consider employing the method addressed in \cite{Lorenzi2016} to the POD-Galerkin Mixed-ROM model, the result will be the following system:
\begin{equation}\label{eq:PODsys_pen}
\begin{cases}
\bm{M} \bm{\dot{a}}= \nu (\bm{B}+\bm{B_T})\bm{a}-\bm{a^T}\bm{C}\bm{a} +\bm{g^T}(\bm{C_{T1}}+\bm{C_{T2}})\bm{a}-\bm{H}\bm{b} + \tau(\sum_{k=1}^{N_{BC}} ({U_{BC}}_k \bm{D^k} - \bm{E^k} \bm{a})),\\
\bm{P} \bm{a}= \bm{0},\\ 
\end{cases}
\end{equation}
where $\tau$ is called the penalization factor, and its value is usually determined by sensitivity analysis. In general the higher the value of $\tau$ is and the stronger is the enforcement of the boundary conditions. The additional boundary terms with respect to system \eqref{eq:PODsys} are defined as follows:
\begin{align}
& (\bm{D^k})_{i} = \left(\bm{\phi_i} \right)_{L^2({\Gamma_D}_k)},\\
& (\bm{E^k})_{ij} = \left(\bm{\phi_i}, \bm{\phi_j} \right)_{L^2({\Gamma_D}_k)}.
\end{align}
In this method, the POD is applied directly on the snapshots matrices for all variables without the homogenization of the fields. This will result in POD modes which don't have homogeneous Dirichlet boundary conditions.

\subsubsection{The lifting function method}\label{sec:lifting}

The lifting or control function method involves the use of the so-called lifting function which handles the non-homogeneous values on the boundaries. The method involves the creation of a new set of snapshots for the velocity field where the non-homogeneous Dirichlet boundary conditions are removed. After that the POD procedure is applied on the newly formed snapshots and this gives POD modes which have homogeneous Dirichlet conditions at the Dirichlet boundary.\par
The procedure of modifying the velocity snapshots is done as follows:
\begin{equation}\label{eq:lifting}
\bm{\tilde{u}_k} = \bm{u_k} - \bm{{\phi_L}} \cdot \bm{U_{BC}},
\end{equation}
where $\bm{{\phi_L}} ~\in~\mathbb{R}^{N_u^h\times N_{BC}}$ is a matrix of the lifting functions $\bm{{\phi_L}_{i,j}}$. Each lifting function $\bm{{\phi_L}_{i,j}}$ has homogeneous Dirichlet boundary conditions in all parts of the Dirichlet boundary except in the $i$-th component at ${\Gamma_D}_j$ where it has unitary value. 
We would like to remark that the same lifting method can be used for pressure fields if the formulation of the problem involves non-homogeneous Dirichlet boundary condition for pressure. In that case if the non-homogeneous pressure value is $p_{out}$ then the new pressure snapshots will be computed as follows:
\begin{equation}\label{eq:lifting}
\tilde{p}_k = p_k - p_{out} {\chi_c},
\end{equation}
where ${\chi_c}$ is the pressure lifting function. The new snapshots matrices for velocity and pressure are denoted, respectively, by $\bm{\mathcal{\tilde{U}}} = [\bm{\tilde{u}_1},\bm{\tilde{u}_2},\dots,\bm{\tilde{u}_{N_s}}]$ and $\bm{\mathcal{\tilde{P}}} = [\tilde{p}_1,\tilde{p}_2,\dots,\tilde{p}_{N_s}]$. These snapshots matrices will be used for computing the reduced order bases for the velocity and pressure POD spaces, respectively.\par
During the online stage, it is required to approximate the velocity and pressure fields for the value of the combined time-parameter vector $\bm{z}^*$ (which might contain new velocity values to be imposed at some parts of the boundary), this could be done as follows:
\begin{equation}
\bm{u}(\bm{x},\bm{z}^*)\approx  \bm{{\phi_L}} \cdot \bm{U_{BC}}^* + \sum_{i=1}^{N_u} a_i(\bm{z}^*) \bm{\phi_i} (\bm{x}), \quad p(\bm{x},\bm{z}^*)\approx  p_{out} {\chi_c} + \sum_{i=1}^{N_p} b_i(\bm{z}^*) {\chi_i} (\bm{x}),
\end{equation}
where $\bm{U_{BC}}^*$ is the vector of boundary velocity values that corresponds to $\bm{z}^*$.
\section{Numerical results}\label{sec:results}
In this section, we present the results obtained applying the proposed POD-Galerkin Mixed-ROM on two turbulent flow problems. The first problem is that of the turbulent flow past a backstep. Such classical benchmark in the turbulence modeling community is here considered in a steady state parametrized setup. The second problem analyzed is that of the turbulent flow past a circular cylinder. In the second problem the case is parametrized with Reynolds number being the parameter. In both cases we will present a comparison of the Mixed-ROM results with the ones obtained by the ROM developed in \cite{Lorenzi2016}. The authors in their work proposed a ROM which considers a single set of reduced coefficients for the velocity, pressure and eddy viscosity field. In such a way it is possible to resort on a unified approach to deal with turbulent flows. Such an approach, that is used as comparison with respect to methodology developed here, will be referred as P-ROM from now on.\par
The finite volume C\texttt{++} library OpenFOAM\textsuperscript{\textregistered} (OF) \cite{weller1998tensorial} is used as the numerical solver at the full order level. At the reduced order level the reduction and resolution of the reduced order system is carried out using the C\texttt{++} based library ITHACA-FV \cite{RoSta17}.\par
\subsection{Steady case}\label{sec:numerical_steady}
\begin{figure}
\centering
\includegraphics[width=0.7\textwidth]{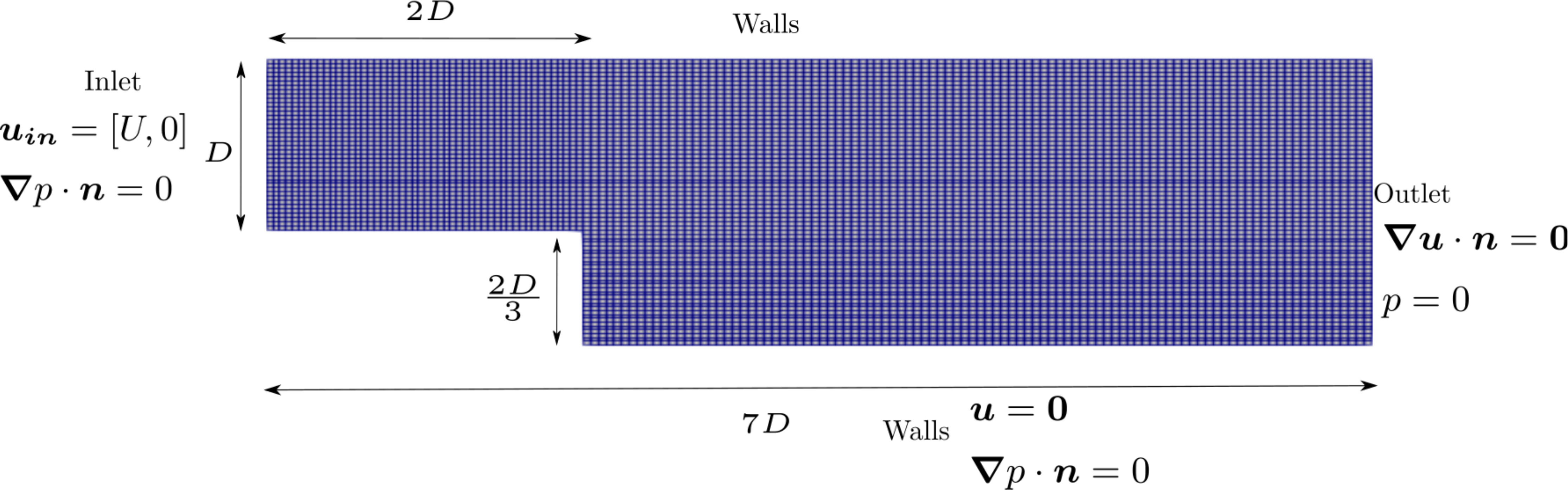}
\caption{The computational domain used in the numerical simulations, all lengths are described in terms of the characteristic length $D$ that is equal to $1$ meter.}\label{fig:comp_domain}
\end{figure}

The steady Mixed-ROM solver has been tested on the backward step benchmark case. \autoref{fig:comp_domain} depicts the layout of the domain with details of the computational mesh. The plot also reports the boundary conditions enforced on every side of the domain. The inflow velocity $U$ has been modified in the simulations to parametrize the problem with respect to the Reynolds number. Thus, the objective of this numerical experiment is to assess the Mixed-ROM solver ability to reproduce flows with high Reynolds number and their dependence on the parameters. To this end, the Mixed-ROM results will be compared both to the full order results and to the results of the P-ROM model. Moreover, to test the Mixed-ROM solver capability to deal with different turbulence models, we tested the model both on full order solutions obtained with $k-\epsilon$ and SST $k-\omega$ models.\par
The $100$ snapshots required for the training during the offline phase were generated by solving the FOM with inlet velocity values ranging from $1$ \si{m\per s} to $25$ \si{m\per s} on an equally spaced distribution. Given the physical viscosity $\nu = 10^{-3}$ \si{m^2\per s} and the characteristic length is $D =1$ \si{m}, this corresponds to a Reynolds number that varies from $1 \times 10^3$ to $2.5 \times 10^4$. In the full order simulations, Gauss linear scheme was selected for the approximation of the gradients and Gauss linear scheme with non-orthogonal correction was selected to approximate the Laplacian terms. A $2$-nd order bounded Gauss upwind scheme was instead used for the approximation of the convective term. Finally, $1$st order bounded Gauss upwind scheme is used to approximate all terms involving the turbulence model parameters $k$, $\epsilon$ and $\omega$.\par

\begin{figure}
\centering
\includegraphics[width=0.7\textwidth]{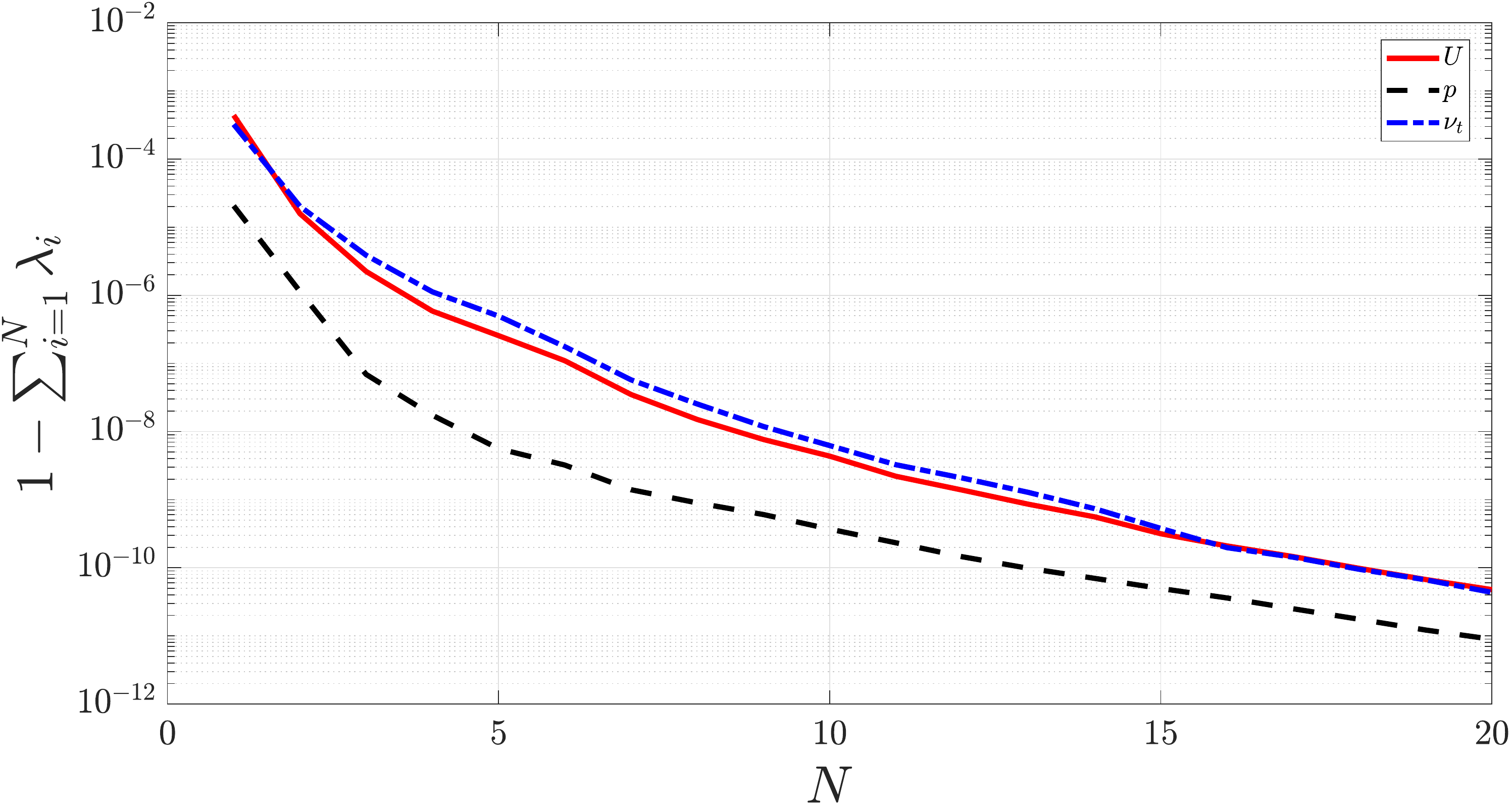}
\caption{Cumulative ignored eigenvalues decay. In the plot, the solid red line refers to the velocity eigenvalues, the dashed black line indicates the pressure eigenvalues and the dash-dotted blue line finally refers to the eddy viscosity eigenvalues.}\label{fig:Cum_backstep}
\end{figure}

The modes of velocity, pressure and eddy viscosity fields have been obtained by POD analysis of the snapshots matrices. \autoref{fig:Cum_backstep} shows the cumulative eigenvalues decay for velocity, pressure and eddy viscosity. As can appreciated relatively small number of modes is sufficient to recover most of the energetic information in the snapshots.\par

Once the reduced model training was carried out and the modes were computed, Mixed-ROM and P-ROM simulations have been carried out on new set of sampling points in the parameter space. More specifically, the online sample values for $U$, denoted by $U_i^*$ where $i=1,...,N_{online-samples}$,  have been chosen as $80$ equally distributed samples in the range of $[3,20]$. This set of samples includes both samples close to those used in the offline stage and samples which lie almost midway between two offline samples. Clearly, the test is aimed at assessing how accurate the reduced approximation is for parameter values that were not in the training set.\par
The enforcement of the correct inflow velocity in the reduced simulations is carried out by means of the penalty method as in \ref{sec:penalty}. In this regard, we must here remark  that the simulations results appeared quite sensitive to the penalization factor $\tau$. Thus, a sensitivity analysis had to be performed to set the value of $\tau$ for both $k-\epsilon$ and SST $k-\omega$ turbulence models considered.\par 

The first step of the online stage is represented by the interpolation of the eddy viscosity coefficients with respect to the values of the considered parameter (the inflow velocity). More specifically, the result of the interpolation is the vector $\bm{g}$, which is used to solve the reduced system \eqref{eq:PODsys} and finally obtain the vectors of coefficients $\bm{a}$ and $\bm{b}$. The interpolation using the RBF in this work has been carried out using the C\texttt{++} library SPLINTER \cite{SPLINTER}.\par 

\begin{figure}
  \centering
  \begin{minipage}[b]{0.5\linewidth}
    \centering
    \includegraphics[width=0.98\linewidth]{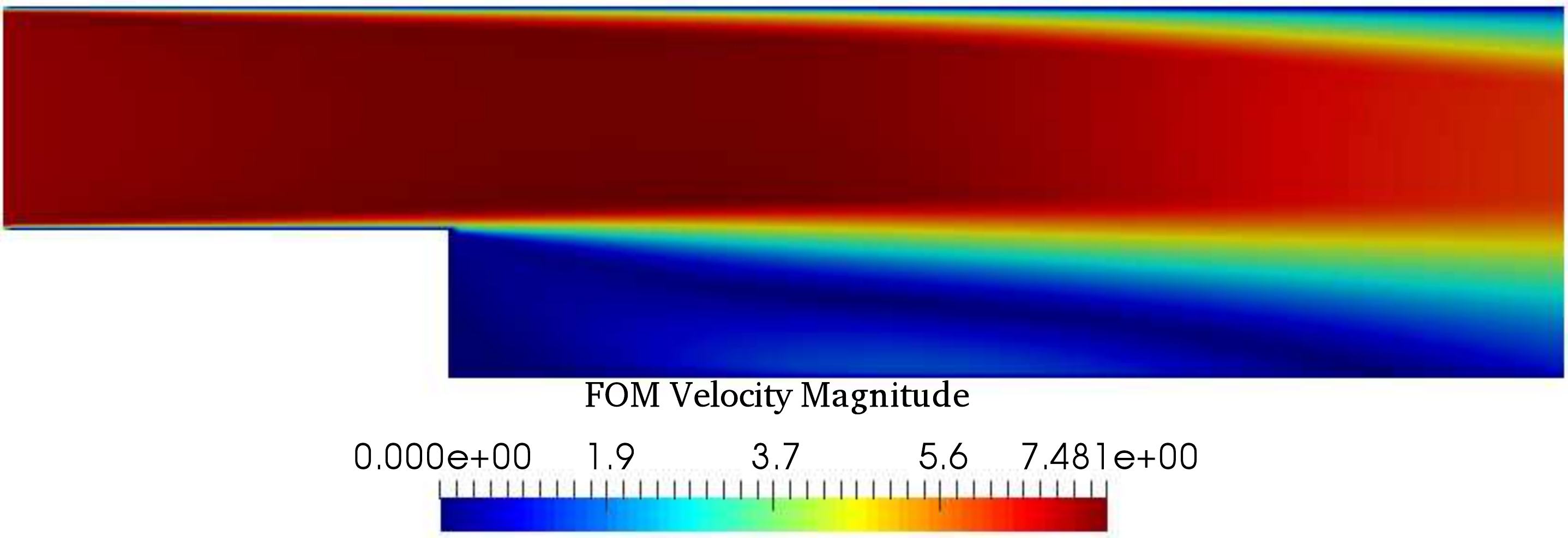}
    \scriptsize(a) 
    \end{minipage}
  \begin{minipage}[b]{0.5\linewidth}
    \centering
    \includegraphics[width=0.98\linewidth]{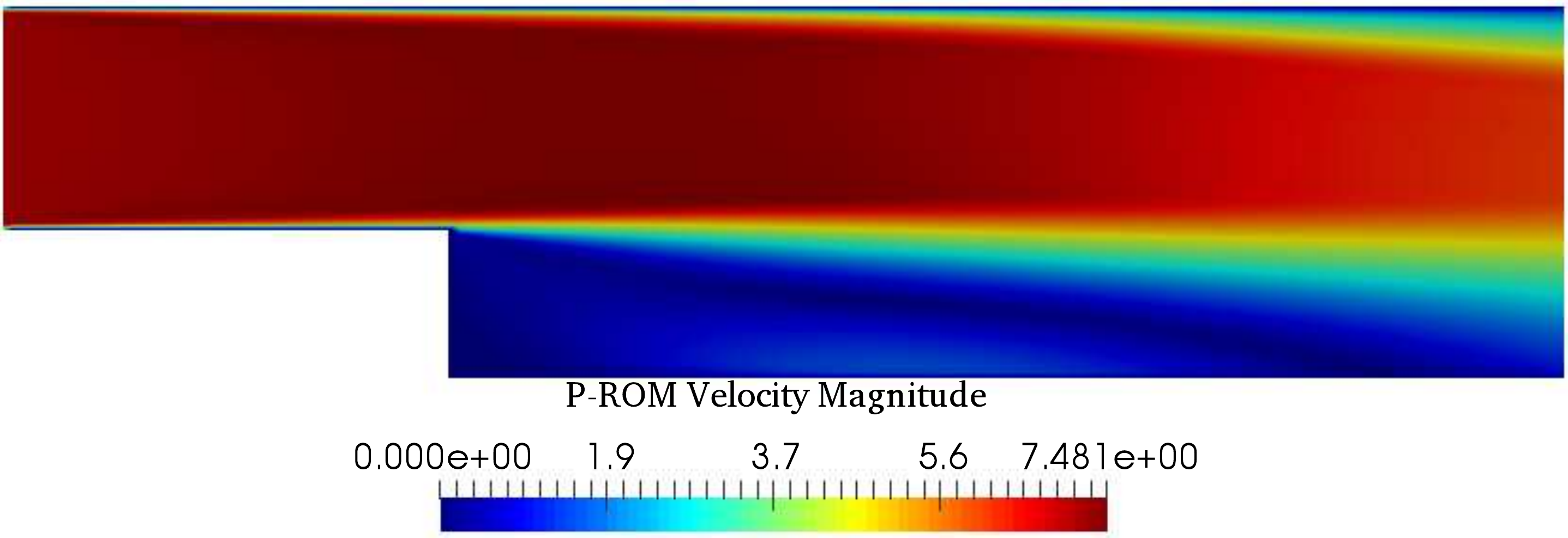}
        \scriptsize(b) 
  \end{minipage}
   \begin{minipage}[b]{0.5\linewidth}
    \centering
        \includegraphics[width=0.98\linewidth]{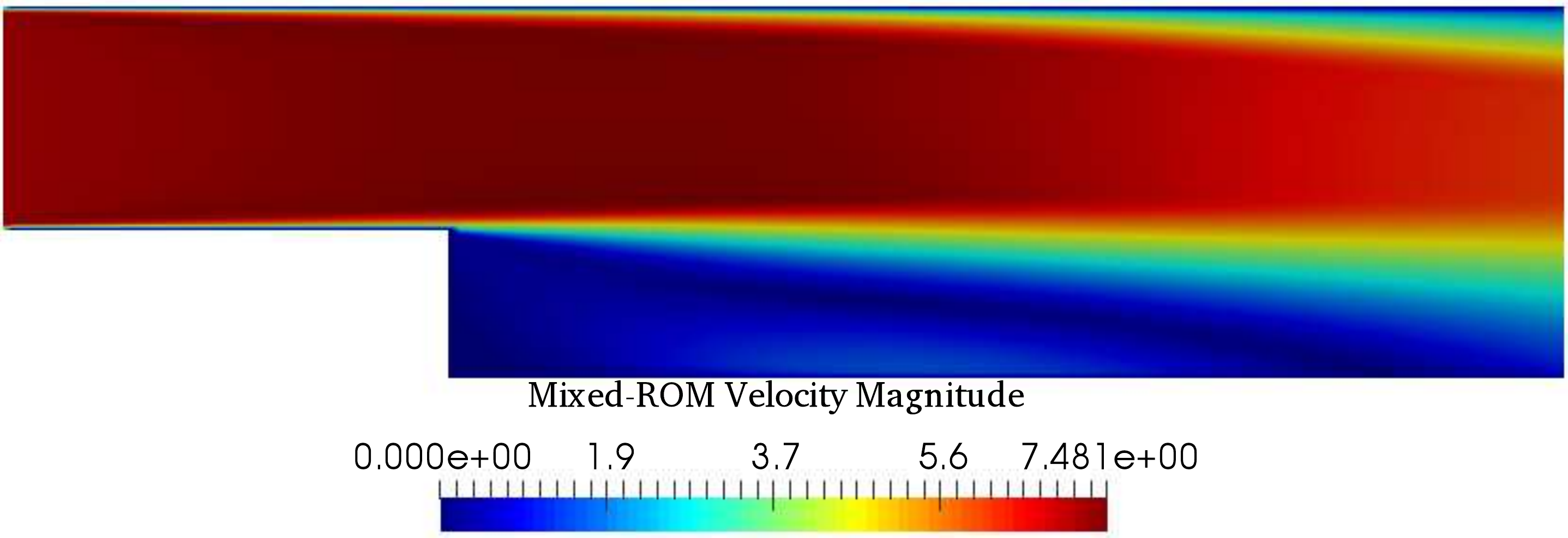} 
        \scriptsize(c) 
  \end{minipage}
  \vspace{-0.2cm}\caption{ $k-\epsilon$ turbulence model case, velocity fields for the value of the parameter $ U = 7.0886$  \si{m\per s}: (a) shows the FOM velocity, while in (b) one can see the P-ROM velocity, and finally in (c) we have the Mixed-ROM velocity.}\label{fig:u_fields}\par
\end{figure}

\begin{figure}
  \centering
 \begin{minipage}[b]{0.5\linewidth}
    \centering
    \includegraphics[width=0.98\linewidth]{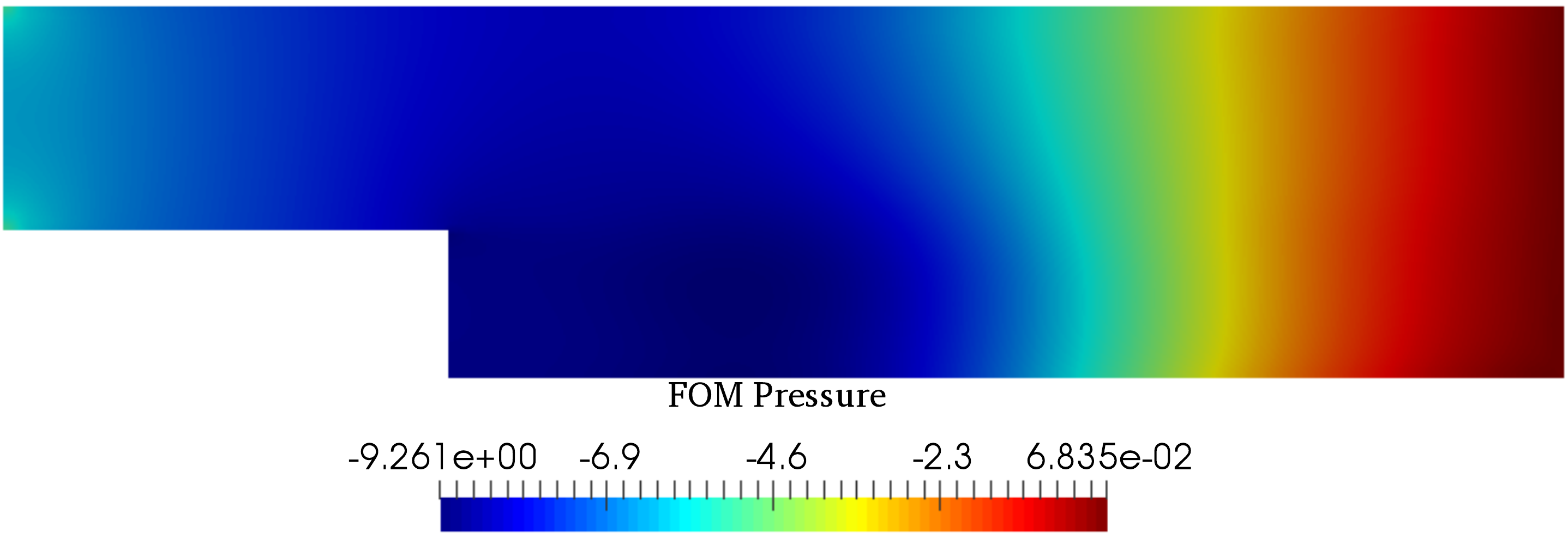}
    \scriptsize(a) 
    \end{minipage}
  \begin{minipage}[b]{0.5\linewidth}
    \centering
    \includegraphics[width=0.98\linewidth]{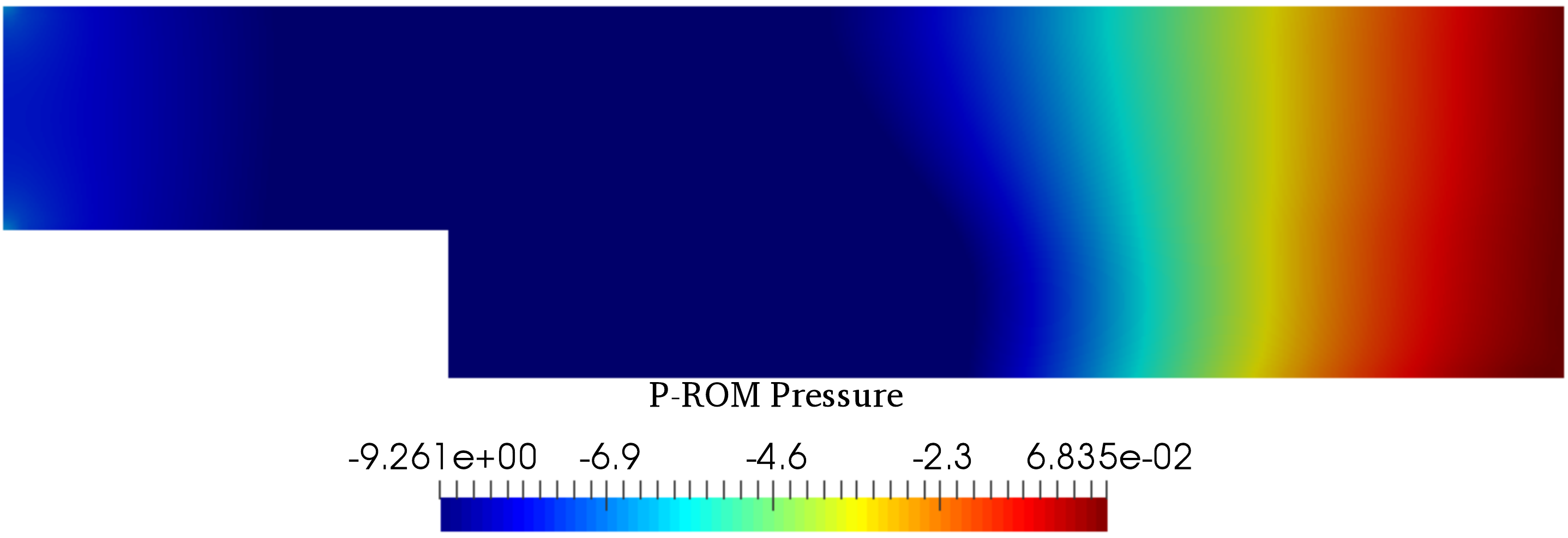}
        \scriptsize(b) 
  \end{minipage}
   \begin{minipage}[b]{0.5\linewidth}
    \centering
        \includegraphics[width=0.98\linewidth]{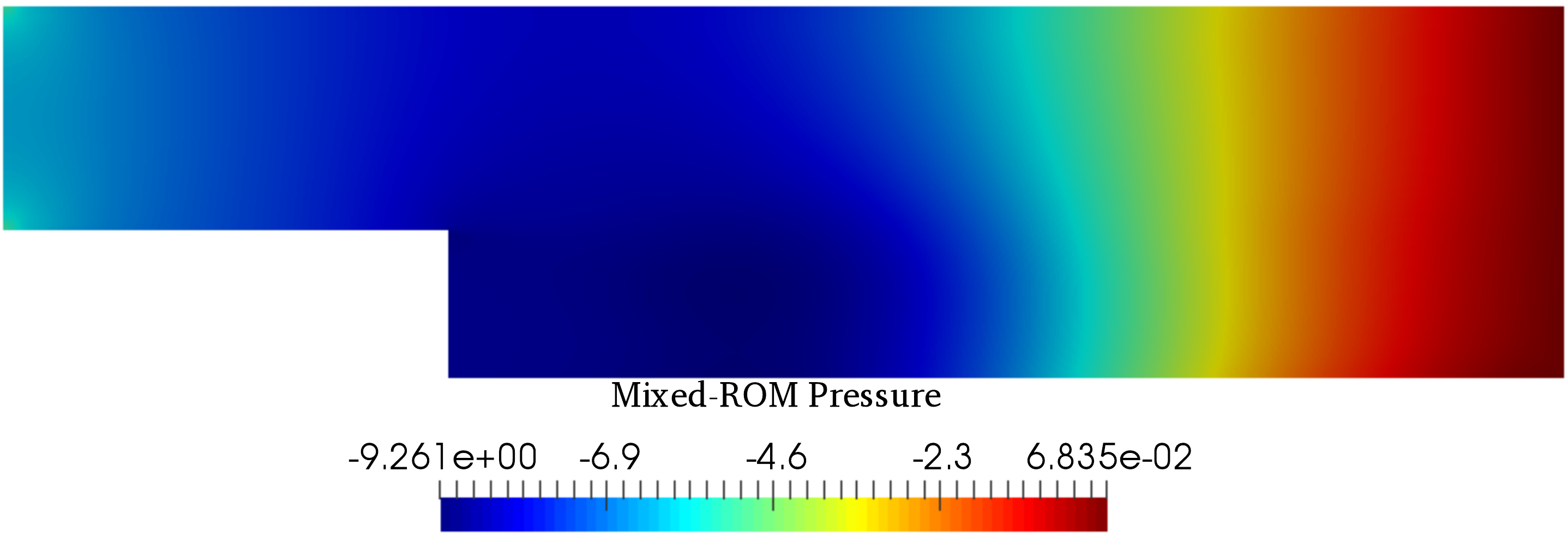} 
    \scriptsize(c)
  \end{minipage}
 \vspace{-0.2cm}\caption{$k-\epsilon$ turbulence model case, pressure fields for the value of the parameter $ U = 7.0886$  \si{m\per s}: (a) shows the FOM pressure, while in (b) one can see the P-ROM pressure, and finally in (c) we have the Mixed-ROM pressure.}\label{fig:p_fields} 
\end{figure}

\autoref{fig:u_fields} depicts the velocity fields corresponding to $ U^* = 7.0886$  \si{m\per s} computed via the FOM, the P-ROM and Mixed-ROM in the case of $k-\epsilon$ turbulence model. A similar comparison is presented in \autoref{fig:p_fields} for the pressure fields. We remark that all the solutions were generated using $10$ velocity, pressure, supremizer and eddy viscosity modes in the online stage for both the Mixed-ROM and the P-ROM. The images clearly indicate that the hybrid projection/data-driven-based approach allows for qualitatively accurate approximations of the FOM solutions. This is clearly not the case when the P-ROM approach is employed since the pressure field does not correctly reproduce its FOM counterpart. To provide a quantitative measurement of both reduced order models performance, we evaluate the relative $L^2$ error for velocity and pressure which, respectively, read
\begin{equation}\label{eq:l2_error}
\epsilon_u = \frac{{\left\lVert \bm{u} - \bm{u}^* \right\rVert}_{L^2(\Omega)}}{{\left\lVert \bm{u}\right\rVert}_{L^2(\Omega)}} \times 100 \%,\\
\epsilon_p = \frac{{\left\lVert p - p^* \right\rVert}_{L^2(\Omega)}}{{\left\lVert p \right\rVert}_{L^2(\Omega)}} \times 100 \%,
\end{equation}
in which $\bm{u}^*$ and $p^*$ are general reduced order velocity and pressure fields, respectively. The relative $L^2$ errors between the
FOM and the Mixed-ROM velocity and pressure fields presented in \autoref{fig:u_fields} and \autoref{fig:p_fields} are respectively $\epsilon_u=0.4444$ $\%$ and  $\epsilon_p=0.3654$ $\%$. As for the P-ROM results, the corresponding errors are $\epsilon_u=0.6522$ $\%$ and $\epsilon_p=20.9441$ $\%$, respectively.\par

\begin{figure}
  \centering
 \begin{minipage}[b]{0.5\linewidth}
    \centering
    \includegraphics[width=0.98\linewidth]{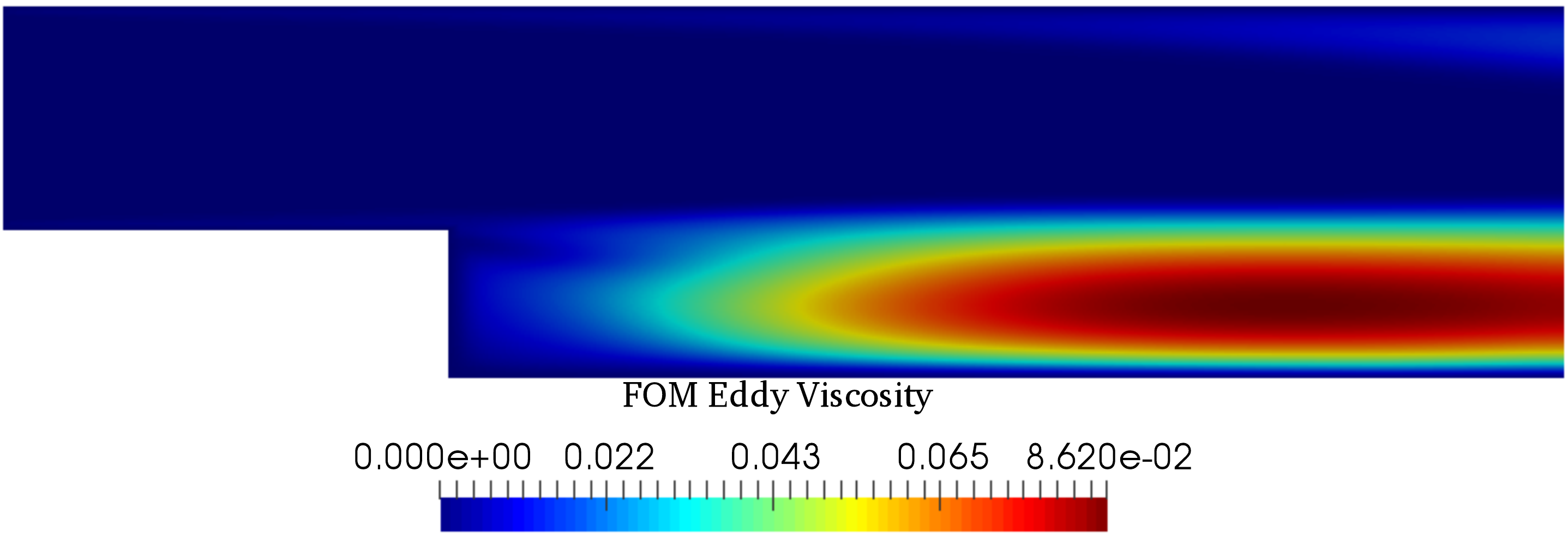} 
    \scriptsize(a) 
    \end{minipage}
   \begin{minipage}[b]{0.5\linewidth}
    \centering
        \includegraphics[width=0.98\linewidth]{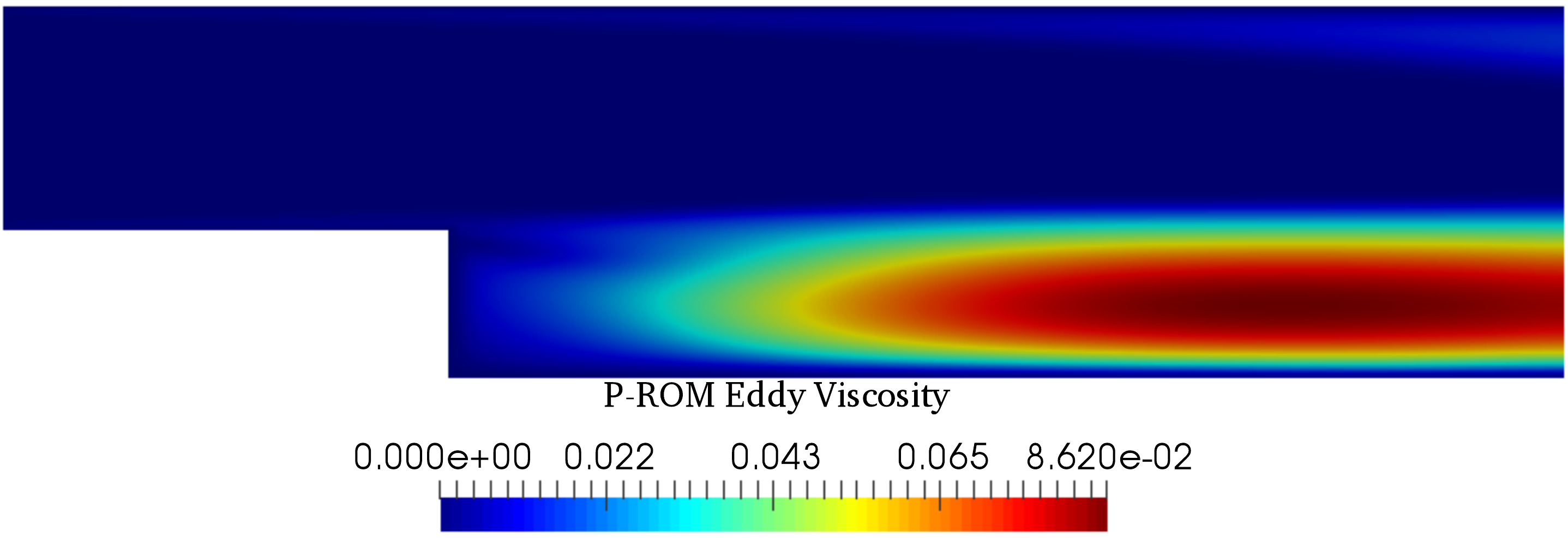} 
    \scriptsize(b)
  \end{minipage}
  \begin{minipage}[b]{0.5\linewidth}
    \centering
        \includegraphics[width=0.98\linewidth]{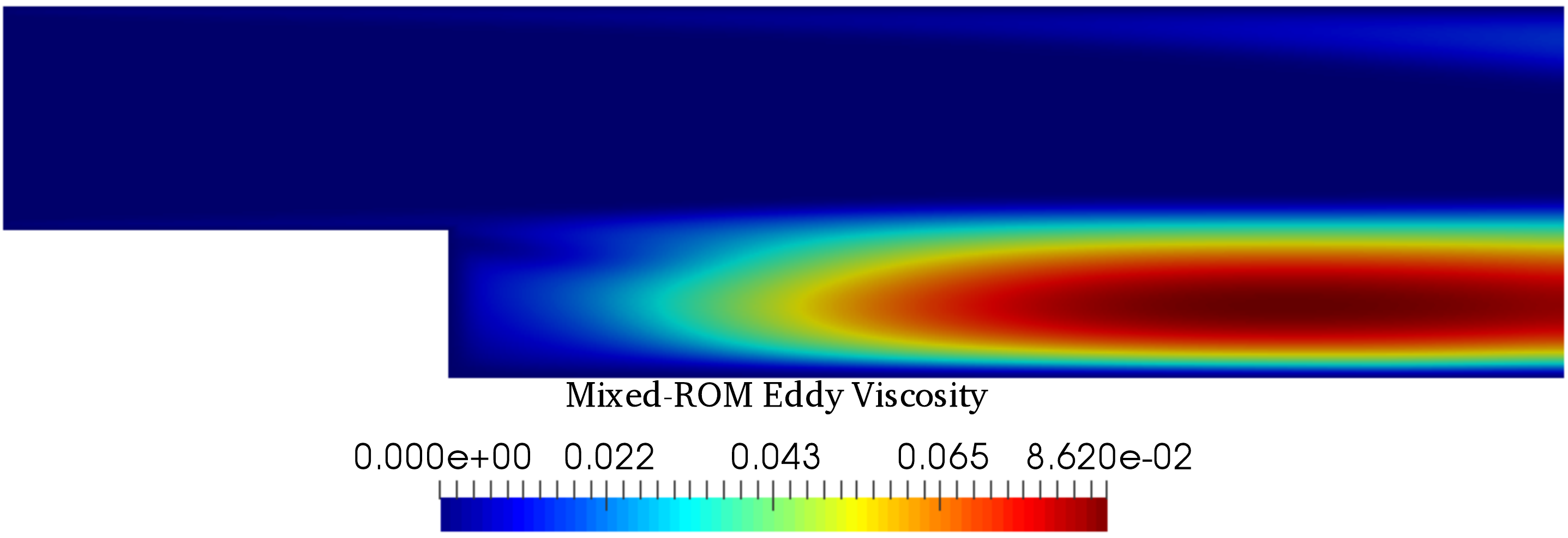} 
    \scriptsize(c)
  \end{minipage}
 \vspace{-0.2cm}\caption{$k-\epsilon$ turbulence model case, eddy viscosity fields: (a) shows the FOM eddy viscosity, while in (b) one can see the P-ROM eddy viscosity, and finally in (c) we have the Mixed-ROM eddy viscosity.}\label{fig:nut_fields} 
\end{figure}
\begin{figure}
  \centering
  \begin{minipage}[b]{0.5\linewidth}
    \centering
    \includegraphics[width=0.98\linewidth]{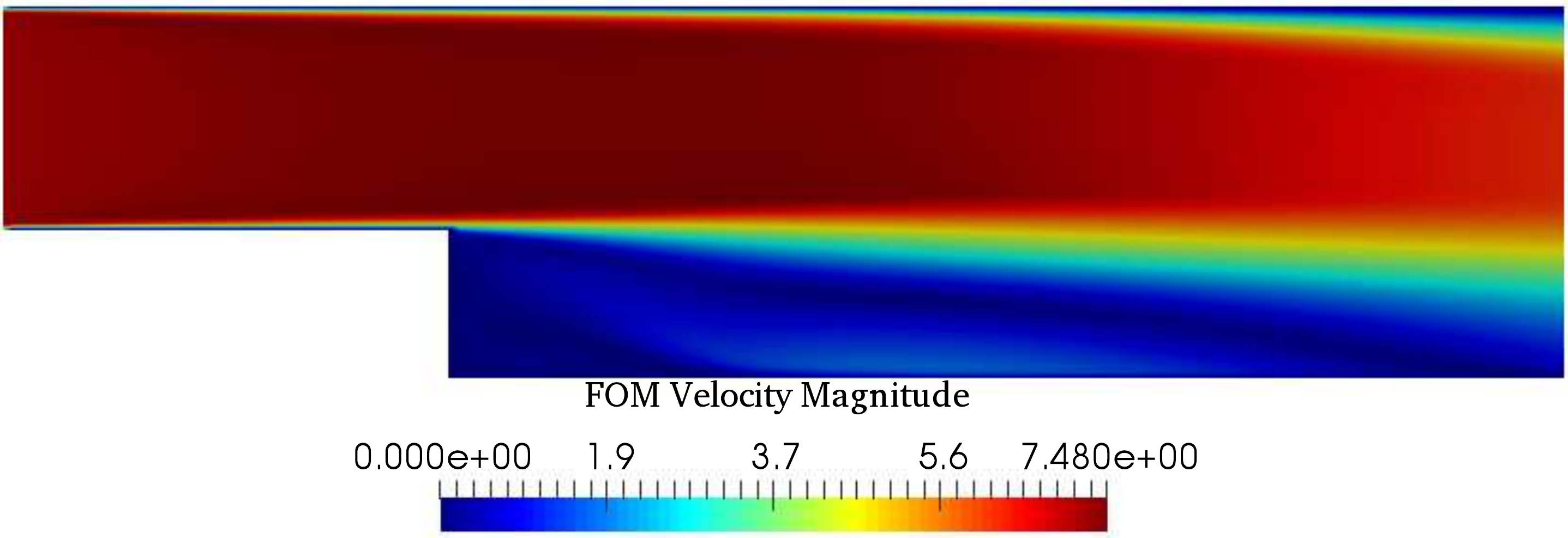}
    \scriptsize(a) 
    \end{minipage}
  \begin{minipage}[b]{0.5\linewidth}
    \centering
    \includegraphics[width=0.98\linewidth]{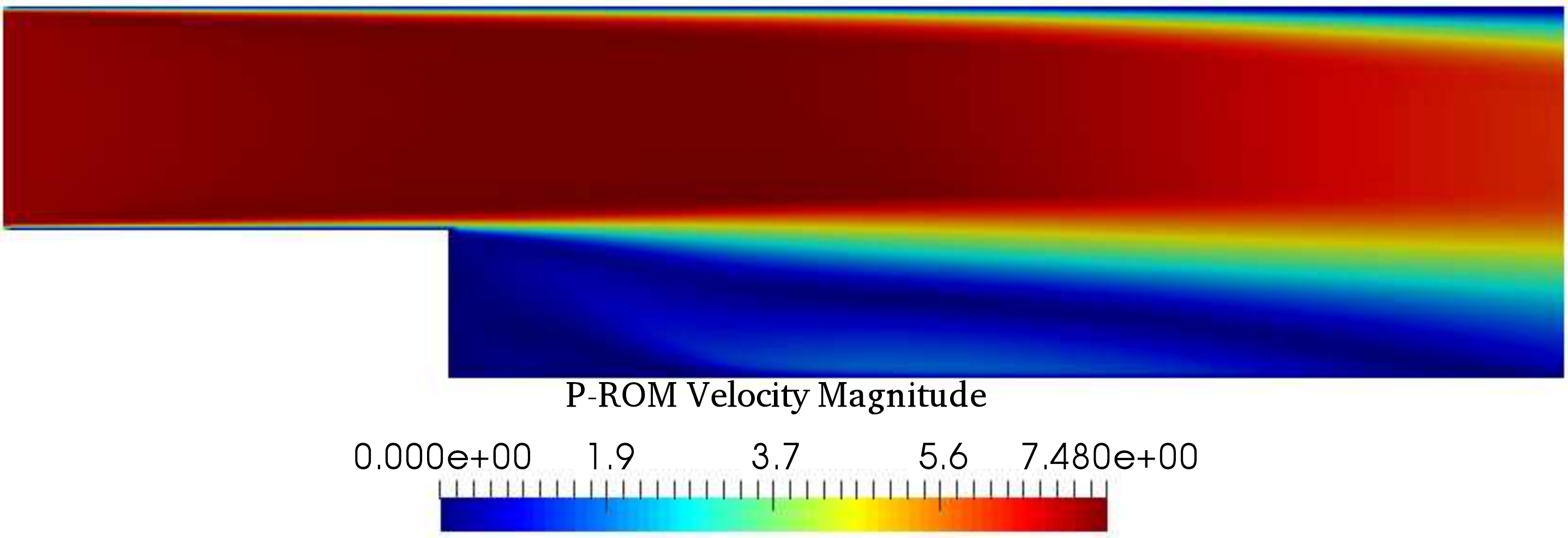}
        \scriptsize(b) 
  \end{minipage}
   \begin{minipage}[b]{0.5\linewidth}
    \centering
        \includegraphics[width=0.98\linewidth]{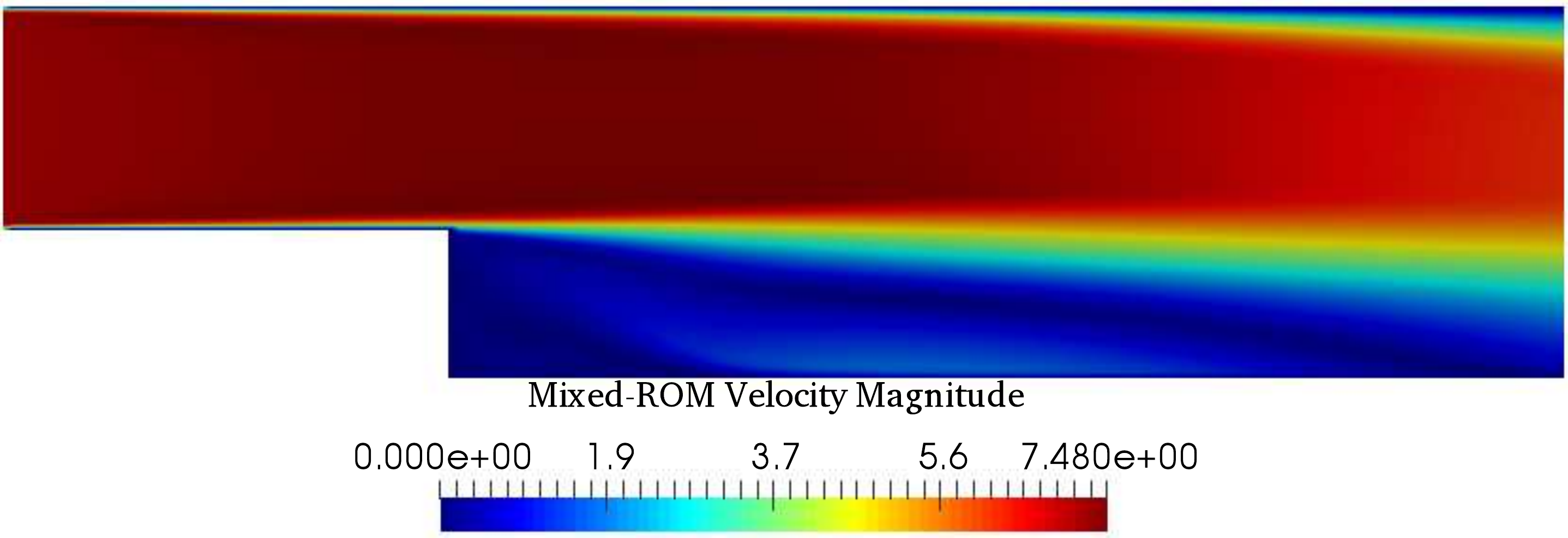} 
        \scriptsize(c) 
  \end{minipage}
  \vspace{-0.2cm}\caption{SST $k-\omega$ turbulence model case, velocity fields for the value of the parameter $ U = 7.0886$  \si{m\per s}: (a) shows the FOM velocity, while in (b) one can see the P-ROM velocity, and finally in (c) we have the Mixed-ROM velocity.}\label{fig:u_fieldsKomega}\par
\end{figure}
\begin{figure}
  \centering
 \begin{minipage}[b]{0.5\linewidth}
    \centering
    \includegraphics[width=0.98\linewidth]{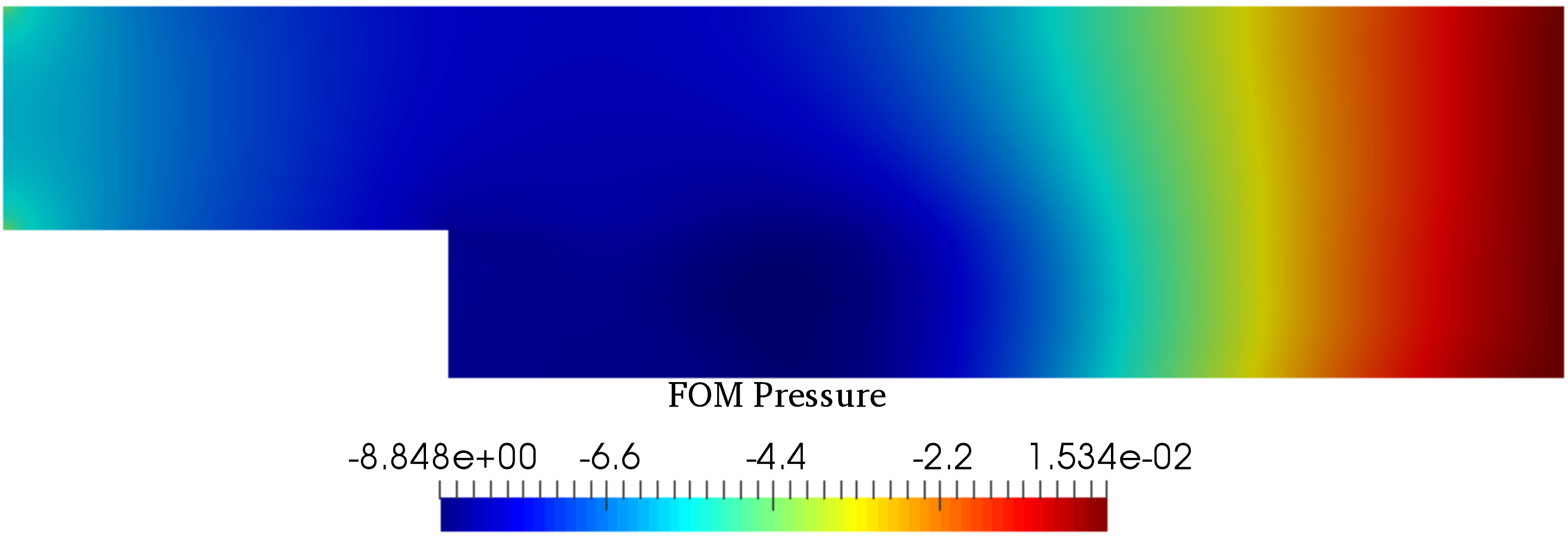}
    \scriptsize(a) 
    \end{minipage}
  \begin{minipage}[b]{0.5\linewidth}
    \centering
    \includegraphics[width=0.98\linewidth]{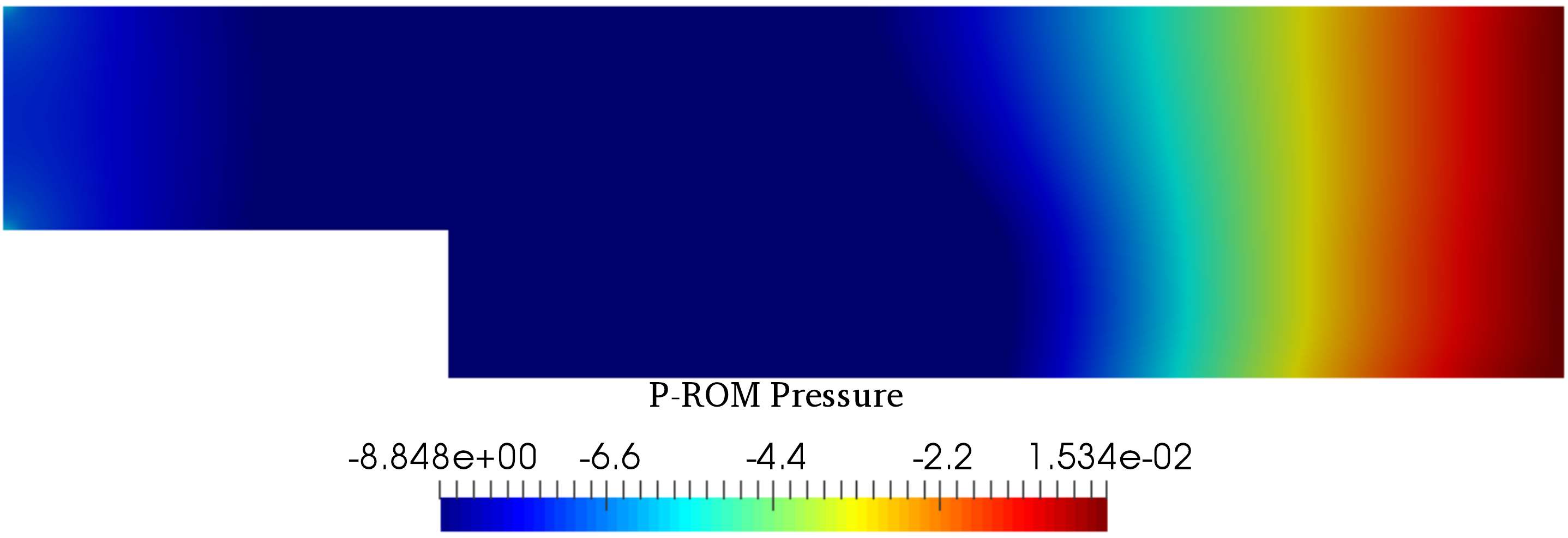}
        \scriptsize(b) 
  \end{minipage}
   \begin{minipage}[b]{0.5\linewidth}
    \centering
        \includegraphics[width=0.98\linewidth]{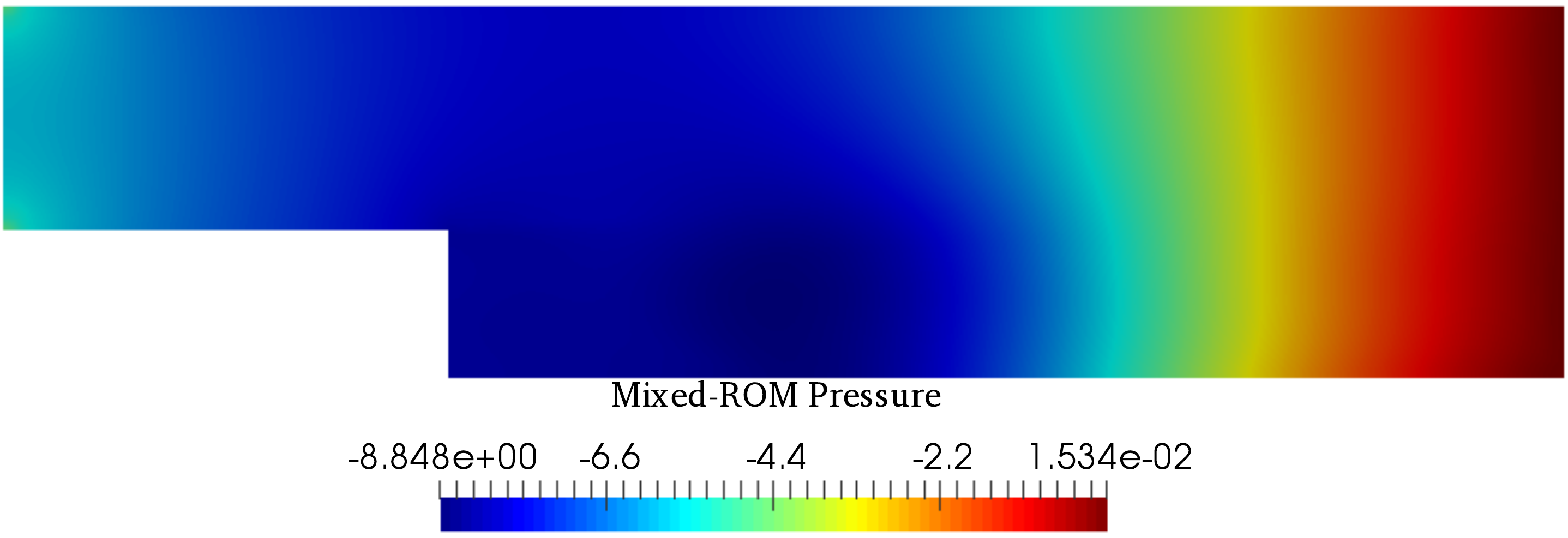} 
    \scriptsize(c)
  \end{minipage}
 \vspace{-0.2cm}\caption{ SST $k-\omega$ turbulence model case, pressure fields for the value of the parameter $ U = 7.0886$  \si{m\per s}: (a) shows the FOM pressure, while in (b) one can see the P-ROM pressure, and finally in (c) we have the Mixed-ROM pressure.}\label{fig:p_fieldsKomega} 
\end{figure}
\begin{figure}
  \centering
 \begin{minipage}[b]{0.5\linewidth}
    \centering
    \includegraphics[width=0.98\linewidth]{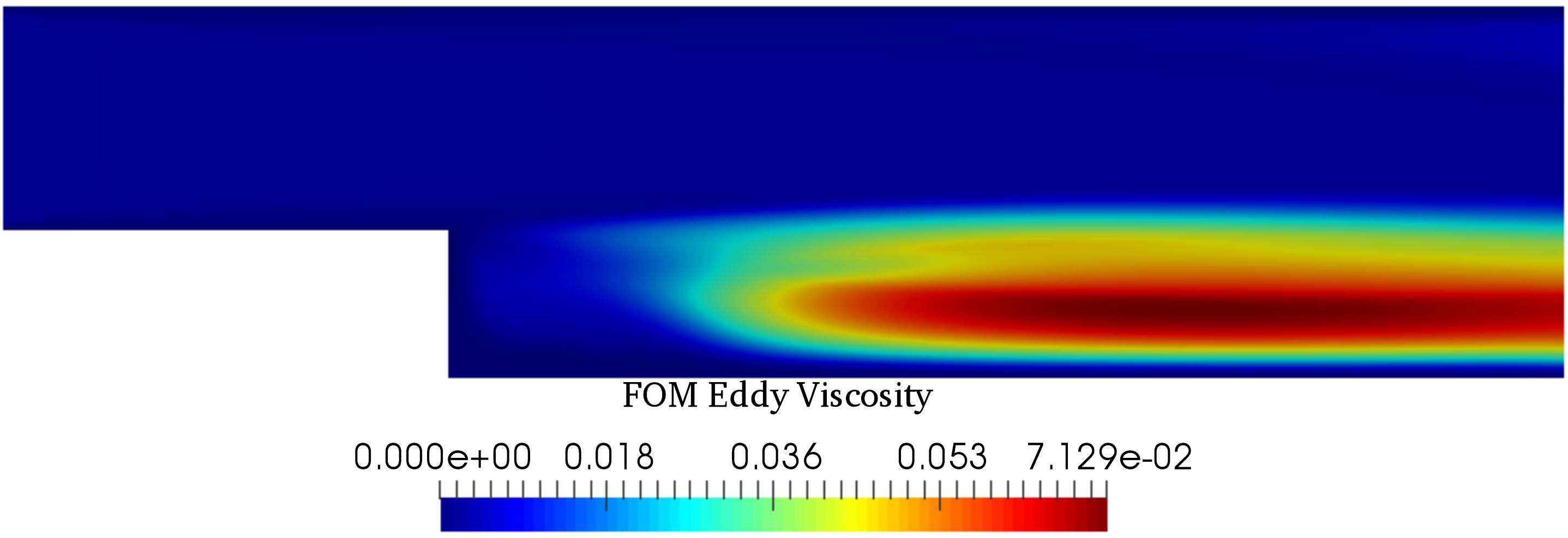}
    \scriptsize(a) 
    \end{minipage}
  \begin{minipage}[b]{0.5\linewidth}
    \centering
    \includegraphics[width=0.98\linewidth]{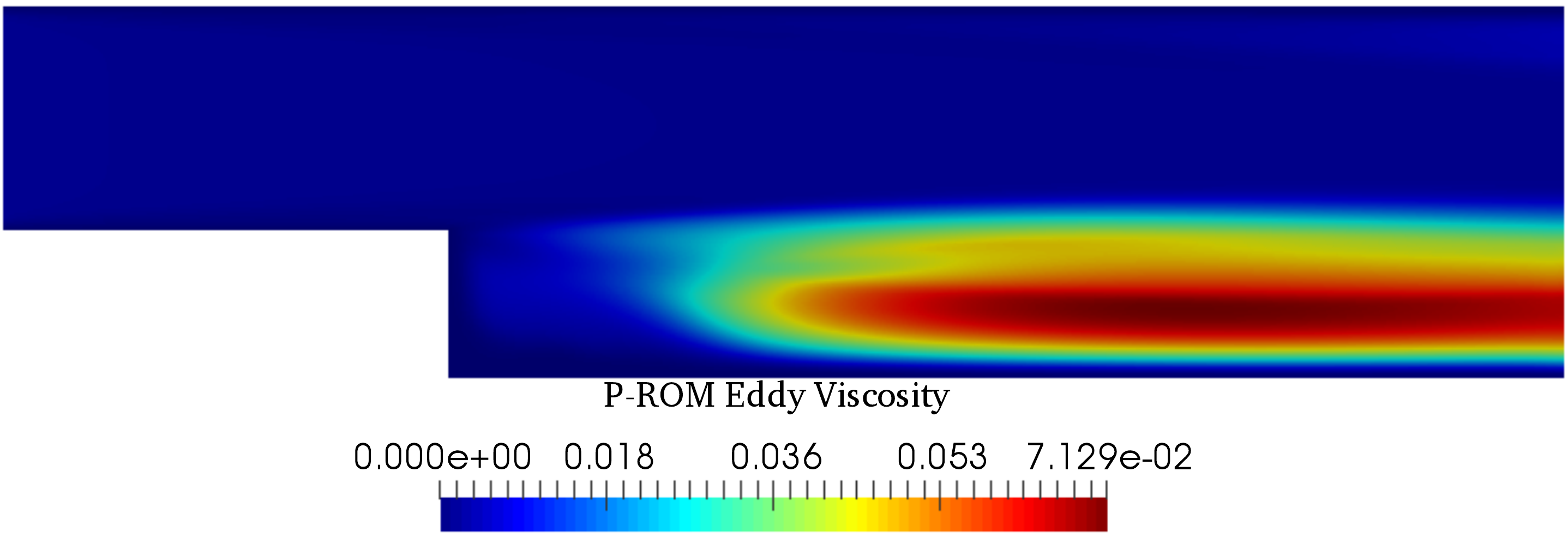}
        \scriptsize(b) 
  \end{minipage}
   \begin{minipage}[b]{0.5\linewidth}
    \centering
        \includegraphics[width=0.98\linewidth]{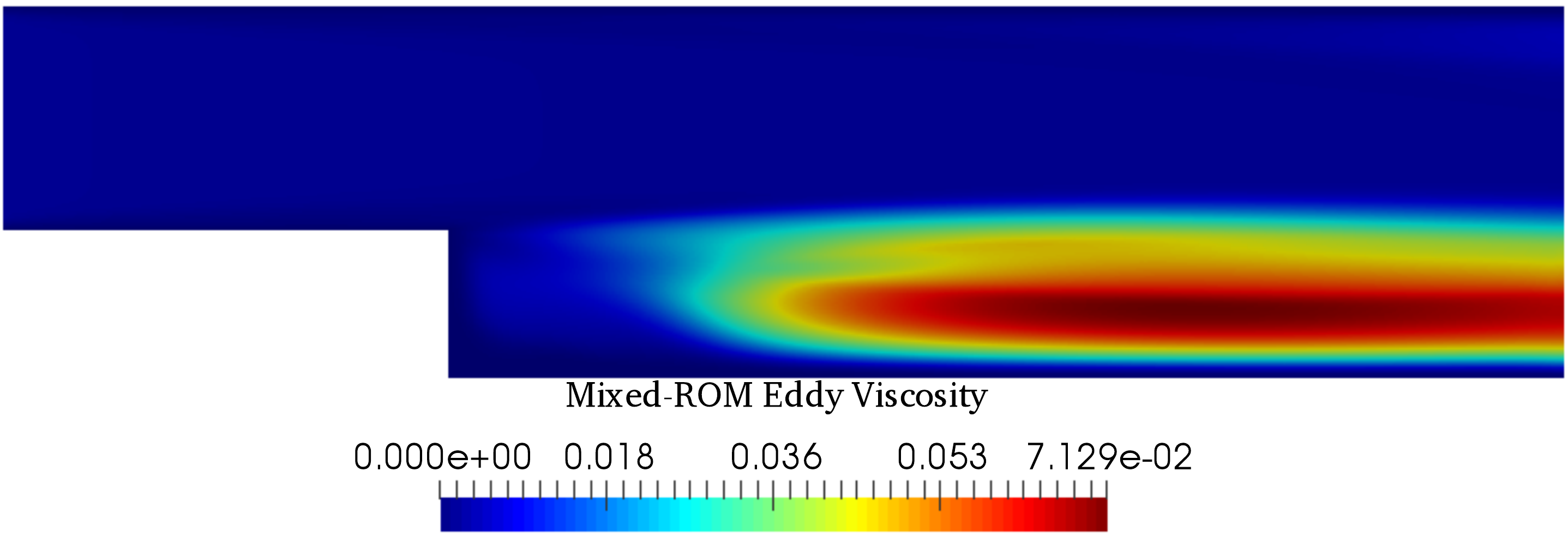} 
    \scriptsize(c)
  \end{minipage}
 \vspace{-0.2cm}\caption{SST $k-\omega$ turbulence model case, eddy viscosity fields: (a) shows the FOM eddy viscosity, while in (b) one can see the P-ROM eddy viscosity, and finally in (c) we have the Mixed-ROM eddy viscosity.}\label{fig:nut_fieldsKomega} 
\end{figure}

A further simulation campaign has been carried out with a different, SST $k-\omega$ turbulence model, to evaluate how responsive the
hybrid Mixed-ROM and the P-ROM results are with respect to the turbulence model employed for the FOM simulations. Thus, a new set of SST $k-\omega$ FOM simulations
has been run using the same inflow velocity values as in $k-\epsilon$ model case. The snapshots generated have been again used to train
both reduced models considered. \autoref{fig:u_fieldsKomega}, \autoref{fig:p_fieldsKomega} and \autoref{fig:nut_fieldsKomega} show the velocity, pressure and eddy viscosity fields obtained by the FOM, the P-ROM and the Mixed-ROM for the inflow velocity value $ U^* = 7.0886$, respectively. Again, the Mixed-ROM results appear in good qualitative agreement with their SST $k-\omega$ FOM counterparts, while the same cannot be claimed for the P-ROM results. By a quantitative standpoint, the $L^2$ relative errors between the FOM and the Mixed-ROM velocity and pressure fields are respectively $\epsilon_u=0.8088$ $\%$ and  $\epsilon_p=0.7329$ $\%$. As for the P-ROM results, the corresponding errors are $\epsilon_u=0.8177$ $\%$ and $\epsilon_p=22.3972$ $\%$, respectively.\par
\begin{figure}
  \centering
  \begin{minipage}[b]{0.5\linewidth}
    \centering
    \includegraphics[width=0.98\linewidth]{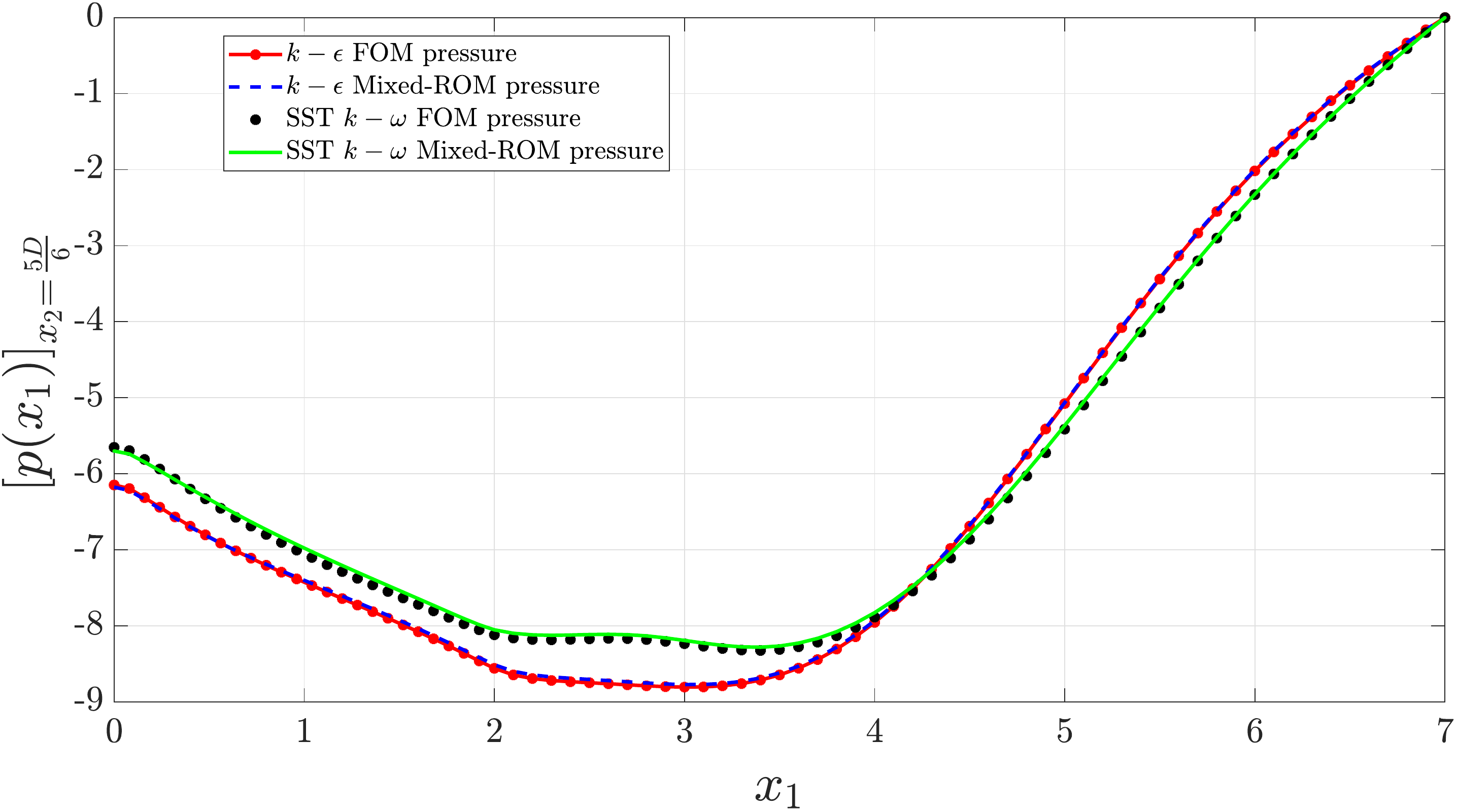}
  \end{minipage}
  \vspace{-0.2cm}\caption{The pressure fields obtained using both $k-\epsilon$ and SST $k-\omega$ turbulence models and the Mixed-ROM ones. The plot is for the pressure value along the $x_1$ direction keeping the value of $x_2$ fixed at half the maximum height.}\label{fig:pressureSlice}\par
\end{figure}

The FOM fields obtained solving the RANS equations with the two different turbulence models have quantitatively speaking different values across the domain (except for the velocity). In order to give a clear idea about how accurate was the reduction performed by the Mixed-ROM regardless of the turbulence model employed at full order level, one may plot the FOM and the Mixed-ROM pressure fields (obtained by the two turbulence models) for a fixed value along the $x_2$ axis (the perpendicular axis) versus the values along the $x_1$ axis (the horizontal one). The last test is done in \autoref{fig:pressureSlice}, where one can see the FOM and the Mixed-ROM pressure fields along the horizontal direction at a fixed height of $x_2 = \frac{5D}{6}$ which is half the height of the domain. The plot is done for both $k-\epsilon$ and SST $k-\omega$. As can be appreciated from the figure, the Mixed-ROM was successful in obtaining pressure field values which are close the FOM ones regardless of the turbulence model utilized at full order level. This accomplishes one of the main goals of the Mixed-ROM developed in this work. \par
\begin{figure}
  \centering
 \begin{minipage}[b]{0.5\linewidth}
    \centering
    \includegraphics[width=0.98\linewidth]{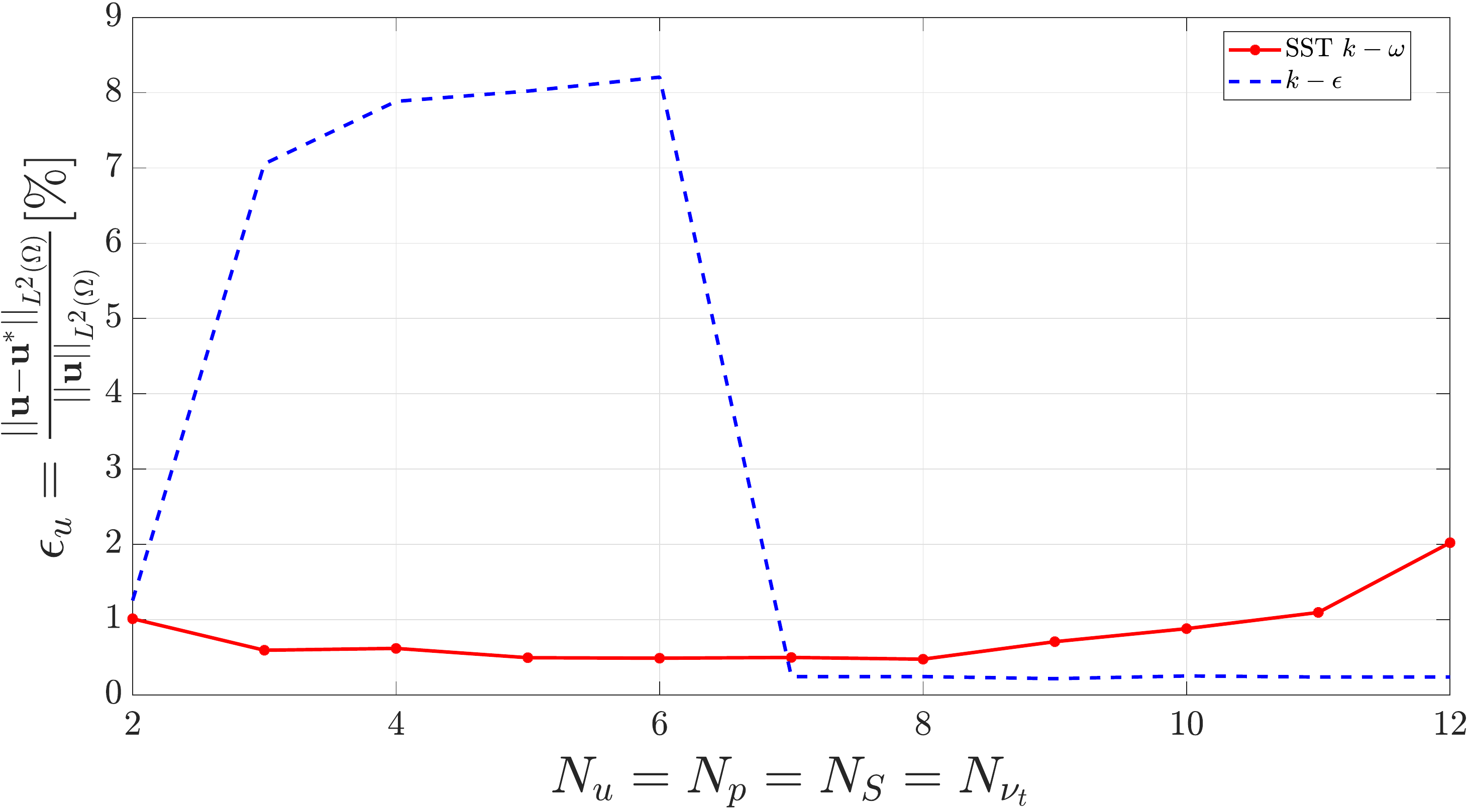} 
    \scriptsize(a) 
    \end{minipage}
  \begin{minipage}[b]{0.5\linewidth}
    \centering
    \includegraphics[width=0.98\linewidth]{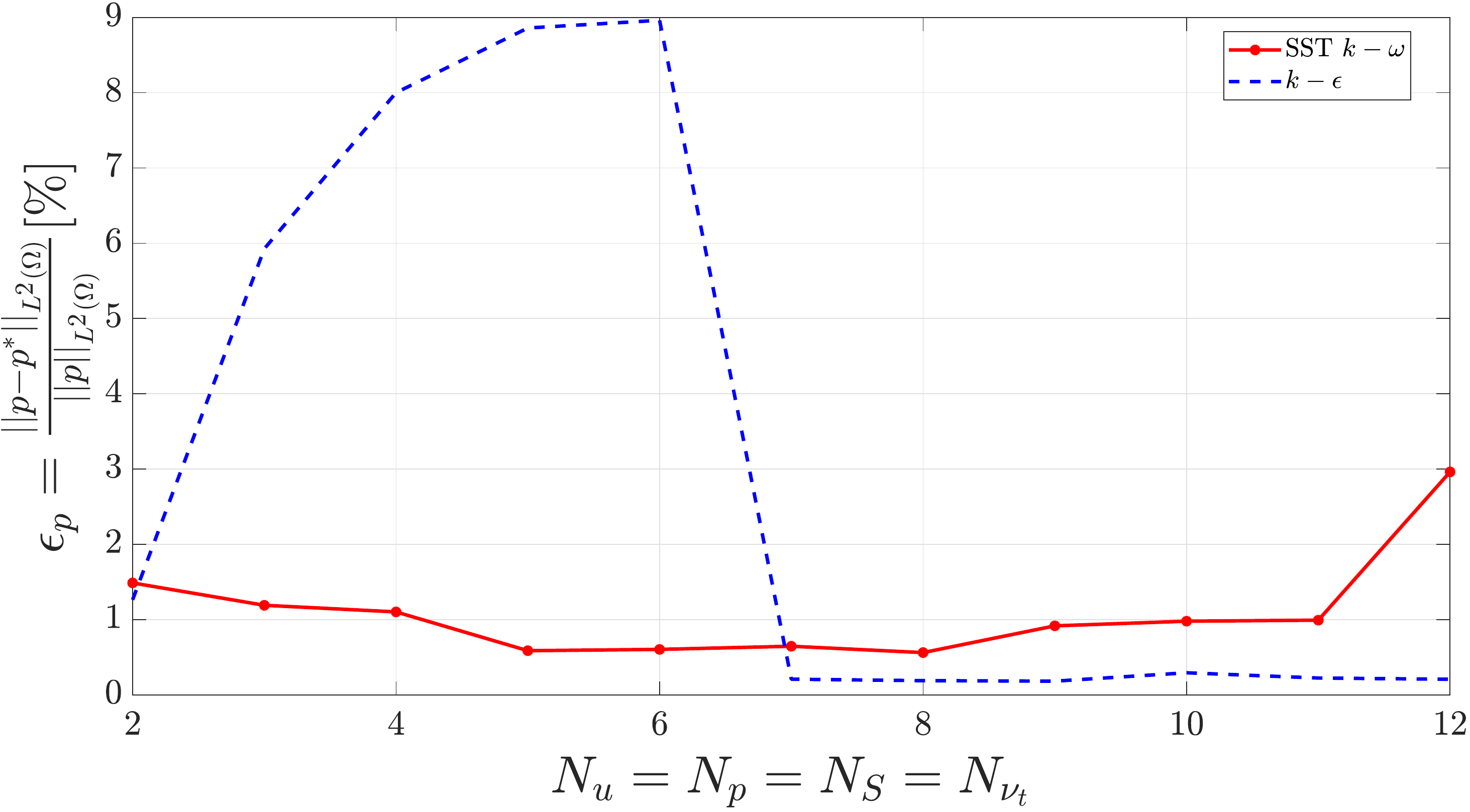}
        \scriptsize(b) 
  \end{minipage}
 \vspace{-0.2cm}\caption{The mean of the $L^2$ relative errors for all the online samples versus the number of modes used in the online stage. The convergence analysis is done for both Mixed-ROM models obtained with two different turbulence models at the full order level which are $k-\epsilon$ and SST $k-\omega$. The errors are reported in percentages, in (a) we have the velocity fields mean error, while in (b) the pressure fields mean error .}\label{fig:conv_ana_steady} 
\end{figure}

Finally, the convergence analysis for the Mixed-ROM results is shown in \autoref{fig:conv_ana_steady}. The plots show the mean $L^2$ relative error for all the $80$ samples used in the cross validation test in the online stage, as a function of the number of modes used. As previously
mentioned, the number of modes used for velocity ($N_u$), pressure ($N_p$), supremizer ($N_S$) and eddy viscosity ($N_{\nu_t}$) was kept
uniform in these preliminary tests. The plots indicate that for the problem considered, the Mixed-ROM results exhibit fast convergence to the FOM 
solution for both $k-\epsilon$ and SST $k-\omega$. Yet, after less then ten modes, the convergence appears to stall, as the error
settles on non zero, but fairly acceptable values. This is likely due to the fact that as the number of modes grow, the gain in accuracy
becomes only marginal compared to the $\nu_t$ field interpolation error.


\subsection{Unsteady case}\label{sec:numerical_unsteady}
\begin{figure}
{
  \centering
  \begin{minipage}[b]{0.57\linewidth}
    \centering
    \includegraphics[width=\linewidth]{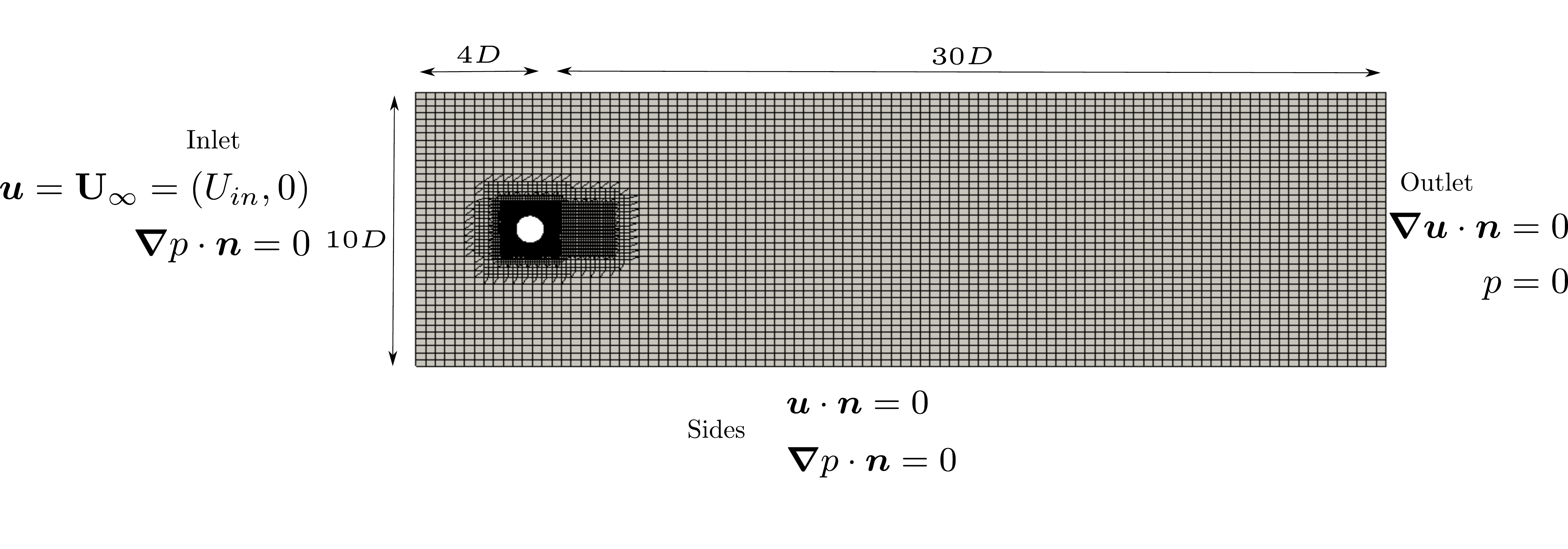}\label{fig:meshUnsteady} 
    \scriptsize(a)
    \end{minipage} 
  \begin{minipage}[b]{0.42\linewidth}
    \centering
    \includegraphics[width=\linewidth]{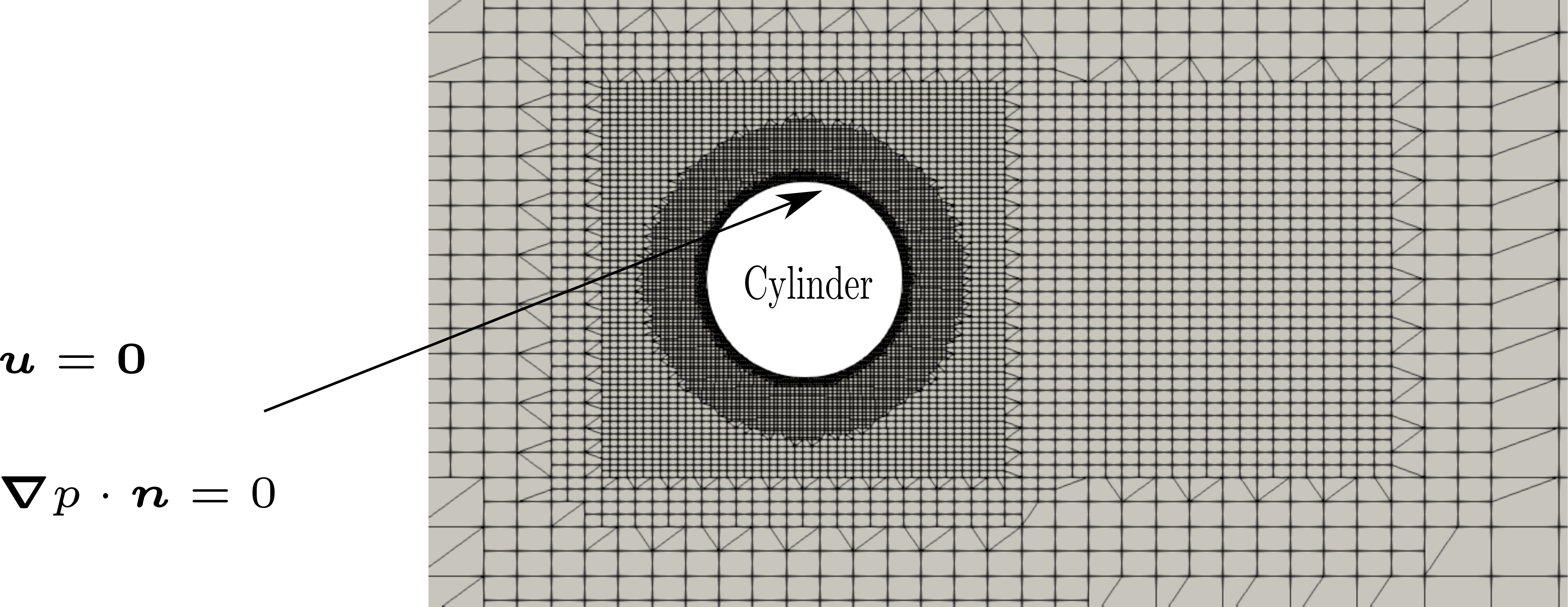}\label{fig:meshzoomUnsteady}
        \scriptsize(b)
        \end{minipage} 
  \vspace{-0.2cm}\caption{(a) The OpenFOAM mesh used in the simulations for the unsteady case of the flow around a circular cylinder. (b) A picture of the mesh zoomed near the cylinder.}\label{fig:MeshUnsteady} 
}
\end{figure}

The present subsection presents the application of the Mixed-ROM on a parametrized non-stationary case. The problem considered is that of the turbulent and unsteady flow around a circular cylinder. For more details on such classical benchmark flow, the reader may refer to \cite{zdravkovichflow,zdravkovich2003flow}. In the framework of the unsteady ROM analysis described in Section \ref{sec:ROM}, time is treated as one of the parameters characterizing the PDE problem. A further parameter considered in this numerical investigation is the Reynolds number associated with the inflow velocity. The domain and the $2$D computational grid used
are depicted in \autoref{fig:MeshUnsteady}, which also reports the boundary conditions imposed in the simulations. In the picture, all the
lengths reported are referred to the problem characteristic length which is the diameter of the cylinder $D=1$ \si{m}. The grid features
$11644$ cells, while the physical viscosity $\nu$ is equal to $10^{-4}$ \si{m^2 \per s}. Uniform and constant horizontal
velocities $\mathbf{U}_\infty = (U_{in},0)$ with $U_{in} ~\in [7.5,12] \quad \si{m\per s}$ (corresponding to Reynolds number in  the range of $[7.5 \times 10^4,1.2 \times 10^5]$) were imposed at the inlet boundary, and the simulations evolve in time
from rest until a final periodic regime solution is reached.\par

In this test, the turbulence model considered is SST $k-\omega$. As for the numerical schemes used to set up the FOM simulations,
time discretization is done using backward Euler scheme, while gradients are approximated using Gauss scheme. The convection term is discretized through a $2$nd order bounded upwind divergence scheme which utilized upwind interpolation weights, with an explicit correction based on the local cell gradient. Finally, the diffusive term is discretized by Gauss linear scheme.\par

The main objective of this numerical test is that of building a reduced order model which can successfully reproduce the flow fields
corresponding to the final periodic regime solution. For such reason it is important to properly select the time window from which
snapshots will be taken and ensure that it contains enough solution cycles ($1.5-2$ cycles at least). The evaluation
of the cycles period length has been carried out through Fourier analysis of the FOM time signal of lift and drag fluid dynamic forces
acting on the cylinder.

As mentioned, the physical parameter varied in the numerical tests, is the horizontal velocity at the inlet called $U_{in}$. Ten samples are taken from the velocity range $[7.5,12] \quad \si{m\per s}$, and for each of such samples the FOM simulator was run and snapshots were acquired at time steps covering approximately two cycles after reaching the final regime solution. It must be remarked that the extent of the time windows in which snapshots were taken was adapted for each velocity sample to track the solution period variations due to the change of the frequency of vortex shedding of the system.

\begin{table}[htp]
\centering
\begin{tabular}{ c | c | c | c | c }
  \textbf{Parameter sample : $U_{in}$ in $\si{m\per s}$}  & \textbf{FOM time step} & \textbf{Snapshot acquiring time}\\
    \hline  			
  $7.5$ & $0.0004$ & $0.008$ \\
  \hline
  $8$ & $0.0004$ & $0.008$ \\
  \hline  
  $8.5$ & $0.00035$ & $0.007$ \\
  \hline  
  $9$ & $0.0003$ & $0.006$ \\
  \hline  
  $9.5$ & $0.0003$ & $0.006$ \\
  \hline  
  $10$ & $0.0003$ & $0.006$ \\
  \hline  
  $10.5$ & $0.0003$ & $0.006$ \\
  \hline  
  $11$ & $0.0003$ & $0.006$ \\
  \hline  
  $11.5$ & $0.00025$ & $0.005$ \\
  \hline  
  $12$ & $0.00025$ & $0.005$ \\
  \hline  
\end{tabular}\caption{Offline parameter samples and the corresponding snapshots data}
\label{tab:parOffline}
\end{table}

\begin{figure}
\centering
\includegraphics[width=0.7\textwidth]{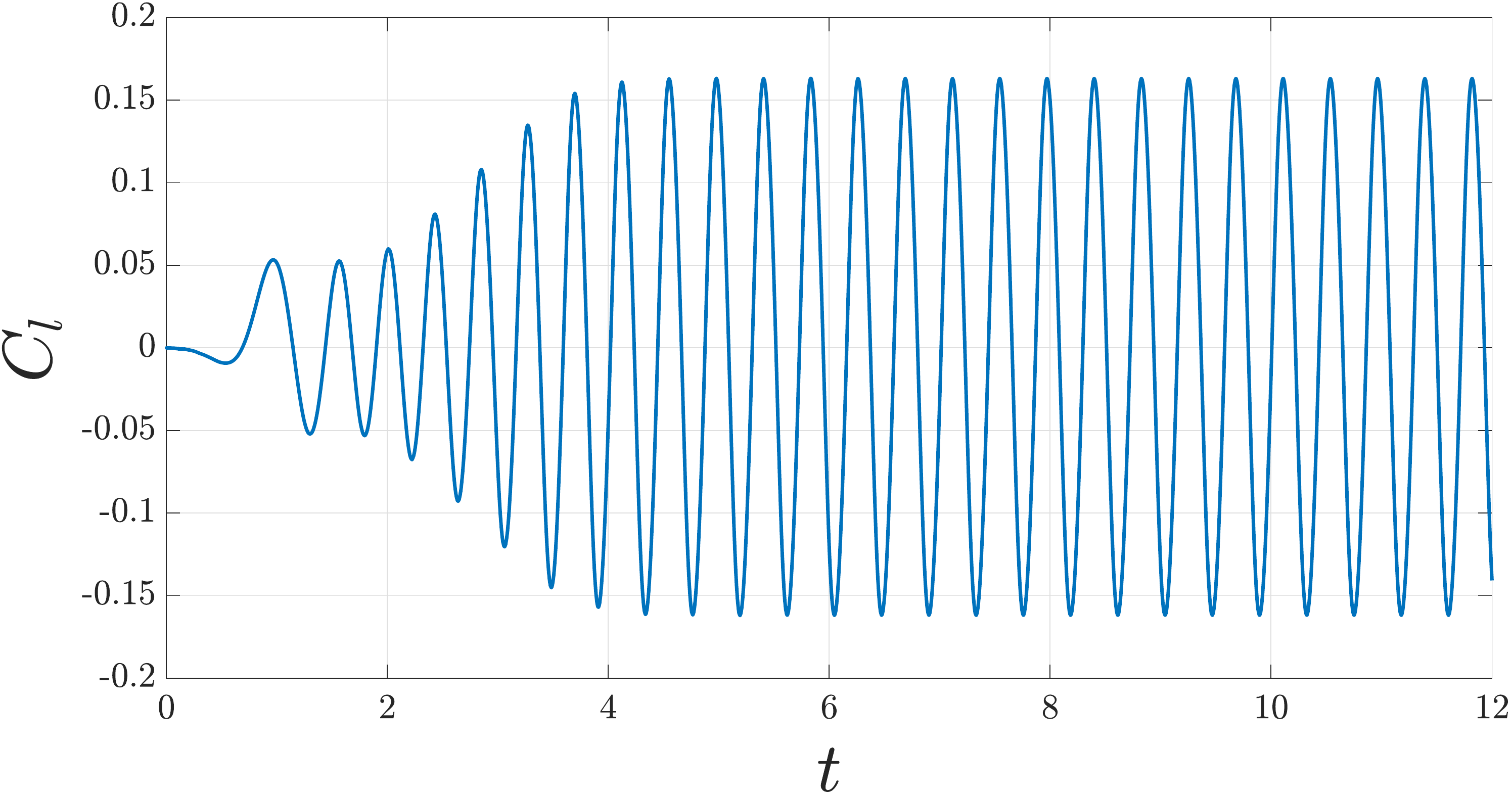}
\caption{The lift coefficient curve for parameter sample $U_{in} = 10$ $\si{m\per s}$.}\label{fig:Lift_cy}
\end{figure}

As an example of this procedure, if we consider the inlet velocity of $10$ $\si{m\per s}$, the FOM simulation has been run for $12$ seconds using the OpenFOAM solver pimpleFoam which adapts the timestep so as to keep the Courant
number $CFL$ \cite{Courant1928,Courant1967} under a prescribed value ${CFL}_{max} = 0.9$. \autoref{fig:Lift_cy} depicts the resulting lift coefficient curve, which is obtained from the lift $L$ as $C_l=\frac{L}{\frac{1}{2}\rho U^2 D}$. The non-uniformly spaced time signal of the lift coefficient has been interpolated on equally distributed time nodes so as to allow the use of Fast Fourier Transform (FFT) for the computation of the time period corresponding to the principal, vortex shedding, frequency. The time period computed is $0.4299$ \si{s}, corresponding to a Strouhal number \cite{Strouhal1878} of St=$0.2326$, which is in line with well assessed experimental value of approximately $0.20$ \cite{Blevins2009}. After this value was available, the simulations have been extended keeping a fixed time step of $0.0003$ \si{s} to start acquiring snapshots which cover two periods at least. More specifically, the simulations were run for $1.2$ \si{s} additional seconds, saving snapshots of the flow field with a $0.006$ \si{s} time rate so as to finally obtain $200$ snapshots. We remark that, to be as consistent as possible, the time step imposed in the resolution of the Mixed-ROM dynamical system \eqref{eq:PODsys} at the reduced level, has been the same one used for the FOM simulations.\par

The offline stage was carried out taking $200$ snapshots for each parameter sample. Table \ref{tab:parOffline} shows the values of the parameters and the corresponding values of the simulation time step and the time interval at which snapshots were acquired.\par

\begin{figure}
\centering
\includegraphics[width=0.7\textwidth]{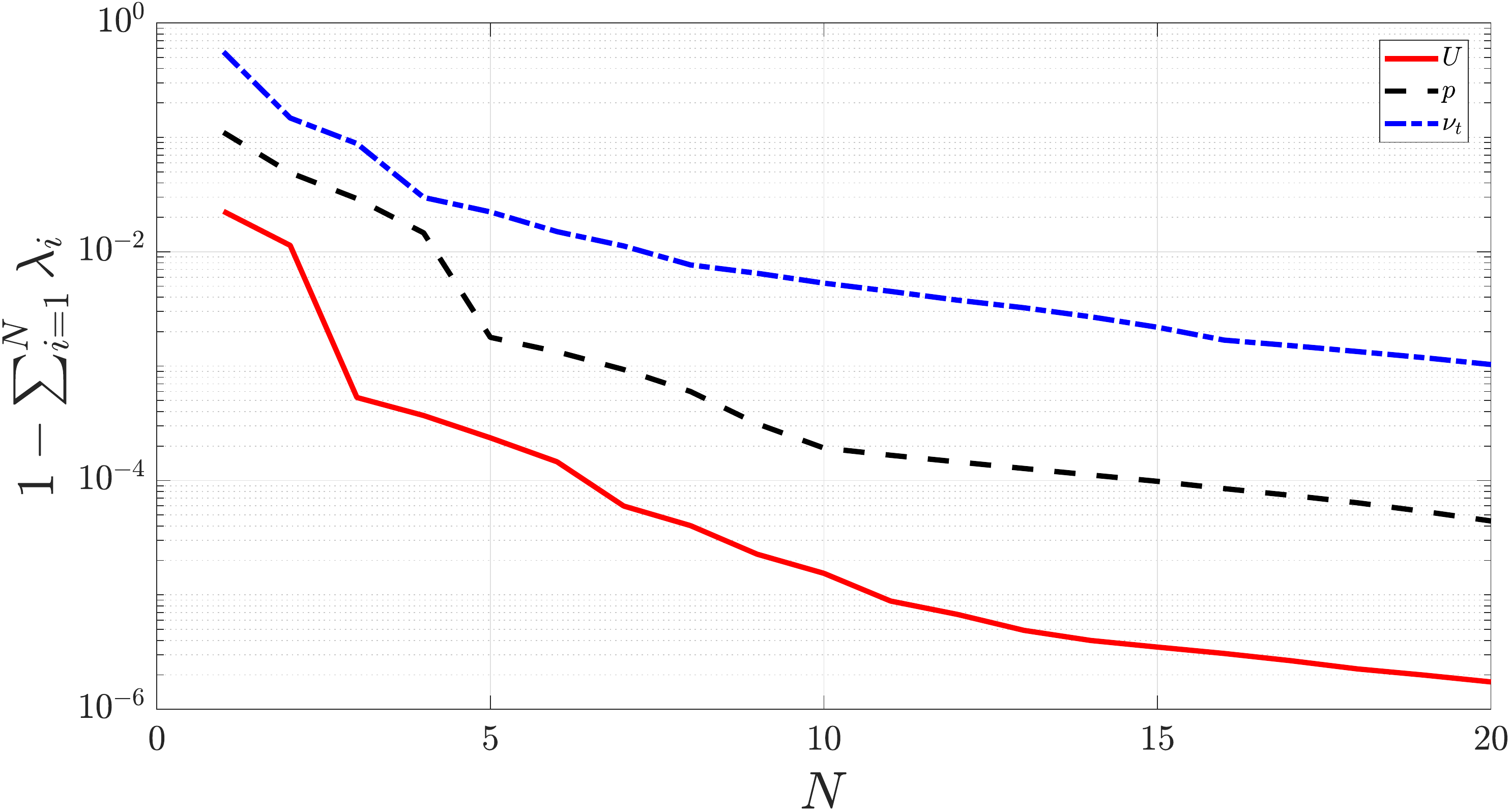}
\caption{Cumulative ignored eigenvalues decay. In the plot, the solid red line refers to the velocity eigenvalues, the dashed black line indicates the pressure eigenvalues and the dash-dotted blue line finally refers to the eddy viscosity eigenvalues.}
\label{fig:Cum_cy}\end{figure}

The non-homogeneous Dirichlet boundary condition at the inlet is enforced with the penalty method. POD modes have been obtained applying POD analysis to snapshots matrices of velocity, pressure and eddy viscosity fields. \autoref{fig:Cum_cy} depicts the decay of the cumulative eigenvalues corresponding to the three correlation matrices. The supremizer problem was then solved for each of the pressure modes, to finally obtain the supremizer modes added to the velocity ones.\par

The online resolution of the Mixed-ROM system requires that an interpolation strategy is used to obtain the eddy viscosity coefficient vector $\bm{g}$ at each time step $t^*$ of the simulation corresponding to the parameter value $U_{in}^*$. More specifically, at each time instant $\bm{g}(t^*,U_{in}^*)$ should be obtained through interpolation --- with respect to the combined time-parameter vector --- from its values corresponding to the snapshots. Yet, while all snapshots are contained in the aforementioned $1.2$ \si{s} time window,  the online time integration must extend for much longer times. This means that for each instant outside the time window of the original snapshots, $\bm{g}$ must be in fact extrapolated. To avoid such problem, the eddy viscosity coefficients are obtained through RBF interpolation from the reduced order velocity coefficients vectors of $\bm{a}$ and $\bm{\dot{a}}$ (\autoref{eq:inter_v}). As the values of the reduced velocity solution vector
components oscillate between minima and maxima over time integration, using $\bm{a(t)}$ instead of $t$ as the RBF interpolation variable has in fact the convenient benefit of avoiding extrapolation. Of course, this is true if the values of the $a$ vector components obtained during the ROM time integration fall within the bounds of the FOM snapshots. For such reason, it is clear that the accuracy of such interpolation outside the offline snapshots window highly depends on how close the current solution vector $\bm{a}$ is to the vectors of the $L^2$ projection coefficients used in the offline stage for training the RBF.\par

As for the dependence on the inlet velocity parameter, the results presented have been obtained by splitting $\nu_t$ into its time average and its fluctuating part. The inlet velocity parameter dependence has been then only enforced on the time average degrees of freedom $\overline{\bm{g}}$, while the aforementioned interpolation based on the reduced velocity vector has been only applied to the fluctuating part $\bm{g}$. This means that in $M=10$ different average eddy viscosity fields were computed by taking the average of the set of snapshots which correspond to one value of the ten inlet velocity samples. The average reduced vector $\overline{\bm{g}}$ has been then obtained from $U_{in}^*$ in the online stage using linear interpolation, while the reduced vector $\bm{g}$ is obtained from RBF interpolation with respect to $\bm{a}$ and $\bm{\dot{a}}$. Finally, the initial values for all vectors $\bm{a}(0,U_{in}^*)$, $\bm{b}(0,U_{in}^*)$ and $\bm{g}(0,U_{in}^*)$ are obtained from the inlet velocity parameter using linear interpolation as well (based on the values of the initial $L^2$ projection vectors of $\bm{a}(0,U_{in})$, $\bm{b}(0,U_{in})$ and $\bm{g}(0,U_{in})$). 

\begin{figure}
  \centering
 \begin{minipage}[b]{0.5\linewidth}
    \centering
    \includegraphics[width=0.98\linewidth]{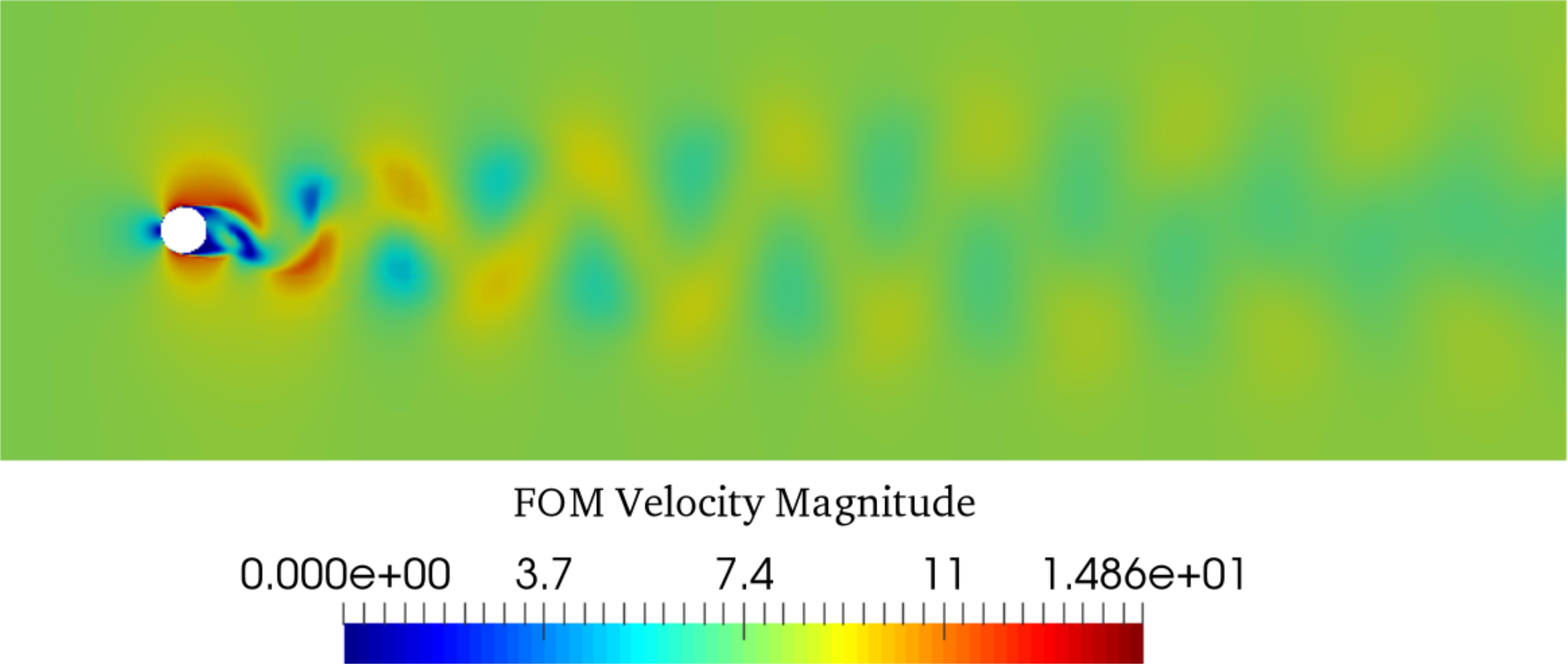} 
    \scriptsize(a) 
    \end{minipage}
  \begin{minipage}[b]{0.5\linewidth}
    \centering
    \includegraphics[width=0.98\linewidth]{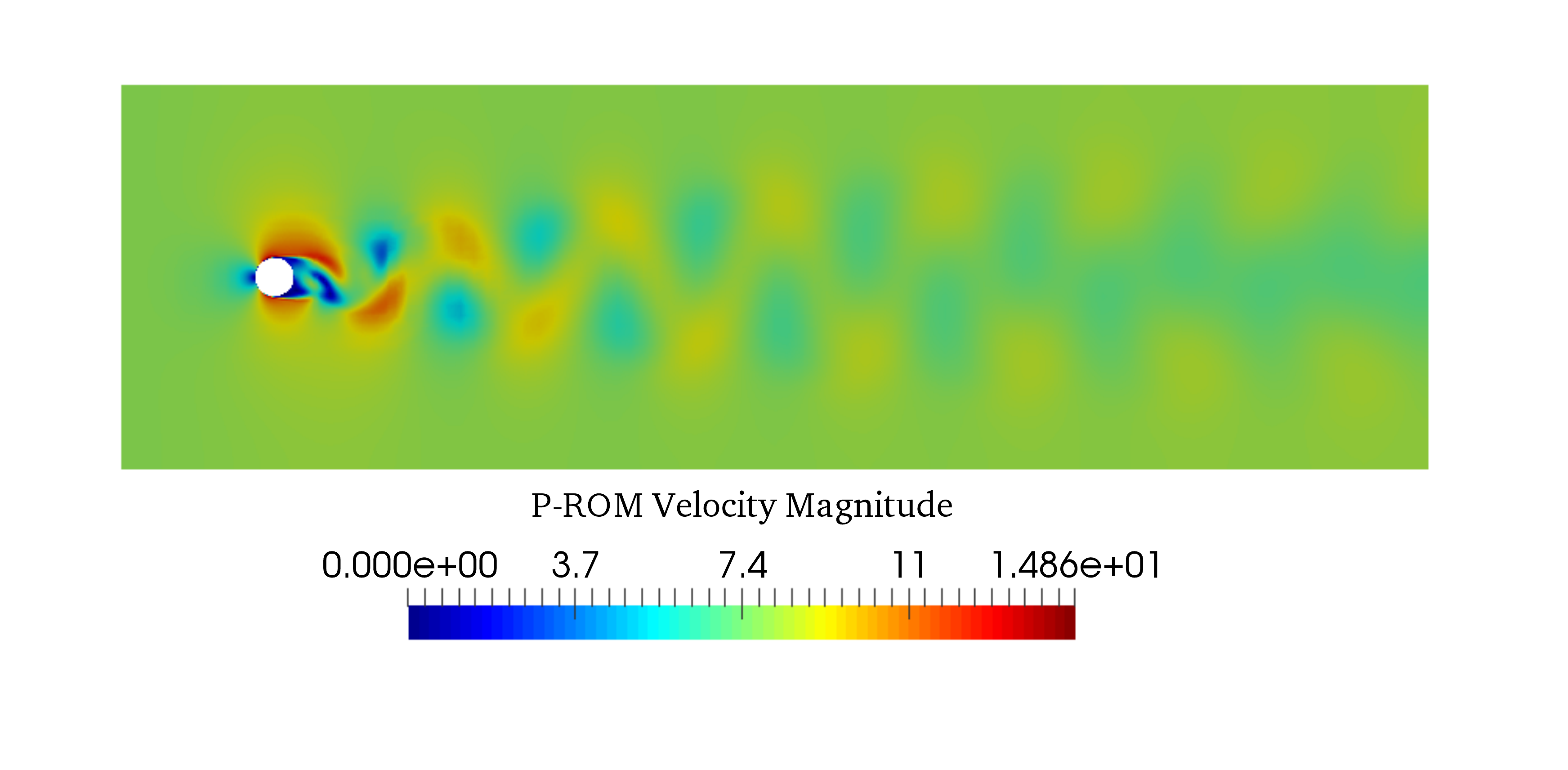}
        \scriptsize(b) 
  \end{minipage}
   \begin{minipage}[b]{0.5\linewidth}
    \centering
        \includegraphics[width=0.98\linewidth]{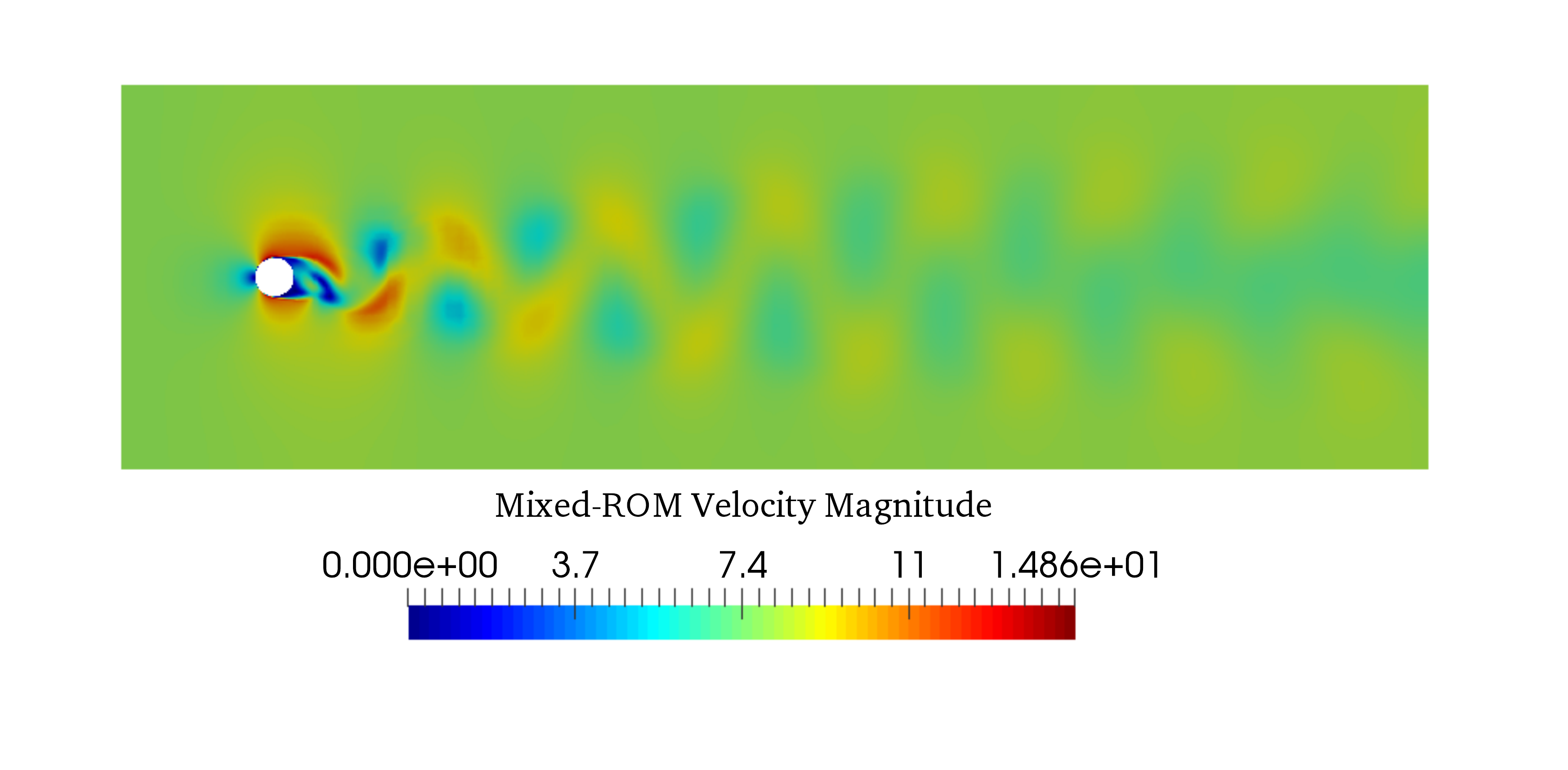} 
    \scriptsize(c)
  \end{minipage}
 \vspace{-0.2cm}\caption{Velocity fields for the parameter value $U_{in} = 7.75$ $\si{m\per s}$ at $ t = 2.8$  \si{s}: (a) shows the FOM velocity, while in (b) one can see the P-ROM velocity, and finally in (c) we have the Mixed-ROM velocity.}\label{fig:u_fields_cy_2.8} 
\end{figure}
\begin{figure}
  \centering
 \begin{minipage}[b]{0.5\linewidth}
    \centering
    \includegraphics[width=0.98\linewidth]{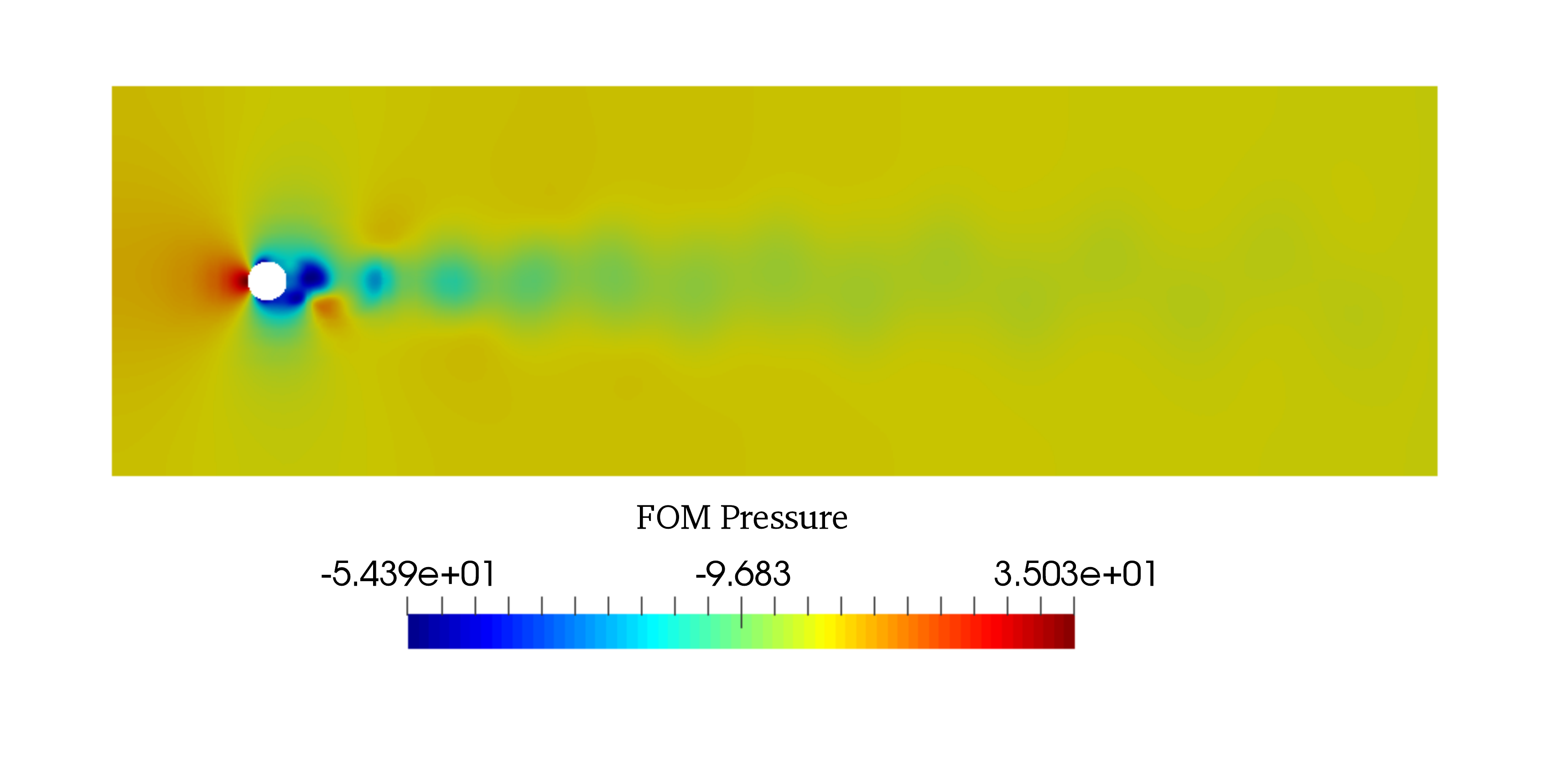} 
    \scriptsize(a) 
    \end{minipage}
  \begin{minipage}[b]{0.5\linewidth}
    \centering
    \includegraphics[width=0.98\linewidth]{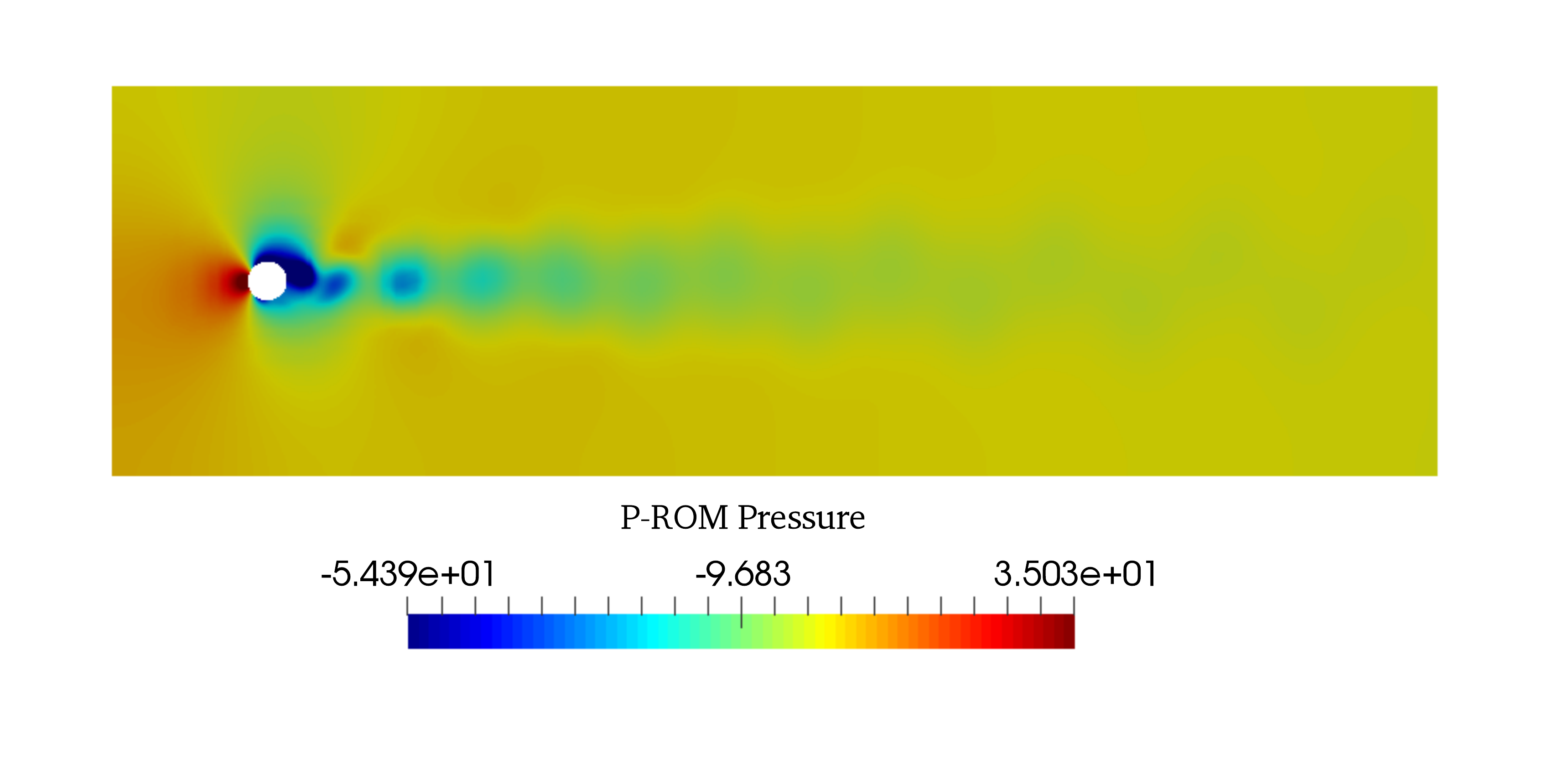}
        \scriptsize(b) 
  \end{minipage}
   \begin{minipage}[b]{0.5\linewidth}
    \centering
        \includegraphics[width=0.98\linewidth]{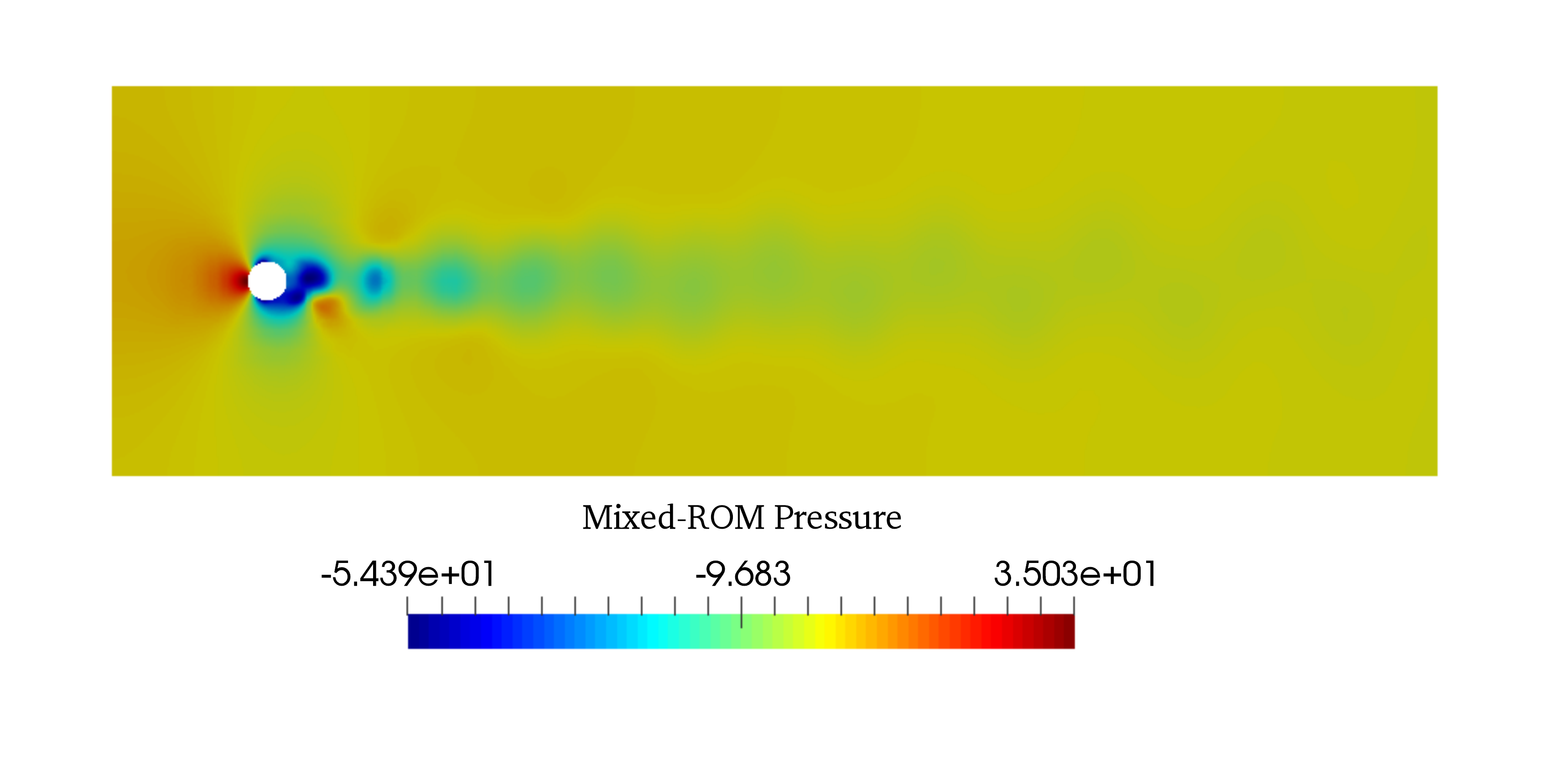} 
    \scriptsize(c)
  \end{minipage}
 \vspace{-0.2cm}\caption{Pressure fields for the parameter value $U_{in} = 7.75$ $\si{m\per s}$ at $ t = 2.8$  \si{s}: (a) shows the FOM pressure, while in (b) one can see the P-ROM pressure, and finally in (c) we have the Mixed-ROM pressure.}\label{fig:p_fields_cy_2.8} 
\end{figure}
\begin{figure}
  \centering
 \begin{minipage}[b]{0.5\linewidth}
    \centering
    \includegraphics[width=0.98\linewidth]{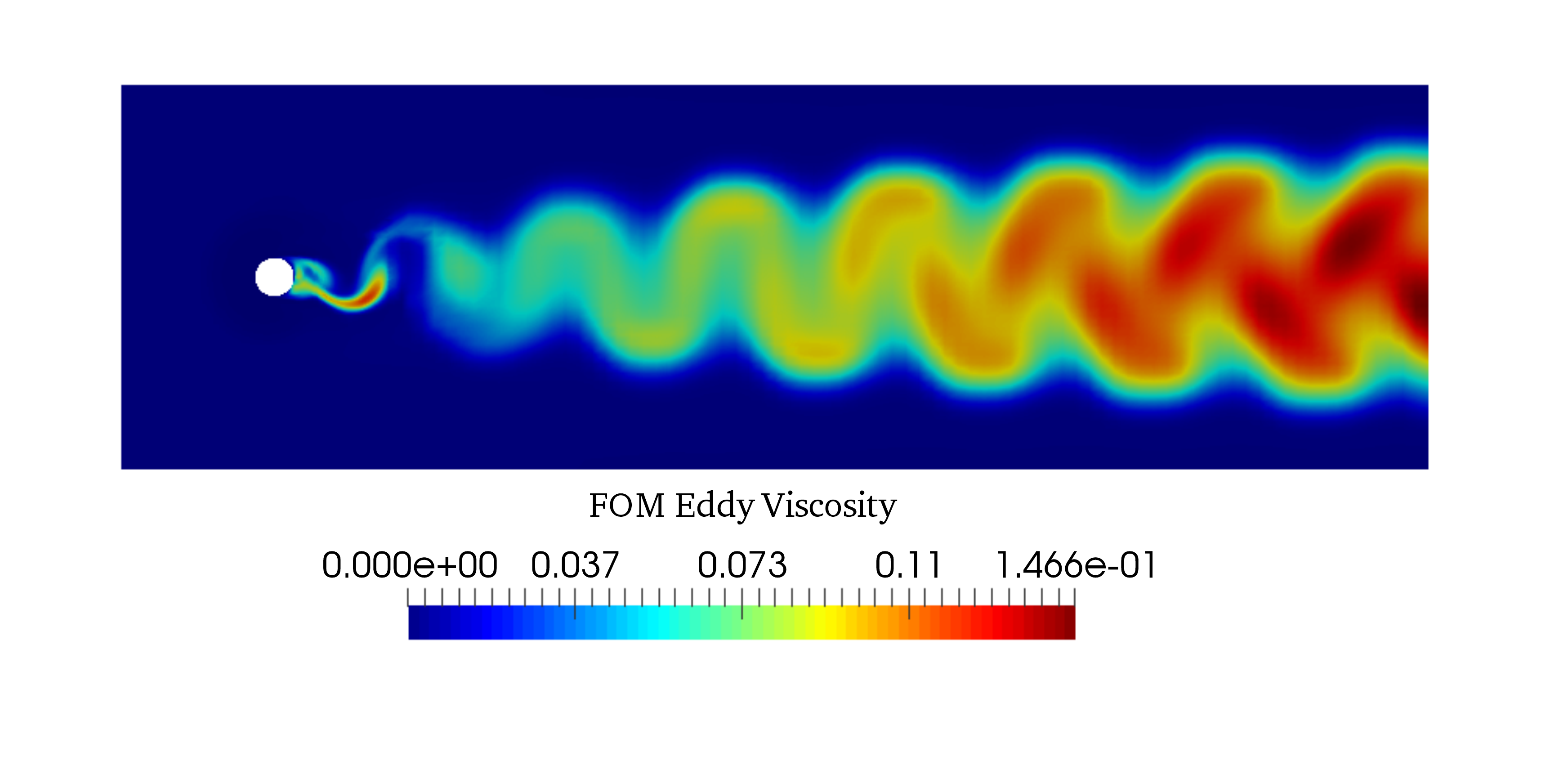} 
    \scriptsize(a) 
    \end{minipage}
  \begin{minipage}[b]{0.5\linewidth}
    \centering
    \includegraphics[width=0.98\linewidth]{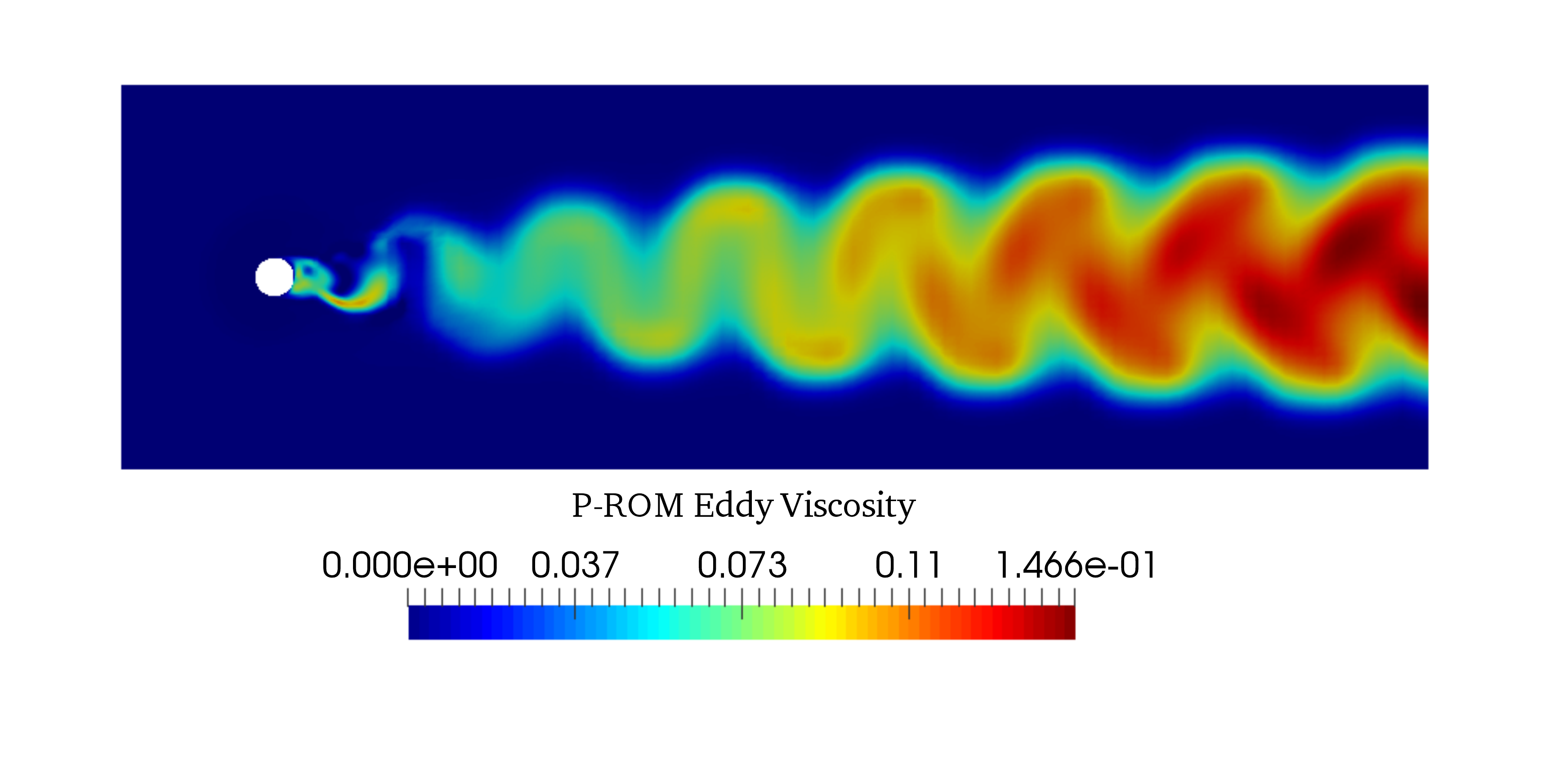}
        \scriptsize(b) 
  \end{minipage}
   \begin{minipage}[b]{0.5\linewidth}
    \centering
        \includegraphics[width=0.98\linewidth]{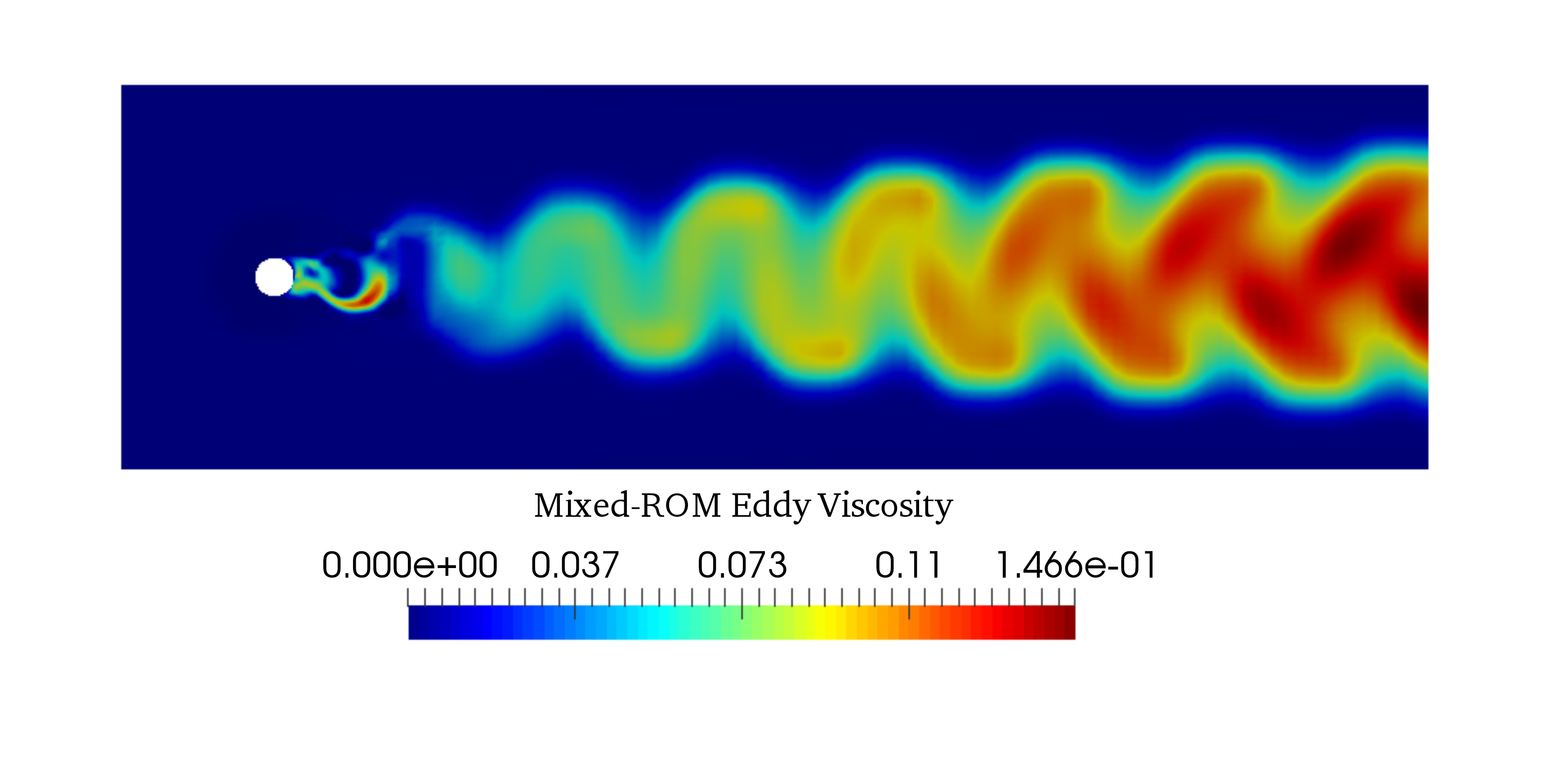} 
    \scriptsize(c)
  \end{minipage}
 \vspace{-0.2cm}\caption{Eddy viscosity fields for the parameter value $U_{in} = 7.75$ $\si{m\per s}$ at $t = 2.8$ \si{s}: (a) shows the FOM eddy viscosity, while in (b) one can see the P-ROM eddy viscosity, and finally in (c) we have the Mixed-ROM eddy viscosity.}\label{fig:nut_fields_cy_2.8} 
\end{figure}

The first numerical test is a cross validation test for the parameter value $U_{in} = 7.75$ $\si{m\per s}$, not contained in the samples set. Once the offline phase was completed with the computation of the reduced order matrices, system \eqref{eq:PODsysAve} was solved for $\bm{a}$ and $\bm{b}$ and the Mixed-ROM solution fields were computed. A comparison is made between the fields obtained by the FOM solver and the ones computed by both the Mixed-ROM and the P-ROM ones. The FOM simulator was run for enough time to reach a periodic regime and then it was launched again with a constant simulation time step of $0.0004$ $\si{s}$ exporting the solution fields every $0.008$ $\si{s}$. The total FOM simulation time for this test was $8$ $\si{s}$ which contained $13$ periods. The starting time of the simulation of the final periodic regime of all tests in this section is set to $0$. The first results shown correspond to the flow fields computed by the FOM, the P-ROM and the Mixed-ROM at $t=2.8$ $\si{s}$. \autoref{fig:u_fields_cy_2.8} shows the velocity fields while \autoref{fig:p_fields_cy_2.8} and \autoref{fig:nut_fields_cy_2.8} present the pressure and eddy viscosity fields, respectively. \autoref{fig:u_fields_cy_2.8} indicates that both the P-ROM and the Mixed-ROM are able to obtain accurate velocity prediction. The $L^2$ relative norm of the error committed by the two models is in fact $1.218$ $\%$ and $0.6921$ $\%$, respectively. On the other hand the pressure fields shown in \autoref{fig:p_fields_cy_2.8} suggest that the P-ROM model fails in giving sufficiently accurate results for the reduction of the pressure fields. In fact, the P-ROM pressure field does not match the FOM one. This is particularly true in the region near the cylinder, which is of course crucial for an accurate reproduction of the forces acting on the body. As for the Mixed-ROM, the reduced pressure field appear to be in closer agreement with the FOM one. This is confirmed by more quantitative assessments, as the $L^2$ relative norm of the Mixed-ROM pressure field error is $4.7894$ $\%$, while that of the P-ROM is $29.5958$ $\%$. We must remark that to obtain the best results with each model developed in this numerical test, the number of modes used in the online stage for the P-ROM is $9$, while in the Mixed-ROM case $20$ velocity modes were used and $10$ modes were employed for pressure, supremizers and eddy viscosity.\par



\begin{figure}
  \centering
 \begin{minipage}[b]{0.5\linewidth}
    \centering
    \includegraphics[width=0.98\linewidth]{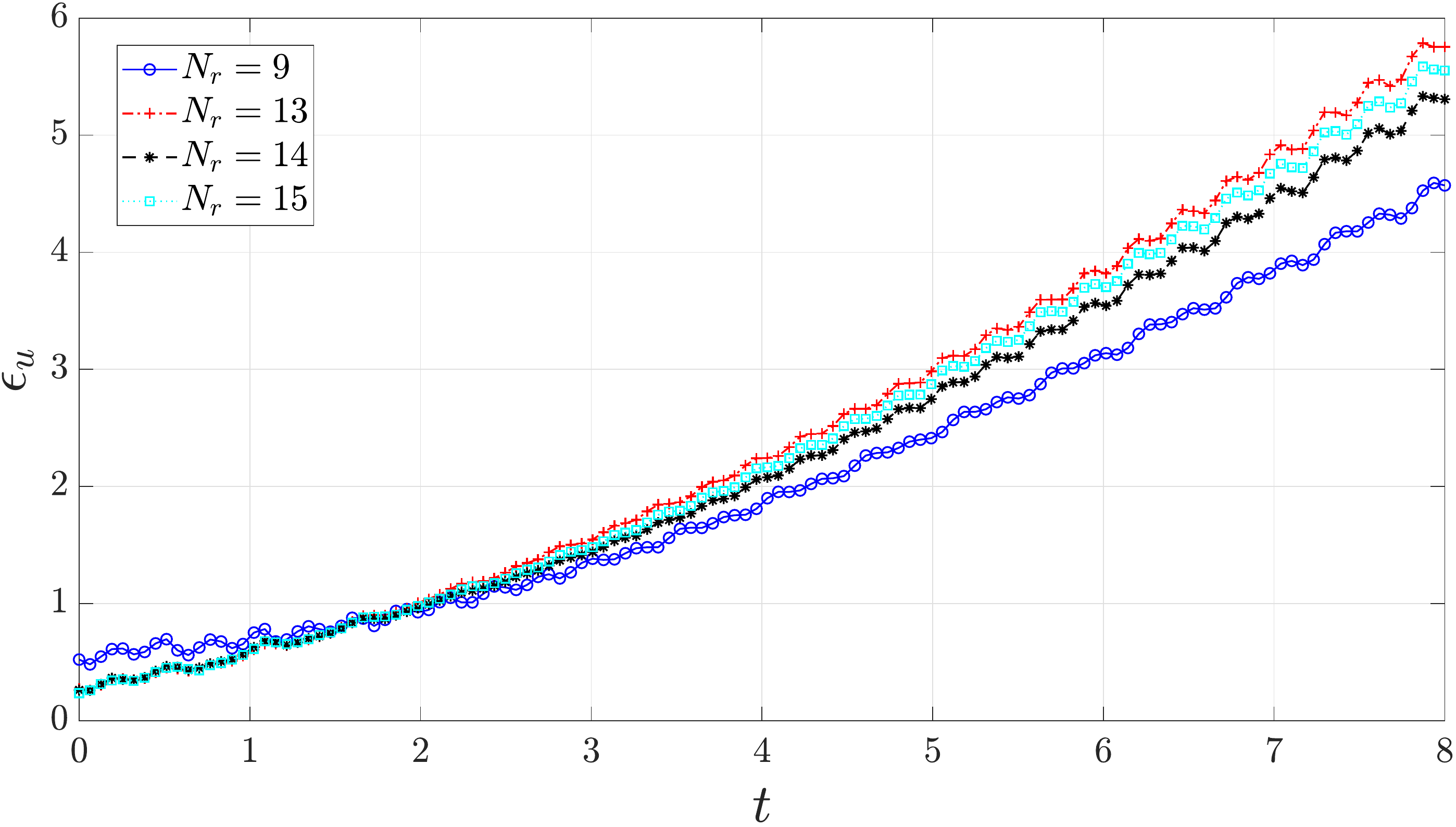} 
    \scriptsize(a) 
    \end{minipage}
  \begin{minipage}[b]{0.5\linewidth}
    \centering
    \includegraphics[width=0.98\linewidth]{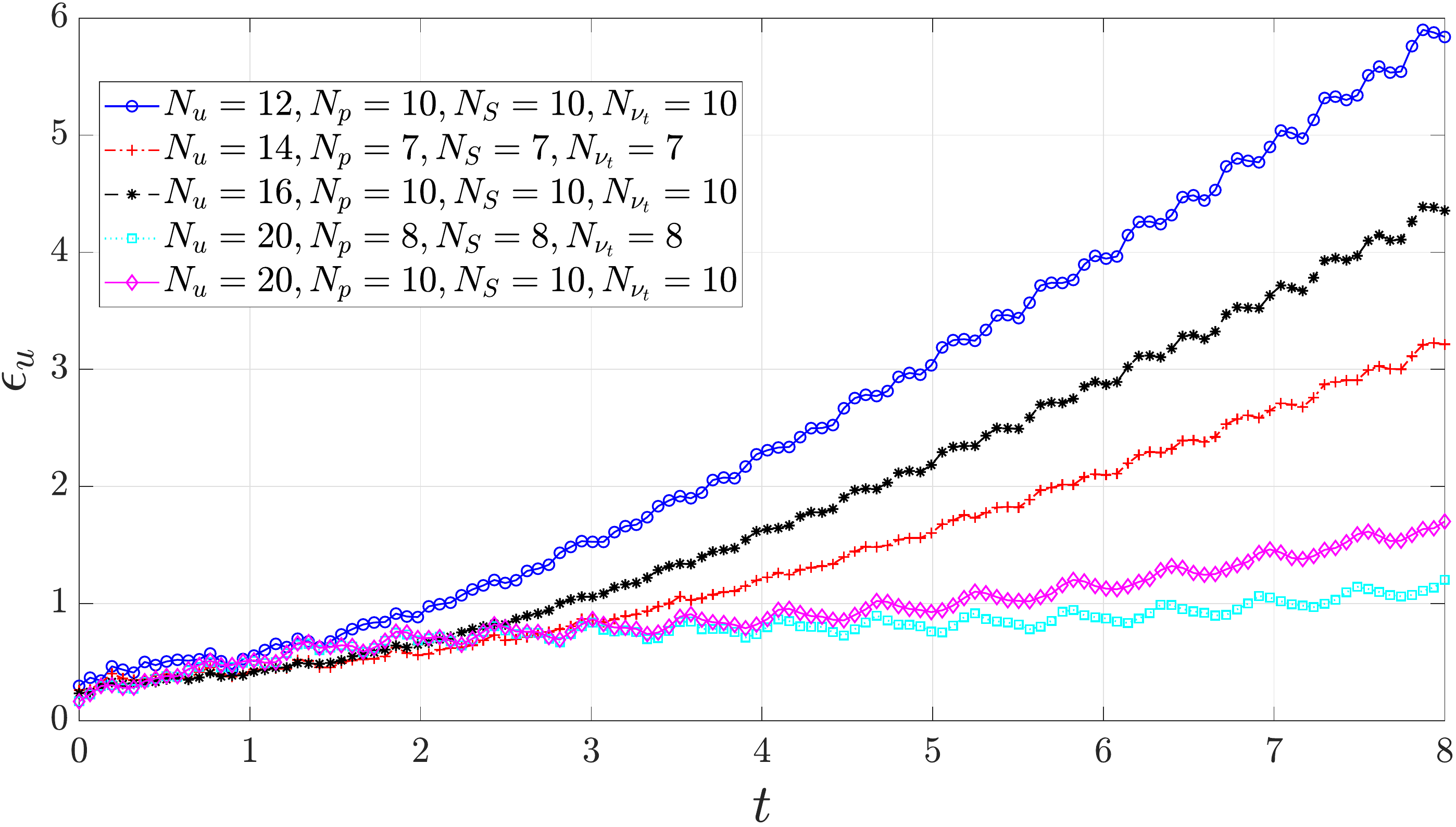}
        \scriptsize(b) 
  \end{minipage}
\caption{The time evolution of the $L^2$ relative errors of the velocity reduced approximations for both the P-ROM and the Mixed-ROM models. The curves correspond to the case run with the parameter value $U_{in} = 7.75$ $\si{m\per s}$ : (a) shows the error curve for the P-ROM model. Figure (b) depicts the case of the Mixed-ROM model. The error values in both graphs are in percentages.}\label{fig:errU_U7.75}
\end{figure}

\begin{figure}
  \centering
 \begin{minipage}[b]{0.5\linewidth}
    \centering
    \includegraphics[width=0.98\linewidth]{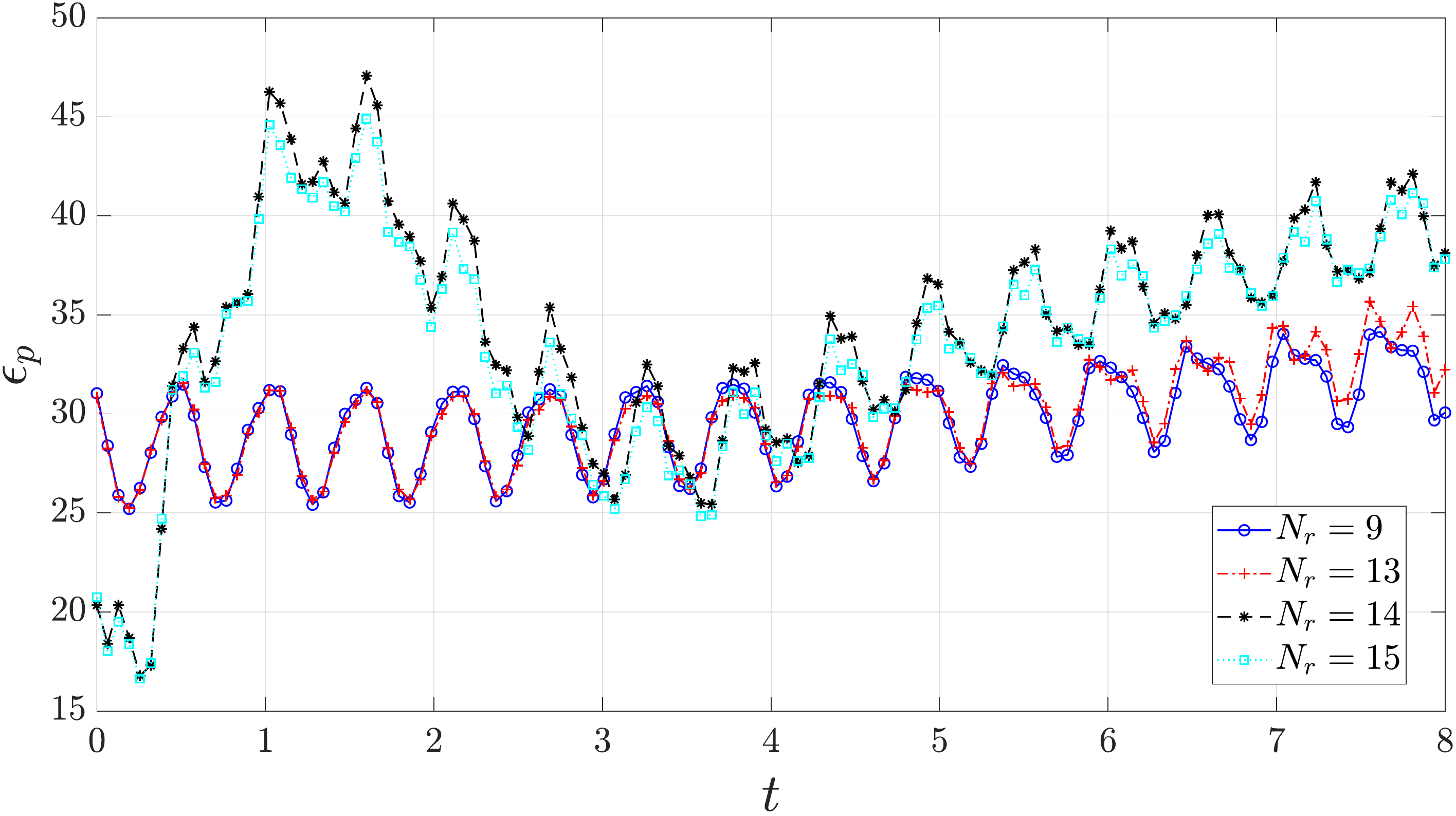} 
    \scriptsize(a) 
    \end{minipage}
  \begin{minipage}[b]{0.5\linewidth}
    \centering
    \includegraphics[width=0.98\linewidth]{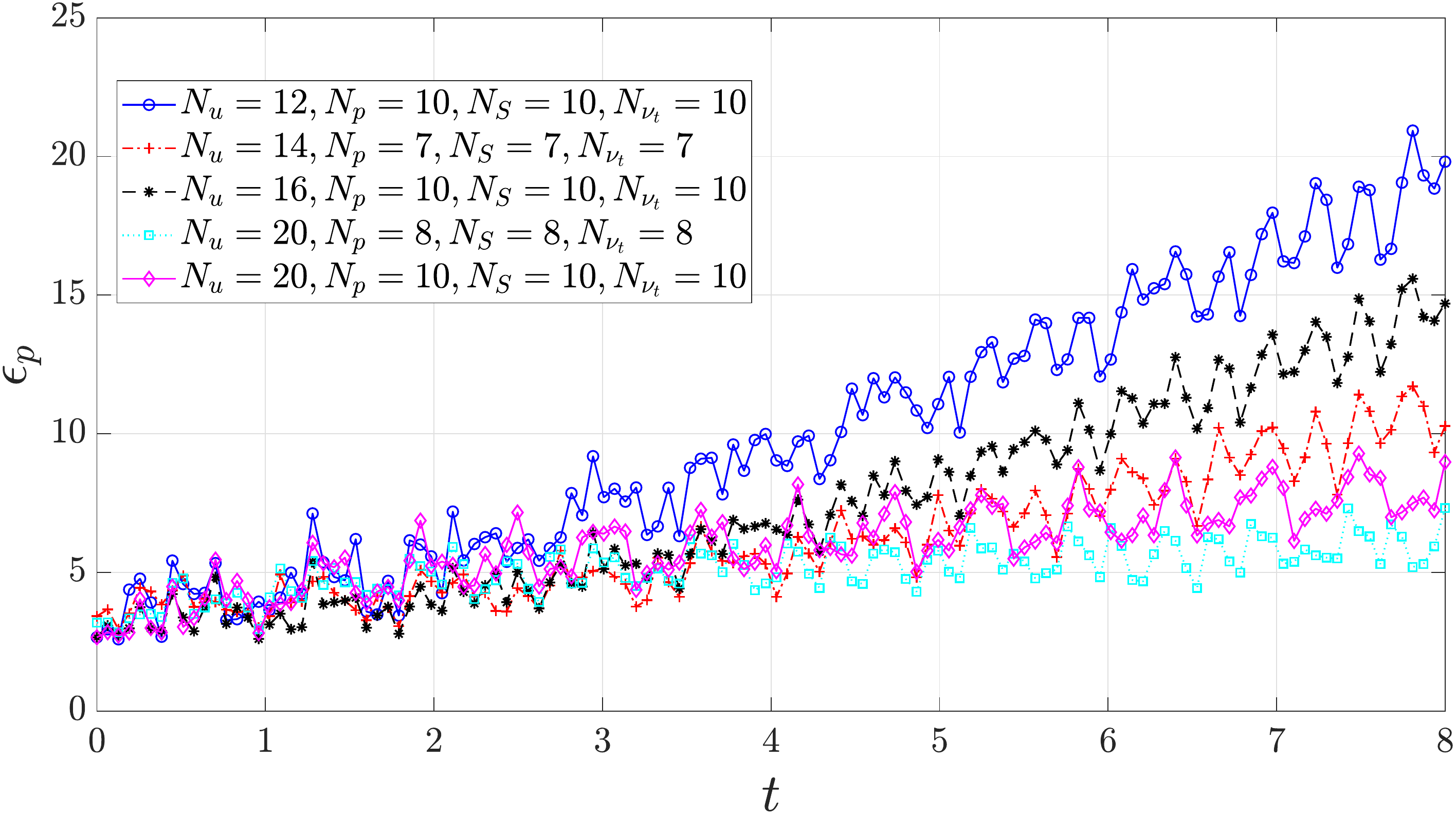}
        \scriptsize(b) 
  \end{minipage}
\caption{The time evolution of the $L^2$ relative errors of the pressure reduced approximations for both the P-ROM and the Mixed-ROM models. The curves correspond to the case run with the parameter value $U_{in} = 7.75$ $\si{m\per s}$ : (a) shows the error curve for the P-ROM model. Figure (b) depicts the case of the Mixed-ROM model. The error values in both graphs are in percentages.}\label{fig:errP_U7.75}
\end{figure}

The behavior of the reduced approximation accuracy over time is analyzed considering the time evolution of the relative $L^2$ error of both the ROM velocity and pressure fields with respect to their FOM counterparts. \autoref{fig:errU_U7.75} depicts values of the velocity error $\epsilon_u$ plotted as a function of time for both the P-ROM and Mixed-ROM models. A similar graph for the pressure field is presented in  \autoref{fig:errP_U7.75}. Both diagrams suggest that the P-ROM model pressure approximation is not as accurate as that obtained with the Mixed-ROM. Again, we must remark that the modal truncation order used in the P-ROM model to generate \autoref{fig:errU_U7.75} and \autoref{fig:errP_U7.75} represent the most accurate choices among all values of $N_r \in [1,30]$, as will be shown in the next results.\par


\begin{figure}
  \centering
 \begin{minipage}[b]{0.5\linewidth}
    \centering
    \includegraphics[width=0.98\linewidth]{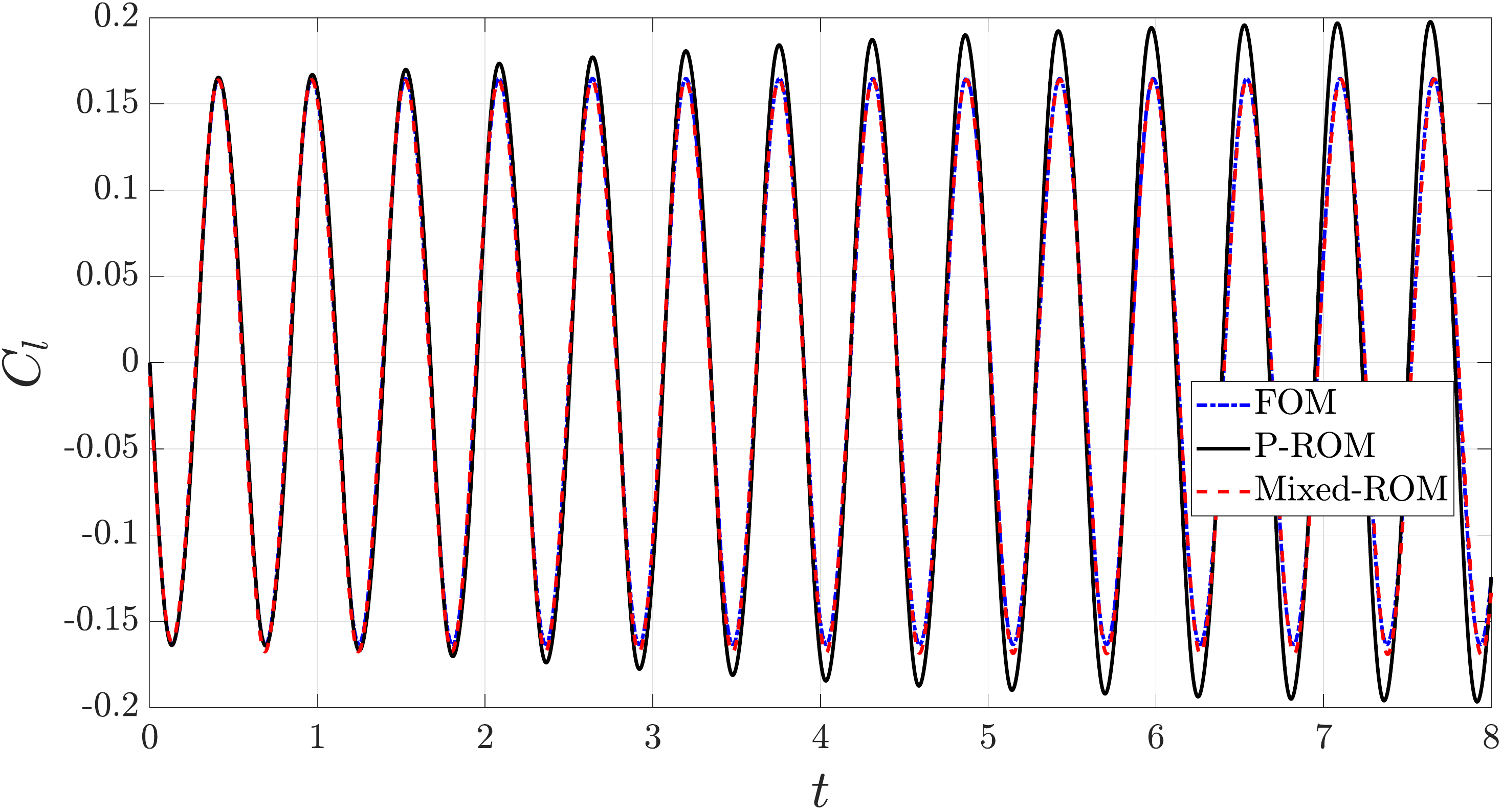} 
    \scriptsize(a) 
    \end{minipage}
  \begin{minipage}[b]{0.5\linewidth}
    \centering
    \includegraphics[width=0.98\linewidth]{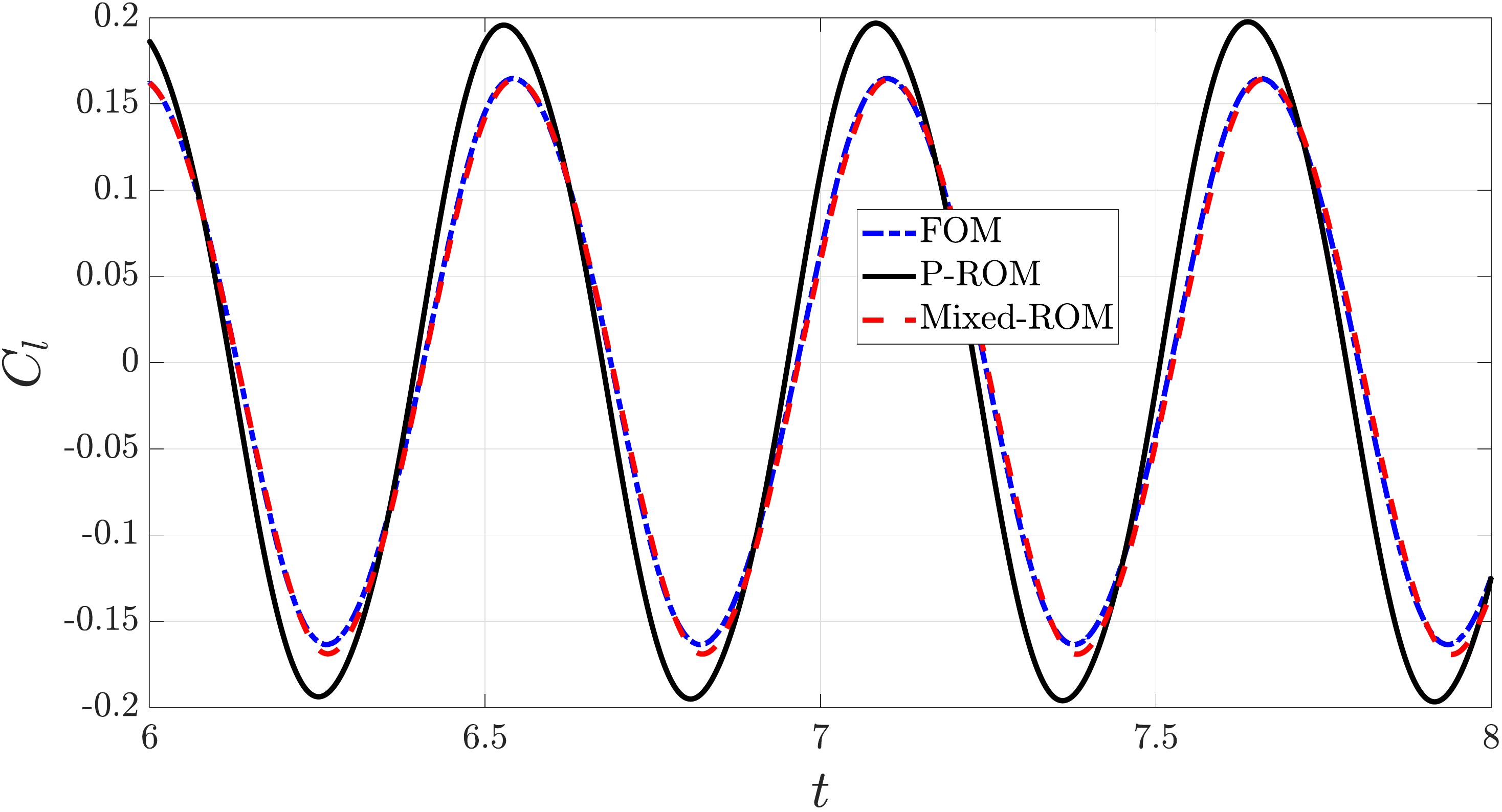}
        \scriptsize(b) 
  \end{minipage}
 \vspace{-0.2cm}\caption{Lift coefficients curves for the cross validation test done for the parameter value $U_{in} = 7.75$ $\si{m\per s}$ for the time range $[0,8]$ $\si{s}$, the figure shows the FOM, the P-ROM and the Mixed-ROM lift coefficients histories : (a) the full range is shown (b) the last $2$ $\si{s}$ $C_l$ is shown.}\label{fig:lift_figures_U7.75} 
\end{figure}

\begin{figure}
  \centering
 \begin{minipage}[b]{0.5\linewidth}
    \centering
    \includegraphics[width=0.98\linewidth]{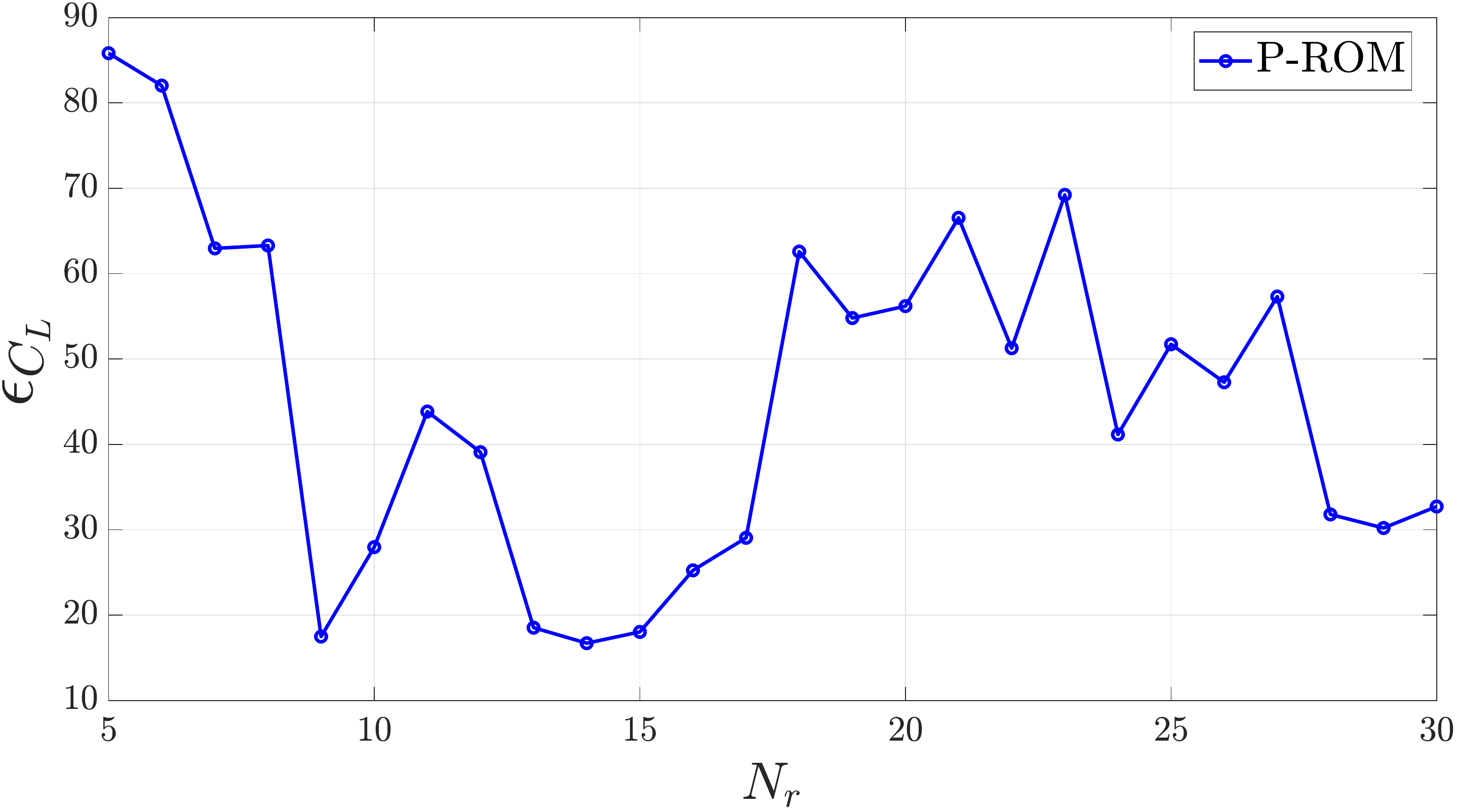} 
    \scriptsize(a) 
    \end{minipage}
  \begin{minipage}[b]{0.5\linewidth}
    \centering
    \includegraphics[width=0.98\linewidth]{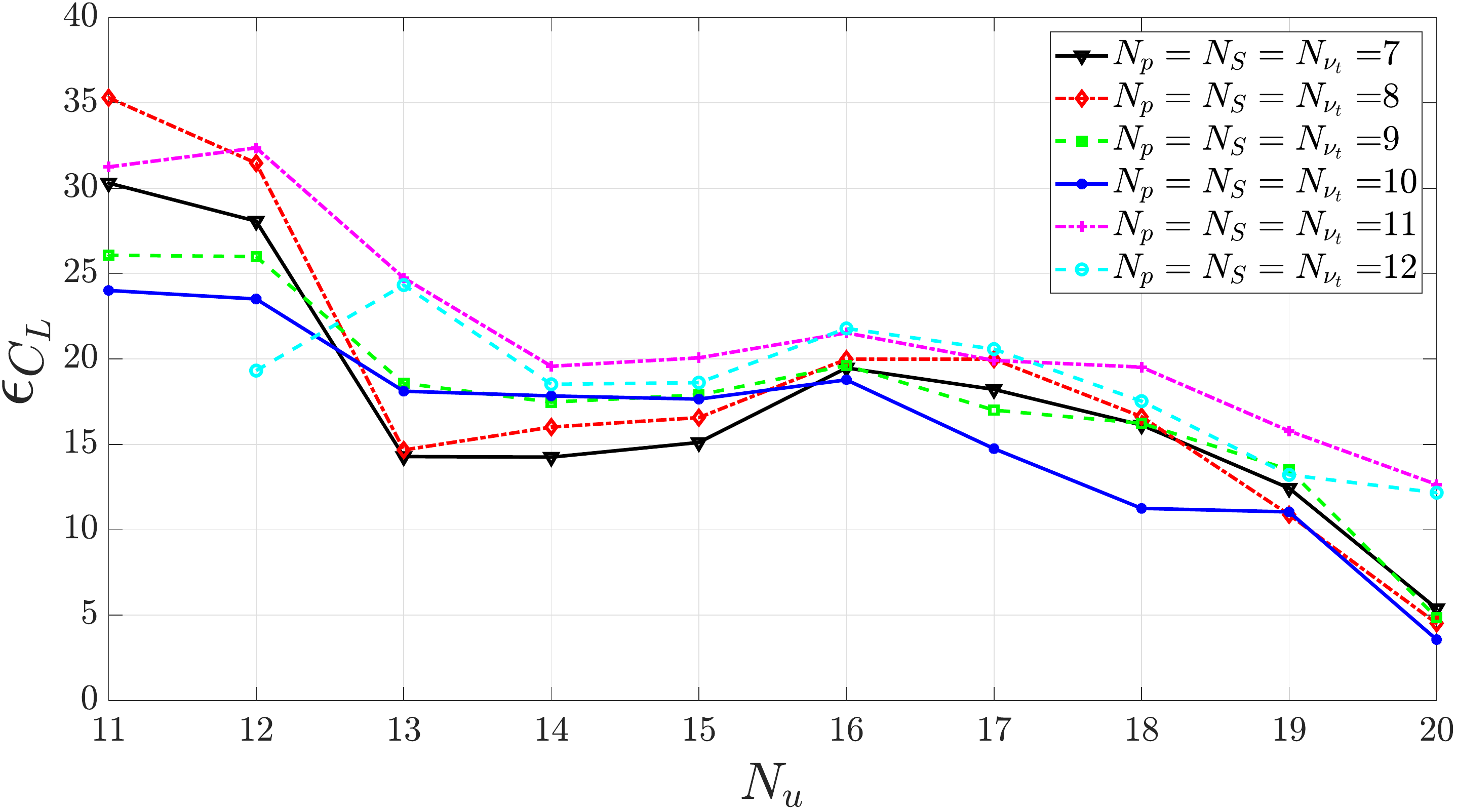}
        \scriptsize(b) 
  \end{minipage}
\caption{The graph of the $L^2$ relative errors for the lift coefficients curve versus number of modes used in the online stage in both cases of the P-ROM and the Mixed-ROM models. The curves correspond to the case run with the parameter value $U_{in} = 7.75$ $\si{m\per s}$. The error is computed between the lift coefficients curve obtained by the FOM solver and the one reconstructed from both the P-ROM and the Mixed-ROM models for the time range $[0,8]$ \si{s} : (a) shows the error curve for the P-ROM model, where $N_r$ is the number of modes used in the online stage for all variables (by construction of the P-ROM it is not possible to choose different number of online modes for the reduced variables). Figure (b) depicts the case of the Mixed-ROM model, where one can see the error values varying the number of modes used for the pure velocity with different fixed settings for the three other variables (the pressure, the supremizers and the eddy viscosity). The error values in both graphs are in percentages.}\label{fig:Lift_cy_err}
\end{figure}

One of the main goals for researchers and engineers studying fluid dynamic problems such as the crossflow cylinder one here considered, is often the evaluation of a force acting on a body or a boundary surface in general. As such forces depend on the local values of the pressure and velocity fields around the body of interest, global error evaluators shown so far might not be good indicators if the aim is that of assessing how well the ROM solvers are able to predict the fluid dynamic forces acting on a body. In the case of the present numerical test for instance, a considerable pressure or velocity error localized in the small region around the cylinder might have a substantial impact on the forces values, while having little effect on the global fields errors. For such reason, the following analysis considers the time evolution of the
lift coefficient $C_l$, \emph{i.e.}: the non-dimensionalized vertical component of the fluid dynamic force acting on the cylinder. It is important to point out that the lift and drag forces exerted by the fluid on the cylinder are not a direct result of the Mixed-ROM computations. The reduced system solution consists in fact in the modal coefficients of the velocity and pressure fields at each time instant, which are in turn used to obtain the Mixed-ROM approximation of the full rank flow field. Such approximation can be obviously used to obtain --- through integration of pressure and skin friction on the cylinder surface --- the reduced order approximation of the fluid dynamic force components and the corresponding force non-dimensional coefficients. Yet, in the reduced order model community this procedure is typically avoided, as it involves a possibly expensive operation such as the evaluation of the full rank flow field. For this reason, the lift and drag coefficients in this work are computed in a fully reduced order fashion, based on the offline computation of suitable matrices which are then used in the online stage. The detailed procedure for online fluid dynamic forces computation is explained in \ref{appendix:B}. To provide an evaluation of the reduced model $C_l$ approximation throughout the whole time integration, \autoref{fig:lift_figures_U7.75} depicts time evolution of the lift coefficient obtained with the FOM, the P-ROM and the Mixed-ROM solvers for the inlet velocity parameter value $U_{in} = 7.75$ $\si{m\per s}$. The left plot (a) shows the $C_l$ values for the full time range under consideration which is $[0,8]$ $\si{s}$. The right diagram (b) represents a detail of  the last three cycles of the time span. The plots clearly indicate that the Mixed-ROM model outperforms the P-ROM model in the $C_l$ approximation. The Mixed-ROM $C_l$ curves seem in fact to closely approximate the FOM lift coefficient ones. The P-ROM $C_l$ approximations are instead not completely accurate and it is evident that the P-ROM suffers from instability issues, as through time integration the P-ROM curve diverges from the FOM one. More quantitative assessment of the lift coefficient accuracy during the time integration is obtained through the evaluation of the $L^2$  relative percentage error, in the integration time interval $[T1,T2]$, between the reduced model approximations of the lift coefficients and their FOM counterparts, namely
\begin{equation}
\epsilon_{C_L} = \frac{{\left\lVert C_l(t) - {C_l}^*(t) \right\rVert}_{L_2(T_1,T_2)}}{{\left\lVert C_l(t) \right\rVert}_{L_2(T_1,T_2)}} \times 100 \%.
\end{equation}
Here, $C_l(t)$ is the time signal of the values of the FOM lift coefficients at all time instants between $T_1$ and $T_2$. On the other hand ${C_l}^*(t)$ is the time evolution of the lift coefficients computed by the reduced order model --- whether P-ROM or Mixed-ROM. \autoref{fig:Lift_cy_err} depicts the $L^2$ relative errors between the reduced model approximation of the lift coefficients and their FOM counterparts, as a function of the online phase modes employed. As expected from the previous figures the convergence plots highlight that the Mixed-ROM model is able to reproduce the FOM force coefficient with significantly greater accuracy than the ROM. In fact, the Mixed-ROM error reaches values as low as $3$ $\%$, while the ROM $C_l$ are consistently above $16$ $\%$ off the FOM values.\par

\begin{figure}
  \centering
 \begin{minipage}[b]{0.5\linewidth}
    \centering
    \includegraphics[width=0.98\linewidth]{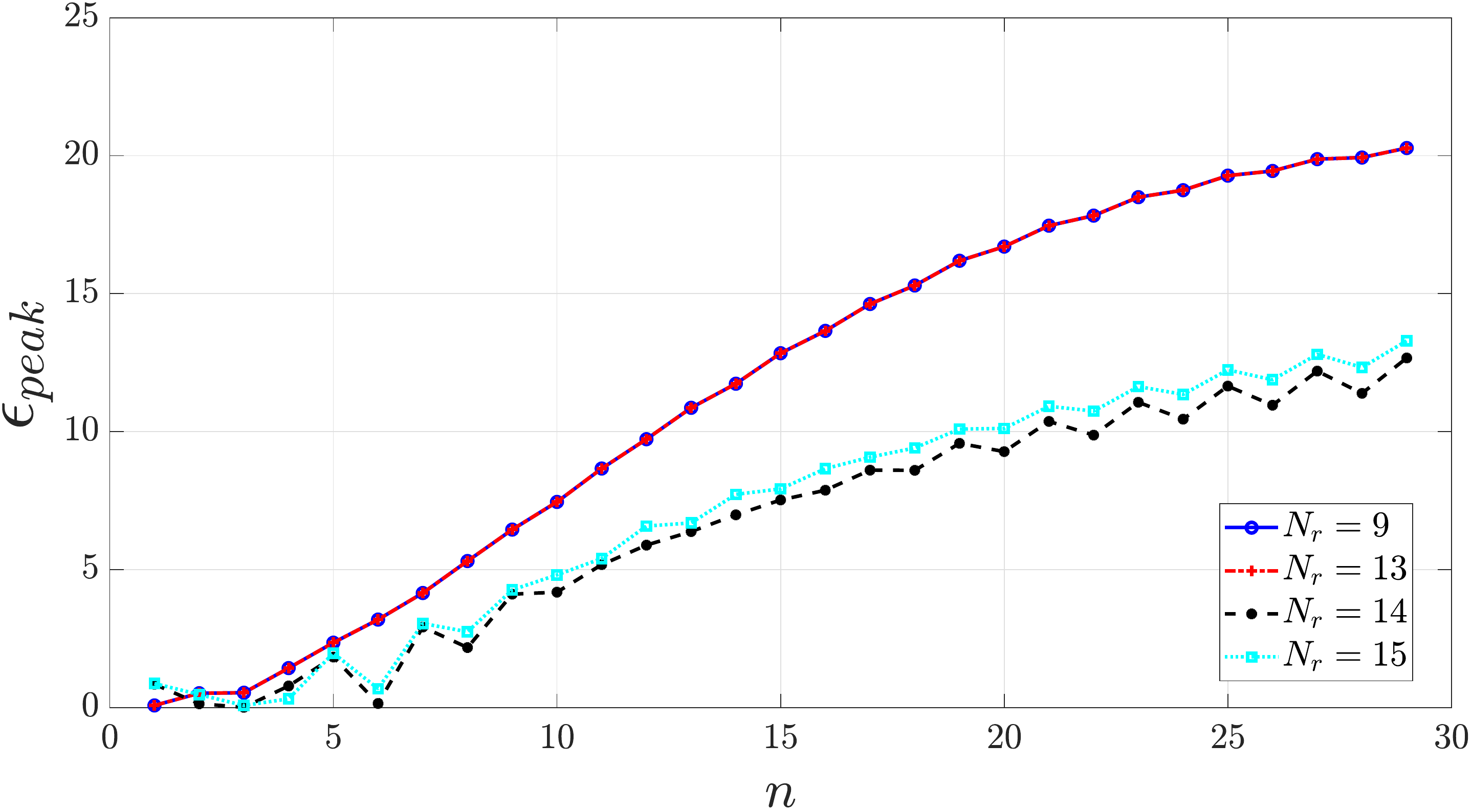} 
    \scriptsize(a) 
    \end{minipage}
  \begin{minipage}[b]{0.5\linewidth}
    \centering
    \includegraphics[width=0.98\linewidth]{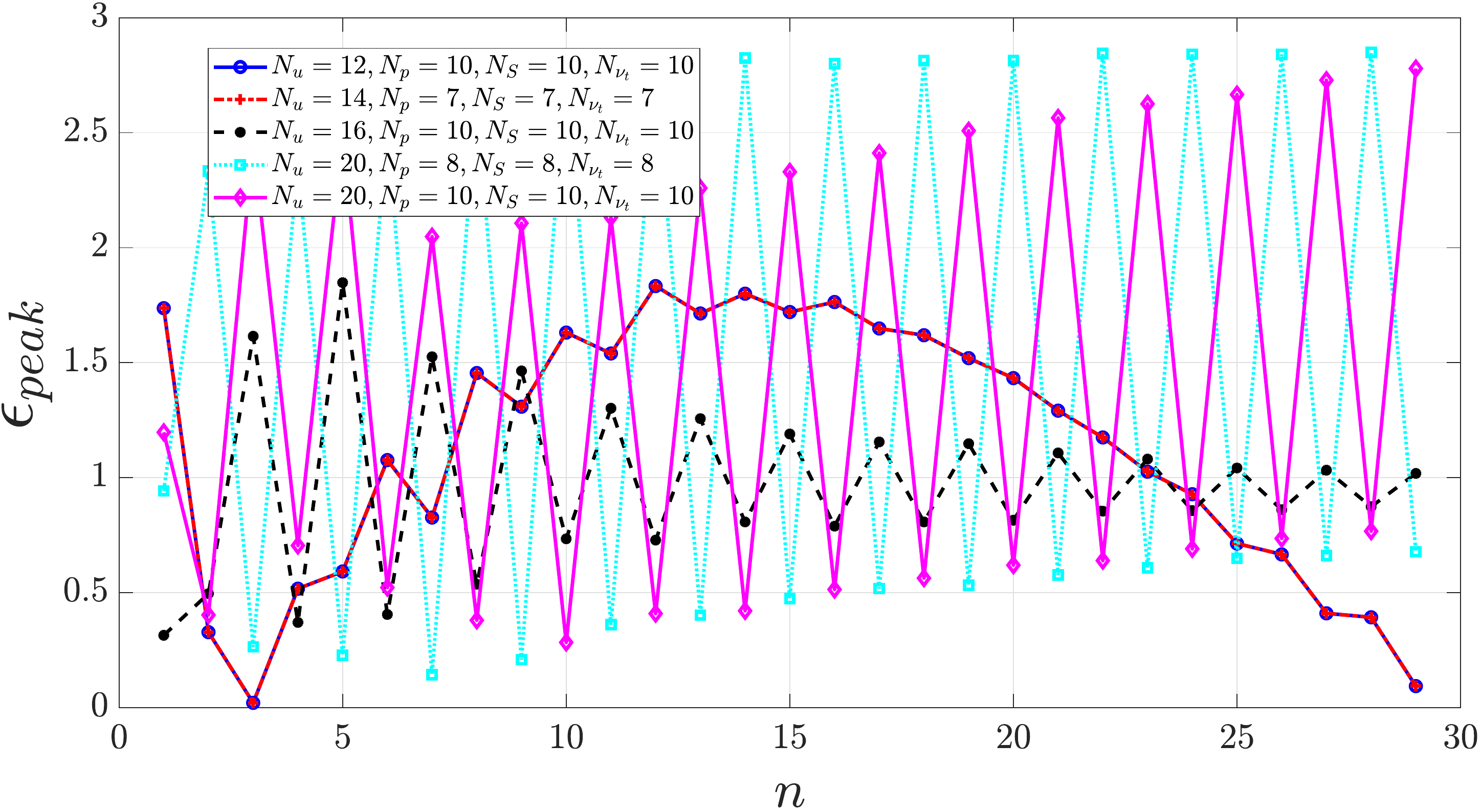}
        \scriptsize(b) 
  \end{minipage}
\caption{The graph of the peaks relative errors for the lift coefficients curves for varied values of the number of modes used in the online stage in both cases of the P-ROM and the Mixed-ROM models. The curves correspond to the case run with the parameter value $U_{in} = 7.75$ $\si{m\per s}$. The error is computed between the peaks values of the lift coefficients curve obtained by the FOM solver and the ones reconstructed from both the P-ROM and the Mixed-ROM models for the time range $[0,8]$ \si{s} : (a) shows the error curve for the P-ROM model, where $N_r$ is the number of modes used in the online stage for all variables (by construction of the P-ROM it is not possible to choose different number of online modes for the reduced variables). Figure (b) depicts the case of the Mixed-ROM model. The error values in both graphs are in percentages.}\label{fig:peak_err}
\end{figure}

To further analyze the results in \autoref{fig:Lift_cy_err}, we also attempt to assess how much the curve $L^2$ error is due to incorrect reproduction of the amplitude or frequency of the lift coefficient oscillations. To evaluate, from quantitative perspective, the accuracy of the lift coefficients peak prediction, we define the relative peak error $\epsilon_{peak}$ as follows:
\begin{equation}
\epsilon_{n,peak} = \frac{PK_{n,FOM}-PK_{n,*}}{PK_{n,FOM}} \times 100 \%,
\end{equation}
where $PK_{n,FOM}$ is the value of the $n-$th FOM $C_l$ peak and $PK_{n,*}$ is the value of the $n-$th P-ROM or Mixed-ROM $C_l$ peak. The relative peak error is plotted for both P-ROM and Mixed-ROM models in \autoref{fig:peak_err}. \autoref{fig:peak_err} presents the relative peak error values obtained for each of the 29 peaks the time interval $[0,8]$ \si{s}. We point out that the values of modal truncation order $N_r$ in correspondence of which the peak errors are computed are the resulted most accurate in the relative $L^2$ lift error analysis presented. \autoref{fig:peak_err} suggests that the P-ROM peaks relative error grows in time and settles around values as high as  $10-20$ $\%$. On the other hand, the corresponding error values for the Mixed-ROM model are less than $3.5$ $\%$ for several different modal truncation order for velocity, pressure, supremizers and eddy viscosity.

A final numerical test is aimed at assessing the accuracy of the Mixed-ROM model for higher Reynolds value. The inlet velocity parameter sample considered in this case is $U_{in} = 11.75$ $\si{m\per s}$. The time interval considered for the reduced order simulations is $t \in [0,10]$ $\si{s}$ which contains around $27$ solution cycles. The results reported for this case are relative to the lift coefficient history, the $L^2$ relative error value for $C_l$, the $C_l$ peaks error and the approximated time period by the Mixed-ROM.\par
\begin{figure}
  \centering
 \begin{minipage}[b]{0.5\linewidth}
    \centering
    \includegraphics[width=0.98\linewidth]{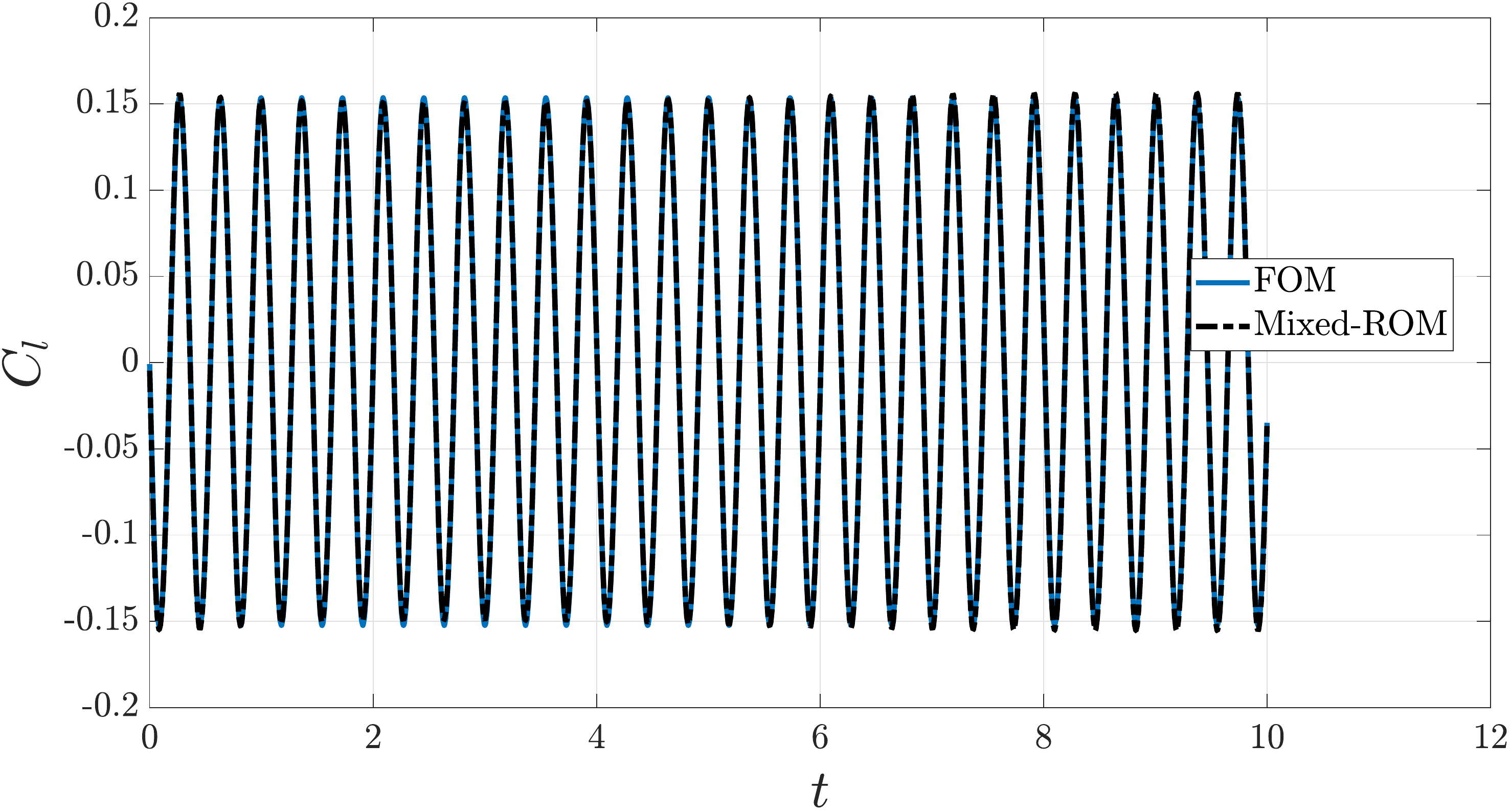} 
    \scriptsize(a) 
    \end{minipage}
  \begin{minipage}[b]{0.5\linewidth}
    \centering
    \includegraphics[width=0.98\linewidth]{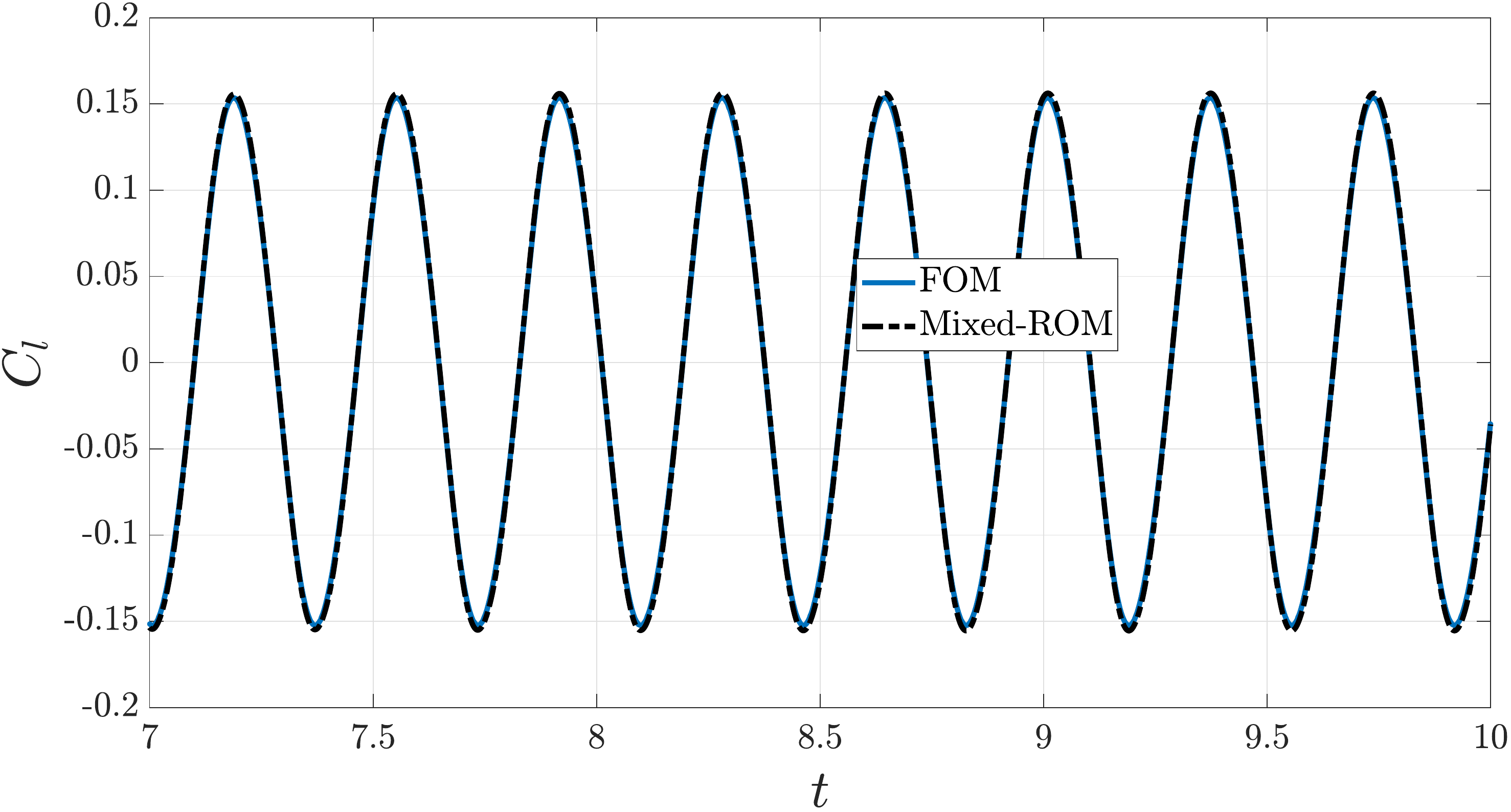}
        \scriptsize(b) 
  \end{minipage}
 \vspace{-0.2cm}\caption{Lift coefficients curves for the cross validation test done for the parameter value $U_{in} = 11.75$ $\si{m\per s}$ for the time range $[0,10]$ $\si{s}$, the figure shows the FOM and the Mixed-ROM lift coefficients histories : (a) the full range is shown (b) the last $3$ $\si{s}$ history of $C_l$ is shown.}\label{fig:lift_figures_U11.75} 
\end{figure}

The Mixed-ROM dynamical system is solved with time step equal to $0.00025$ $\si{s}$. The Mixed-ROM fields were reproduced using $12$ modes for velocity and $10$ for each of pressure, supremizers and eddy viscosity. The $C_l$ curves obtained with the FOM solver and the Mixed-ROM model are presented in \autoref{fig:lift_figures_U11.75}.  The results in \autoref{fig:lift_figures_U11.75} prove that the Mixed-ROM was successful in reducing the problem with satisfactory accuracy. The $L^2$ relative error between FOM and Mixed-ROM solution is in fact $1.9654\%$. As for the relative peak error, the highest value detected in the $[0,10]$ \si{s} time interval is $2.0672\%$. Finally the average time period computed in by the FOM solver is about $0.3641$ $\si{s}$, while the average time period computed by the Mixed-ROM is roughly $0.3642$ $\si{s}$. The corresponding 0.2\% relative error suggests that the main source of error in the ROM predictions is due to the amplitude inaccuracies rather then to incorrect frequency reconstructions.\par

\begin{figure}
  \centering
  \begin{minipage}[b]{0.5\linewidth}
    \centering
    \includegraphics[width=0.98\linewidth]{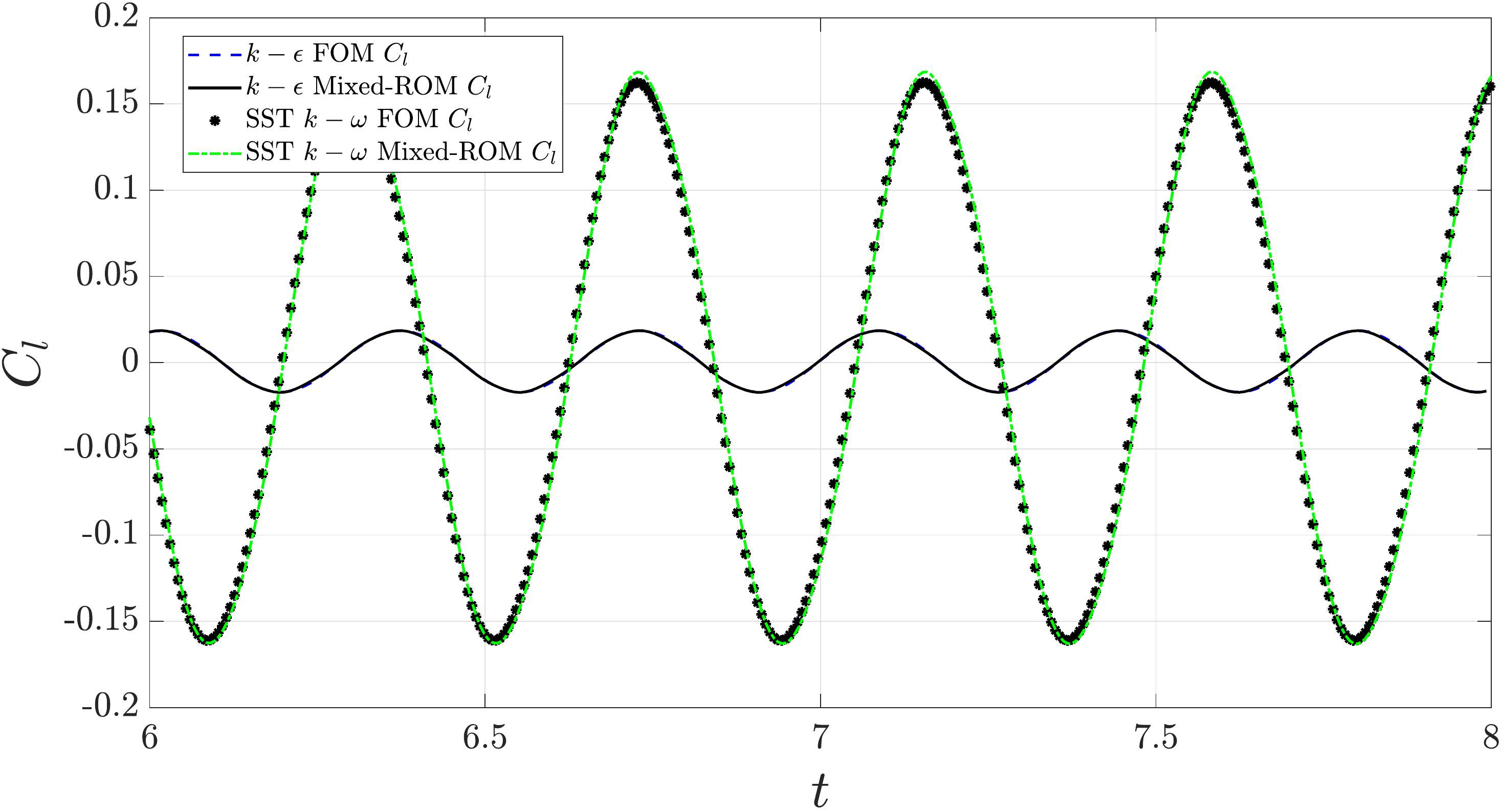}
  \end{minipage}
  \vspace{-0.2cm}\caption{The lift coefficient curves obtained using both $k-\epsilon$ and SST $k-\omega$ turbulence models and the Mixed-ROM ones. The case considered is a non-parametrized one with $U_{in} = 10$ $\si{m\per s}$ corresponding to $Re=10^5$. The plot is for the time range $t \in [6,8]$, the Mixed-ROM achieved relative $L^2$ errors (over the range $t \in [0,8]$) which are less than $5$ $\%$ in both cases.}\label{fig:LiftCylinderDiffTurbModels}\par
\end{figure}

The final test in this section is meant to assess one of the main objectives of this work, that is to test the presented reduction approach for the variance of the FOM turbulence model. We considered the non-parametrized case of $Re= 10^5$, the FOM was run for both $k-\epsilon$ and SST $k-\omega$ models. After having reached the fully periodic regime, snapshots were taken for the first $1.2$ $\si{s}$ and $1.6$ $\si{s}$ for $k-\epsilon$ and SST $k-\omega$ models, respectively. The reduction was done extrapolating in time, where the Mixed-ROM simulations were run for $8$ $\si{s}$. The lift coefficient curves are shown in \autoref{fig:LiftCylinderDiffTurbModels}, where one can see both FOM $C_l$ signals for the two different turbulence models and their reduced counterparts. It is evident from the graph that the Mixed-ROM proves sensitive to the specific turbulence model used in the FOM solver, although no additional PDEs for the turbulent quantities are solved at the reduced level.

\section{Conclusions and Outlook} 
This work presents a hybrid data-driven/projection-based approach to reduce turbulent flows. The approach developed in this work called Mixed-ROM is based on introducing a non-intrusive reduced order version of the eddy viscosity field to the formulation of the reduced order model. The Mixed-ROM employs interpolation using radial basis function in the online stage for the computation of the reduced order eddy viscosity coefficients. This interpolation can be done with the independent variable being the combined time-parameter vector or the combined vector of the velocity $L^2$ projection coefficients and their vector time derivatives. The Mixed-ROM proved to be accurate in reconstructing the fluid dynamics fields in both cases of steady and unsteady flows with a Reynolds number on the order of $10^5$. In the unsteady case considered in this work which is the flow around a circular cylinder, the Mixed-ROM showed that it is capable of reconstructing an important variable of interest that is the lift coefficient time history which mainly depends on local flow features around the cylinder. In that same example, the Mixed-ROM gives satisfactory results when it comes to the extrapolation in time. The Mixed-ROM has been able to obtain accurate predictions with an acceptable computational cost, showing with a speed up of $SU = 10$ in the unsteady case and around $SU = 1000$ in the steady case.\par


As for potential future work, data-driven techniques can be used in building reduced order models and other methodologies could help in approximating certain maps which are needed for the ultimate goal of reducing CFD problems. Such methodologies include the Artificial Neural Networks (ANN) with which one could potentially improve the accuracy of the approximation of the eddy viscosity coefficients conducted in this work. Another idea is to use DMD for the extrapolation problem for the unsteady flows. In addition, there is a need to find stabilization techniques for the long time integration problem for unsteady flows \cite{Fick2018} and for multi-physics problems \cite{GeoStaRoBlu2019, GeoStaRoBlu2018, 8627766, Busto2020}.

\section*{Acknowledgements} 
We acknowledge the support provided by the European Research Council Executive Agency by the Consolidator Grant project AROMA-CFD ``Advanced Reduced Order Methods with Applications in Computational Fluid Dynamics'' - GA 681447, H2020-ERC CoG 2015 AROMA-CFD and INdAM-GNCS projects. In addition, this work was partially performed in the framework of the project SOPHYA - ``Seakeeping Of Planing Hull YAchts'' and of the project PRELICA - ``Advanced methodologies for hydro-acoustic design of naval propulsion'', both supported by Regione FVG, POR-FESR 2014-2020, Piano Operativo Regionale Fondo Europeo per lo Sviluppo Regionale.
\section*{Appendix A. List of abbreviations and symbols}
\printnomenclature
\Urlmuskip=0mu plus 1mu\relax
\section*{Appendix B. Lift and drag forces offline/online computations}
\label{appendix:B}
This section introduces the computations done in both the offline and the online stages for obtaining the surface forces acting on a part of the domain called $\partial \Omega_f$.\par 
The total viscous and pressure forces $\bm{F}$ acting on $\partial \Omega_f$ are given by the following integral:
\begin{equation}\label{eq:forces}
\bm{F} = \int_{\partial \Omega_f} (2\mu \bm{\nabla} \bm{u} - p \bm{I}) \bm{n}ds.
\end{equation}
In many application in fluid dynamics it is very important to efficiently compute the forces acting on certain objects inside the domain. For instance the problem of flow past a circular cylinder considered in this work is one of them. One should avoid resorting to the full order mesh for computing the integral above because this makes the approach not entirely a reduced one.\par 
The first step in developing an offline/online decoupling approach for computing the forces is to insert the approximation \eqref{eq:decompose} into \eqref{eq:forces}, this yields the following:
\begin{equation}
\bm{F} = \int_{\partial \Omega_f} (2\mu \bm{\nabla} (\sum_{i=1}^{N_u} a_i(t;\bm{\mu}) \bm{\phi}_i(\bm{x})) - \sum_{i=1}^{N_p} b_i (t;\bm{\mu}) {\chi_i} \bm{I})  \bm{n}ds,
\end{equation}
\begin{equation}
\bm{F} = \int_{\partial \Omega_f} 2\mu \sum_{i=1}^{N_u} a_i(t;\bm{\mu}) \bm{\nabla} \bm{\phi}_i(\bm{x})  \bm{n}ds - \int_{\partial \Omega_f} \sum_{i=1}^{N_p} b_i (t;\bm{\mu}) {\chi_i} \bm{n}ds,
\end{equation}
\begin{equation}
\bm{F} = \sum_{i=1}^{N_u} a_i(t;\bm{\mu}) \int_{\partial \Omega_f} 2\mu \bm{\nabla} \bm{\phi}_i(\bm{x})  \bm{n}ds -  \sum_{i=1}^{N_p} b_i (t;\bm{\mu})  \int_{\partial \Omega_f}{\chi_i} \bm{n}ds.
\end{equation}
After reaching this point one can define the following quantities

\begin{align}\label{eq:term1}
& \bm{\delta}_i=\int_{\partial \Omega_f} 2\mu \bm{\nabla} \bm{\phi}_i(\bm{x})  \bm{n}ds, \quad \text{for} \quad i=1,...,N_u, \\
& \bm{\theta}_j=\int_{\partial \Omega_f}{\chi_j} \bm{n}ds,  \quad \text{for} \quad j=1,...,N_p,\label{eq:term2} 
\end{align}
where each term of $\bm{\nabla} \bm{\phi}_i(\bm{x})$ and ${\chi_j}$ can be seen as velocity and pressure field, respectively. This will make the computations of \eqref{eq:term1} and \eqref{eq:term2} possible in the offline stage and they will be stored in order to be later used in the online stage.\par
In the online stage when a new time-parameter vector $\bm{z}^*$ is introduced, the forces are computed as follows:
\begin{equation}
\bm{F}^* = \int_{\partial \Omega_f} (2\mu \bm{\nabla} \bm{u}(\bm{z}^*,\bm{x}) - p(\bm{z}^*,\bm{x}) \bm{I})  \bm{n}ds,
\end{equation}
which simplify to
\begin{equation}
\bm{F}^* = \sum_{i=1}^{N_u} a_i(\bm{z}^*) \bm{\delta}_i -  \sum_{j=1}^{N_p} b_j (\bm{z}^*)  \bm{\theta}_j.
\end{equation}
\bibliographystyle{amsplain_mod}
\bibliography{bib/bibfile} 
\end{document}